\documentclass[a4paper,11pt,oneside,reqno]{amsart}
\usepackage{amsfonts,amsmath, amssymb,amsthm,amscd}
\usepackage[usenames,dvipsnames]{color}
\usepackage{bm}
\usepackage{mathrsfs}
\usepackage{yfonts}
\usepackage[a4paper,mpexclude,DIV13]{typearea}
\usepackage{verbatim}
\usepackage{graphicx}
\usepackage[latin1]{inputenc}
\usepackage{latexsym}
\usepackage{lscape}
\usepackage{epsfig}
\usepackage{subfigure}
\include{bibtex}
\usepackage{setspace}

\begin{document}
\bibliographystyle{alpha}

\newcommand{\var}[1]{\mathrm{\overline{var}}(#1)}

\newcommand{\cn}[1]{\overline{#1}}
\newcommand{\e}[0]{\epsilon}
\newcommand{\eps}[0]{\epsilon}
\newcommand{\PP}{\ensuremath{\mathbb{P}}}

\newcommand{\hfluc}{\ensuremath{h_{\epsilon}^{\textrm{fluc}}}}

\newcommand{\re}{\ensuremath{\mathrm{Re}}}
\newcommand{\im}{\ensuremath{\mathrm{Im}}}
\newcommand{\N}{\ensuremath{\mathbb{N}}}
\newcommand{\R}{\ensuremath{\mathbb{R}}}
\newcommand{\C}{\ensuremath{\mathbb{C}}}
\newcommand{\Z}{\ensuremath{\mathbb{Z}}}
\newcommand{\Q}{\ensuremath{\mathbb{Q}}}
\newcommand{\T}{\ensuremath{\mathbb{T}}}
\newcommand{\E}{\ensuremath{\mathbb{E}}}
\def \Ai {{\rm Ai}}

\newcommand{\OO}[0]{\Omega}
\newcommand{\F}[0]{\mathfrak{F}}
\newcommand{\G}[0]{\mathfrak{G}}
\newcommand{\ta}[0]{\theta}
\newcommand{\w}[0]{\omega}
\newcommand{\ra}[0]{\rightarrow}
\newtheorem{theorem}{Theorem}[section]
\newtheorem{partialtheorem}{Partial Theorem}[section]
\newtheorem{conj}[theorem]{Conjecture}
\newtheorem{lemma}[theorem]{Lemma}
\newtheorem{proposition}[theorem]{Proposition}
\newtheorem{corollary}[theorem]{Corollary}
\newtheorem{claim}[theorem]{Claim}
\newtheorem{experiment}[theorem]{Experimental Result}

\def\todo#1{\marginpar{\raggedright\footnotesize #1}}
\def\change#1{{\color{green}\todo{change}#1}}
\def\note#1{\textup{\textsf{\color{blue}(#1)}}}

\theoremstyle{definition}
\newtheorem{rem}[theorem]{Remark}

\theoremstyle{definition}
\newtheorem{com}[theorem]{Comment}

\theoremstyle{definition}
\newtheorem{definition}[theorem]{Definition}

\theoremstyle{definition}
\newtheorem{definitions}[theorem]{Definitions}

\title{The Kardar-Parisi-Zhang equation and universality class}
\author[I. Corwin]{Ivan Corwin}
\address{I. Corwin\\
  Courant Institute\\
  New York University\\
  251 Mercer Street\\
  New York, NY 10012, USA}
\email{corwin@cims.nyu.edu}

\maketitle

\begin{abstract}
Brownian motion is a continuum scaling limit for a wide class of random processes, and there has been great success in developing a theory for its properties (such as distribution functions or regularity) and expanding the breadth of its universality class. Over the past twenty five years a new universality class has emerged to describe a host of important physical and probabilistic models (including one dimensional interface growth processes, interacting particle systems and polymers in random environments) which display characteristic, though unusual, scalings and new statistics. This class is called the Kardar-Parisi-Zhang (KPZ) universality class and underlying it is, again, a continuum object -- a non-linear stochastic partial differential equation -- known as the KPZ equation.

The purpose of this survey is to explain the context for, as well as the content of a number of mathematical breakthroughs which have culminated in the derivation of the exact formula for the distribution function of the KPZ equation started with {\it narrow wedge} initial data.  In particular we emphasize three topics: (1) The approximation of the KPZ equation through the weakly asymmetric simple exclusion process; (2) The derivation of the exact one-point distribution of the solution to the KPZ equation with narrow wedge initial data; (3) Connections with directed polymers in random media.

As the purpose of this article is to survey and review, we make precise statements but provide only heuristic arguments with indications of the technical complexities necessary to make such arguments mathematically rigorous.
\end{abstract}

\tableofcontents

\section{Introduction}

The overall goal of this review is to give a clear overview of over twenty-five years of work -- starting even before the seminal paper of Kardar-Parisi-Zhang \cite{KPZ} in 1986 -- which has culminated in early 2010 with the discovery of the probability distribution for the solution to the KPZ stochastic PDE (not to be confused with the KPZ of quantum gravity) and the understanding that the KPZ equation is the universal continuum mechanism for the crossover between the two universality classes associated with growth models: The KPZ class and the EW (Edwards Wilkinson) class \cite{EW}. The clarity which I aim for here is very much a benefit of hindsight, as is the program outlined for how everything fits together. Emphasis will be on the important conceptual steps and calculations involved in this program -- not on the many technical challenges and difficulties which arise along the way. While many of the important developments in the study of the KPZ equation will be explained or alluded to, this literature is vast, and will not be surveyed entirely herein.

The introduction to this article is meant to quickly provide the reader with the main results as well as philosophies associated with the KPZ equation and its universality class. While the main focus of this review is on growth processes and their relationship to the KPZ equation, in the last section the connections to the study of directed polymers in random media will be illuminated.

This article is focused primarily on the mathematically rigorous side of results pertaining to the KPZ equation and universality class. However, care has been taken to accurately represent the accomplishments of physicists in this area and many of the relevant references are cited (along with brief discussions). Still, it is nearly impossible to be exhaustive in recording the accomplishments of the countless people who have worked in this area and thus the author begs the pardon of anyone whose contribution to the field is overlooked in this article.

\subsection{Beyond the Gaussian universality class}
Since its discovery two hundred years ago the Gaussian distribution has come to represent one of mathematics greatest societal contributions -- a robust theory explaining and analyzing much of the randomness inherent in the physical world. Physical and mathematical systems accurately described in terms of Gaussian statistics are said to be in the {\it Gaussian universality class}. This class, however, is not all encompassing. For example, classical extreme value statistics or Poisson statistics better capture the randomness and severity of events ranging from natural disasters to emergency room visits.

Recently, significant research efforts have been focused on understanding systems which are not well described in terms of {\it any} of the classically developed statistical universality classes. The failure of these systems to conform with classical descriptions is generally due to the non-linear relationship between natural observables and underlying sources of random inputs and noise. A variety of models for such complex systems have been actively studied for over forty years in both mathematics and physics.

History should give great credit to those giants of probability, PDEs, representation theory and mathematical/statistical physics who understood that it was worthwhile to dig deep in the study of these systems. They defined these models, identified links between them and with other fields, pinpointed many mysteries and seriously attacked the most important problems here with a full battery of methods. In light of this work and a few other recent breakthroughs I seek to identify and frame the two problems which are of central importance in the study of these systems: (1) Universality of the scalings, statistics and limit objects, and (2) solvability and integrability which enables exact formulas and computations for quantities associated with these systems. These two complementary goals are behind significant amounts of on-going research within probability, mathematical physics and related fields.

This review will focus on three types of systems:
\begin{enumerate}
\item Random growth interfaces and interacting particle systems,
\item Directed polymers in random media,
\item Non-linear SPDEs such as the KPZ equation, stochastic heat equation (SHE) with multiplicative noise and stochastic Burgers equation.
\end{enumerate}
A number of books and review articles have been written about these models in mathematics and physics, including \cite{BaSt, CSY, HHZ, Holl, Kes, KL, KK, KS, Lig, Mea, SM, Sp, Spo}. Though not directly addressed here, the study of these systems is closely related to and influenced by problems in random matrix theory, non-intersecting path ensembles, random tilings and certain combinatorial problems involving asymptotic representation theory \cite{BorodinDet, FSreview, KJSurvey, Sp}. In particular, many (but not all) of the statistics which arise in these systems were first analytically discovered in the context of random matrix theory.

In all of these types of models there are certain observables of interest (such as height functions, particle currents, free energy, or eigenvalues) which have random fluctuations that have been shown, through experimental evidence, numerical simulations and in some cases mathematical proof, to grow as a characteristic power law in the systems size or time. In the case of the classical central limit theorem this power law would have an exponent of $1/2$ and the limiting distribution for the scaled fluctuations would be Gaussian. This, however, is not the case here. Surprisingly for all of the above listed models the observables of interest fluctuate with a scaling exponent of $1/3$ and display the same set of non-Gaussian limiting distributions. That is to say that they form a new universality class -- the {\it KPZ universality class}.

\subsubsection{The Gaussian versus KPZ universality classes}
Let us now briefly illustrate the difference between the Gaussian and KPZ universality classes. All of the models considered here are $1+1$ dimensional, meaning that there is one space and one time dimension. In the context of growth models this means that the growing interface is given by a curve (as opposed to a surface -- a problem of great interest but few predictions and even fewer results).

The {\it random deposition} growth model involves square blocks being stacked in columns (indexed by $\Z$). Growth is driven by independent Poisson processes (one for each column) of blocks falling from above and accumulating on top of the growing stacks of blocks. Due to the independence, each column evolves independently (i.e., with no spatial correlation) while the height at time $t$ fluctuates like $t^{1/2}$ and is governed by the Gaussian distribution (due to the classical central limit theorem).

The {\it ballistic deposition} growth model (shown in Figure 0b) involves the same independent process of blocks falling in each column, however a block sticks to the first block it touches. So if a block falls in column $x$ but the height of column $x+1$ is larger than that of column $x$, the block will stick to the side of column $x+1$ and the height of column $x$ will suddenly jump to equal that of $x+1$. This important modification breaks the independence of the column heights and importantly introduces spatial correlation. The salient features of this model (and all growth models in the KPZ universality class) are threefold:
\begin{enumerate}
\item {\bf Smoothing:} Deep holes are rapidly filled to smooth the interface.
\item {\bf Rotationally invariant, slope dependent, growth speed:} The larger the absolute value of the slope, the larger the height tends to increase.
\item {\bf Space-time uncorrelated noise:} The blocks fall spatially and temporally independently.
\end{enumerate}
Based on these three features it is predicted that in time $t$ the column height will fluctuate around its mean like $t^{1/3}$ and that the correlation of these fluctuations will be non-trivial on a spatial scale of $t^{2/3}$ (in sharp contrast to random deposition). This prediction has not been proved. However, in the case of the {\it corner growth model} (see Figure 0a) the same predictions have been proved \cite{KJ,PS2}. We too will focus on the corner growth model -- not ballistic deposition.

\begin{figure}
\setlength{\unitlength}{1.2pt}
\begin{picture}(200,100)(0,-10)\label{corner}
\put(0,0){\line(1,1){80}}
\put(0,0){\line(-1,1){80}}
\put(-10,10){\line(1,1){10}}
\put(10,10){\line(-1,1){10}}
\put(-20,20){\line(1,1){10}}
\put(0,20){\line(-1,1){10}}
\put(-30,30){\line(1,1){10}}
\put(-10,30){\line(-1,1){10}}
\put(0,20){\line(1,1){10}}
\put(20,20){\line(-1,1){10}}

\put(10,30){\line(1,1){10}}
\put(30,30){\line(-1,1){10}}

\put(20,40){\line(1,1){10}}
\put(40,40){\line(-1,1){10}}

\put(-80,0){\line(1,0){160}}
\put(0,-5){\makebox(0,0){$\textrm{\tiny{0}}$}}
\put(-10,-5){\makebox(0,0){$\textrm{\tiny{-1}}$}}
\put(-20,-5){\makebox(0,0){$\textrm{\tiny{-2}}$}}
\put(-30,-5){\makebox(0,0){$\textrm{\tiny{...}}$}}
\put(10,-5){\makebox(0,0){$\textrm{\tiny{1}}$}}
\put(20,-5){\makebox(0,0){$\textrm{\tiny{2}}$}}
\put(30,-5){\makebox(0,0){$\textrm{\tiny{...}}$}}

\put(-80,-40){\line(1,1){10}}
\put(-80,-40){\line(-1,1){10}}

\put(-50,-35){\makebox(0,0){$\stackrel{\textrm{\tiny{rate q}}}{\longrightarrow}$}}

\put(-20,-40){\line(1,1){10}}
\put(-20,-40){\line(-1,1){10}}
\put(-30,-30){\line(1,1){10}}
\put(-10,-30){\line(-1,1){10}}


\put(30,-40){\line(1,1){10}}
\put(30,-40){\line(-1,1){10}}
\put(20,-30){\line(1,1){10}}
\put(40,-30){\line(-1,1){10}}

\put(60,-35){\makebox(0,0){$\stackrel{\textrm{\tiny{rate p}}}{\longrightarrow}$}}

\put(90,-40){\line(1,1){10}}
\put(90,-40){\line(-1,1){10}}

\put(-80,-60){\makebox(0,0)[l]{Figure 0a: The corner growth model}}
\put(-80,-70){\makebox(0,0)[l]{with growth rate $q$, death rate $p$}}
\put(-80,-80){\makebox(0,0)[l]{and total asymmetry $q-p=\gamma>0$}}
\end{picture}


\setlength{\unitlength}{1.3pt}
\begin{picture}(200,100)(-120,-125)
\linethickness{1pt}
\put(0,0){\line(1,0){140}}

\put(10,0){\line(0,1){10}}
\put(20,0){\line(0,1){10}}
\put(10,10){\line(1,0){10}}

\put(10,10){\line(0,1){10}}
\put(20,10){\line(0,1){10}}
\put(10,20){\line(1,0){10}}

\put(10,20){\line(0,1){10}}
\put(20,20){\line(0,1){10}}
\put(10,30){\line(1,0){10}}

\put(10,30){\line(0,1){10}}
\put(20,30){\line(0,1){10}}
\put(10,40){\line(1,0){10}}

\put(10,40){\line(0,1){10}}
\put(20,40){\line(0,1){10}}
\put(10,50){\line(1,0){10}}

\put(20,0){\line(0,1){10}}
\put(30,0){\line(0,1){10}}
\put(20,10){\line(1,0){10}}

\put(20,10){\line(0,1){10}}
\put(30,10){\line(0,1){10}}
\put(20,20){\line(1,0){10}}

\put(20,20){\line(0,1){10}}
\put(30,20){\line(0,1){10}}
\put(20,30){\line(1,0){10}}

\put(30,0){\line(0,1){10}}
\put(40,0){\line(0,1){10}}
\put(30,10){\line(1,0){10}}

\put(30,10){\line(0,1){10}}
\put(40,10){\line(0,1){10}}
\put(30,20){\line(1,0){10}}

\put(30,20){\line(0,1){10}}
\put(40,20){\line(0,1){10}}
\put(30,30){\line(1,0){10}}

\put(30,30){\line(0,1){10}}
\put(40,30){\line(0,1){10}}
\put(30,40){\line(1,0){10}}

\put(30,40){\line(0,1){10}}
\put(40,40){\line(0,1){10}}
\put(30,50){\line(1,0){10}}

\put(40,0){\line(0,1){10}}
\put(50,0){\line(0,1){10}}
\put(40,10){\line(1,0){10}}

\put(40,10){\line(0,1){10}}
\put(50,10){\line(0,1){10}}
\put(40,20){\line(1,0){10}}

\put(40,20){\line(0,1){10}}
\put(50,20){\line(0,1){10}}
\put(40,30){\line(1,0){10}}

\put(50,52){\line(0,1){10}}
\put(60,52){\line(0,1){10}}
\put(50,62){\line(1,0){10}}
\put(50,52){\line(1,0){10}}

\put(54,57){\makebox(0,0)[l]{$\downarrow$}}
\put(63,57){\makebox(0,0)[l]{\tiny{Poisson process of falling blocks}}}
\color{Gray}
\put(50,30){\line(0,1){10}}
\put(60,30){\line(0,1){10}}
\put(50,40){\line(1,0){10}}
\put(50,30){\line(1,0){10}}
\color{black}
\put(60,0){\line(0,1){10}}
\put(70,0){\line(0,1){10}}
\put(60,10){\line(1,0){10}}

\put(60,10){\line(0,1){10}}
\put(70,10){\line(0,1){10}}
\put(60,20){\line(1,0){10}}

\put(60,20){\line(0,1){10}}
\put(70,20){\line(0,1){10}}
\put(60,30){\line(1,0){10}}

\put(60,30){\line(0,1){10}}
\put(70,30){\line(0,1){10}}
\put(60,40){\line(1,0){10}}

\put(70,0){\line(0,1){10}}
\put(80,0){\line(0,1){10}}
\put(70,10){\line(1,0){10}}

\put(80,0){\line(0,1){10}}
\put(90,0){\line(0,1){10}}
\put(80,10){\line(1,0){10}}

\put(80,10){\line(0,1){10}}
\put(90,10){\line(0,1){10}}
\put(80,20){\line(1,0){10}}

\put(80,20){\line(0,1){10}}
\put(90,20){\line(0,1){10}}
\put(80,30){\line(1,0){10}}

\put(80,30){\line(0,1){10}}
\put(90,30){\line(0,1){10}}
\put(80,40){\line(1,0){10}}

\put(80,40){\line(0,1){10}}
\put(90,40){\line(0,1){10}}
\put(80,50){\line(1,0){10}}

\put(90,0){\line(0,1){10}}
\put(100,0){\line(0,1){10}}
\put(90,10){\line(1,0){10}}

\put(90,10){\line(0,1){10}}
\put(100,10){\line(0,1){10}}
\put(90,20){\line(1,0){10}}

\put(90,20){\line(0,1){10}}
\put(100,20){\line(0,1){10}}
\put(90,30){\line(1,0){10}}

\put(100,0){\line(0,1){10}}
\put(110,0){\line(0,1){10}}
\put(100,10){\line(1,0){10}}

\put(100,10){\line(0,1){10}}
\put(110,10){\line(0,1){10}}
\put(100,20){\line(1,0){10}}

\put(100,20){\line(0,1){10}}
\put(110,20){\line(0,1){10}}
\put(100,30){\line(1,0){10}}

\put(110,0){\line(0,1){10}}
\put(120,0){\line(0,1){10}}
\put(110,10){\line(1,0){10}}

\put(110,10){\line(0,1){10}}
\put(120,10){\line(0,1){10}}
\put(110,20){\line(1,0){10}}

\put(110,20){\line(0,1){10}}
\put(120,20){\line(0,1){10}}
\put(110,30){\line(1,0){10}}

\put(110,30){\line(0,1){10}}
\put(120,30){\line(0,1){10}}
\put(110,40){\line(1,0){10}}

\put(110,40){\line(0,1){10}}
\put(120,40){\line(0,1){10}}
\put(110,50){\line(1,0){10}}

\put(120,0){\line(0,1){10}}
\put(130,0){\line(0,1){10}}
\put(120,10){\line(1,0){10}}
\put(55,-5){\makebox(0,0){$x$}}
%


\put(10,-60){\line(1,0){30}}

\put(10,-60){\line(0,1){10}}
\put(20,-60){\line(0,1){10}}
\put(10,-50){\line(1,0){10}}

\put(10,-50){\line(0,1){10}}
\put(20,-50){\line(0,1){10}}
\put(10,-40){\line(1,0){10}}

\put(10,-40){\line(0,1){10}}
\put(20,-40){\line(0,1){10}}
\put(10,-30){\line(1,0){10}}

\put(30,-60){\line(0,1){10}}
\put(40,-60){\line(0,1){10}}
\put(30,-50){\line(1,0){10}}

\put(20,-20){\line(0,1){10}}
\put(30,-20){\line(0,1){10}}
\put(20,-10){\line(1,0){10}}
\put(20,-20){\line(1,0){10}}
\color{Gray}
\put(24,-15){\makebox(0,0)[l]{$\downarrow$}}

\put(20,-40){\line(0,1){10}}
\put(30,-40){\line(0,1){10}}
\put(20,-30){\line(1,0){10}}
\put(20,-40){\line(1,0){10}}

\color{black}


\put(60,-60){\line(1,0){30}}

\put(60,-60){\line(0,1){10}}
\put(70,-60){\line(0,1){10}}
\put(60,-50){\line(1,0){10}}

\put(60,-50){\line(0,1){10}}
\put(70,-50){\line(0,1){10}}
\put(60,-40){\line(1,0){10}}

\put(70,-60){\line(0,1){10}}
\put(80,-60){\line(0,1){10}}
\put(70,-50){\line(1,0){10}}

\put(70,-50){\line(0,1){10}}
\put(80,-50){\line(0,1){10}}
\put(70,-40){\line(1,0){10}}

\put(80,-60){\line(0,1){10}}
\put(90,-60){\line(0,1){10}}
\put(80,-50){\line(1,0){10}}

\put(70,-20){\line(0,1){10}}
\put(80,-20){\line(0,1){10}}
\put(70,-10){\line(1,0){10}}
\put(70,-20){\line(1,0){10}}
\color{Gray}
\put(74,-15){\makebox(0,0)[l]{$\downarrow$}}

\put(70,-40){\line(0,1){10}}
\put(80,-40){\line(0,1){10}}
\put(70,-30){\line(1,0){10}}
\put(70,-40){\line(1,0){10}}

\color{black}


\put(110,-60){\line(1,0){30}}

\put(110,-60){\line(0,1){10}}
\put(120,-60){\line(0,1){10}}
\put(110,-50){\line(1,0){10}}

\put(110,-50){\line(0,1){10}}
\put(120,-50){\line(0,1){10}}
\put(110,-40){\line(1,0){10}}

\put(130,-60){\line(0,1){10}}
\put(140,-60){\line(0,1){10}}
\put(130,-50){\line(1,0){10}}

\put(130,-50){\line(0,1){10}}
\put(140,-50){\line(0,1){10}}
\put(130,-40){\line(1,0){10}}

\put(130,-40){\line(0,1){10}}
\put(140,-40){\line(0,1){10}}
\put(130,-30){\line(1,0){10}}

\put(120,-20){\line(0,1){10}}
\put(130,-20){\line(0,1){10}}
\put(120,-10){\line(1,0){10}}
\put(120,-20){\line(1,0){10}}
\color{Gray}
\put(124,-15){\makebox(0,0)[l]{$\downarrow$}}

\put(120,-40){\line(0,1){10}}
\put(130,-40){\line(0,1){10}}
\put(120,-30){\line(1,0){10}}
\put(120,-40){\line(1,0){10}}

\color{black}
\put(10,-70){\makebox(0,0)[l]{Figure 0b: Ballistic deposition model.}}
\put(10,-80){\makebox(0,0)[l]{Growth occurs when blocks stick to}}
\put(10,-90){\makebox(0,0)[l]{first point of contact (denoted in \color{Gray} gray\color{black}).}}
\end{picture}\label{fig1}
\end{figure}
\vskip -.4in

\subsubsection{Kardar, Parisi and Zhang's prediction}\label{KPZpred}

The KPZ universality class was introduced in the context of studying the motion of growing interfaces in a 1986 paper of Kardar, Parisi and Zhang \cite{KPZ} which has since been cited thousands of times in both the mathematics and physics literature. Employing Forster, Nelson and Stephen's \cite{FNS} 1977 dynamical renormalization group techniques (highly non-rigorous from a mathematical perspective), \cite{KPZ} predicted that scaling exponent of $1/3$ and $2/3$ should describe the fluctuations and correlations for a large class of models, stable under varying model parameters such as underlying probability distributions or local rules.

These exponents had first been identified in 1977 by \cite{FNS} in the study of the stochastic Burgers equation (stochastic space-time noise, not initial data as studied in \cite{Burgers}). In 1985 the model of a directed polymer in random media was first formulated in \cite{HuHe} and in the same year the scaling exponents for the driven lattice gas (the asymmetric simple exclusion process) were discovered by \cite{vBKS}.\footnote{In fact, \cite{vBKS} computed in approximation the stationary two-point function for the ASEP though in simulations full distributions were out of reach at this time.} The connections between polymers and lattice gases were understood quickly \cite{Klet,HuHeFi}. The big step of Kardar, Parisi and Zhang was therefore to relate these models and calculations to interface motion in arbitrary dimensions. This both broadened the universality class of models which shared characteristic exponents as well as provided obvious physical realizations of the universality class. Their observations caught the imagination of many and the paper deserves the full credit of having triggered a large number of investigations. Examples of physical phenomena modeled by the KPZ class include turbulent liquid crystals \cite{TS}, crystal growth on a thin film \cite{WZG}, facet boundaries \cite{DSC}, bacteria colony growth \cite{WIMM,MWI}, paper wetting \cite{KHO}, crack formation \cite{EJHS}, and burning fronts \cite{MMA,MMK,MMM}.

The scaling exponent predictions of \cite{KPZ} were based on studying a continuum stochastically growing height function $\mathcal{H}(T,X)$ given in term of a stochastic PDE (ill-posed however) which is now known as the {\it KPZ equation}. The time derivative of the height function depends on three factors (the same as highlighted for the ballistic deposition model): smoothing (the Laplacian), rotationally invariant, slope dependent, growth speed (the square of the gradient), noise (space-time white noise).
\begin{equation}
\partial_T \mathcal{H} = \nu\partial_X^2 \mathcal{H} +\tfrac{1}{2}\lambda(\partial_X \mathcal{H})^2 +\sqrt{D}\dot{\mathscr{W}}.
\end{equation}
Here $\dot{\mathscr{W}}$ is space-time white noise and $\nu,\lambda$ and $D$ are non-zero parameters which can often be (heuristically) computed for a particular growth model directly from the microscopic dynamics. When one of these parameters is computed to be zero, this is indicative of the growth model not being in the KPZ class (for instance when $\lambda=0$ this is the EW class). We will specialize the parameters to the values $\nu=\tfrac{1}{2}$, $\lambda=-1$ and $D=1$ and thus refer to the KPZ equation as
\begin{equation}\label{KPZ}
\partial_T \mathcal{H} = \tfrac{1}{2}\partial_X^2 \mathcal{H} - \tfrac{1}{2}(\partial_X \mathcal{H})^2 +\dot{\mathscr{W}}.
\end{equation}
However, nothing is lost in doing this since the general solution $\mathcal{H}_{\nu,\lambda,D}$ (the subscripts emphasize the coefficients) to the KPZ equation can be recovered by a combination of change of variables and rescaling as
\begin{equation}
\mathcal{H}_{\nu,\lambda,D} = -\frac{\lambda^2 D}{(2\nu)^3} \mathcal{H}_{\tfrac{1}{2},-1,1} \left(\frac{\lambda^2 D}{(2\nu)^3} x, - \frac{\lambda^4 D^2}{(2\nu)^5} t \right).
\end{equation}

It has been argued in the physics literature that the long time behavior of this equation should be the same as that of a variety of physical systems as well as mathematical models which share these three features -- hence the belief in the wide universality of the KPZ class. Besides universality of the scaling exponents, there quickly developed a belief that the asymptotic long time scaling limit probability distributions (also called amplitudes) were universal (within certain geometry dependent subclasses).

In the years following \cite{KPZ} it became a very hot subject in the physics literature to try to both demonstrate the scaling exponents as well as compute or estimate properties of these universal scaling limit probability distributions. Many (mathematically non-rigorous) approaches -- both analytical and numerical -- were employed to gather more information about the KPZ universality class and its inhabitants. These included replica methods, matrix models, field theory methods, saddle point equations, mode coupling equations, perturbation theory, and Monte-Carlo simulations (see \cite{HHZ} for references and examples of many of these approaches).

These early papers worked to establish in broad terms and via brute force numerical methods, the full notation of KPZ universality. Most of these focused on flat geometries (as it was easier to simulate) and calculated statistics for the height function fluctuations for various models. Via Monte-Carlo methods the skew and kurtosis were approximated in \cite{THH,KMHH} and in \cite{KBM} a rough plot of the distribution function was given for both a growth model and a polymer model (and they agreed).

\subsubsection{Exact statistics breakthrough for the KPZ class}

Still, an exact and analytic description of the statistics (probability distributions) associated with the KPZ class went entirely unknown until, in the late 1990s, a group of mathematicians determined the exact formula for the one-point statistics of the KPZ class (in the wedge growth geometry corresponding to Figure 0a). That seminal work of Baik, Deift and Johansson \cite{BDJ, KJ} dealt with two closely related discrete models (polynuclear, and corner growth) predicted to have the KPZ scaling.\footnote{While \cite{KJ} dealt directly with the corner growth model, \cite{BDJ} considered the fluctuations of the longest increasing subsequence of a random permutation (Ulam's problem). The Poissonize version of Ulam's problem is related to an interacting particle system known as the Hammersley process \cite{Ham,AD} (see \cite{OdlyzkoRains,BaerBrock} for relevant numerical experiments). It was Pr\"{a}hofer and Spohn \cite{PS0} who recognized the connection of poissonized Ulams problem to the polynuclear growth model.} By using exact formulas (arising from combinatorics and representation theory) and then by studying asymptotics of the resulting expressions, \cite{BDJ,KJ} were able to calculate the statistics of the one-point fluctuations for these KPZ class models.\footnote{\cite{BDJ} studied the combinatorial problem in terms of Toeplitz matrices (a reduction due to Gessel \cite{Gessel}) and the associated Riemann-Hilbert problem method. Their asymptotics therefore were those of Riemann-Hilbert steepest descent and their limiting formula's were in terms of Riemann-Hilbert problems and hence Painlev\'{e} expressions. \cite{KJ} derives a Fredholm determinant and thus is able to perform classical steepest descent directly on the integral kernel, resulting in Fredholm determinant formulas for the limiting statistics.} Surprisingly these statistics had been discovered in the early 90s by Tracy and Widom in the context of random matrix theory \cite{TW}. The previous numerical work of the 1990s agreed with the values readily computed from the exact formulas \cite{PS0}. Experimental work has shown that the scalings and the statistics for the KPZ class are excellent fits for certain physical phenomena. Of note is the recent work of \cite{TS} on liquid crystal growth (see also \cite{WZG, MMK, MMM} for other experimental evidence).

Even before the work of Baik, Deift and Johansson, the KPZ class and equation were topics of interest in mathematics. However, after their work the field blossomed (for instance there have been two semester long MSRI programs focused in  this direction). In the past ten years there has been a significant amount of refinement of the theory of the KPZ universality class. For instance, it is now understood that while the $1/3$ and $2/3$ scaling exponents occur for all models, the exact statistics of KPZ class models fall into different subclasses based on the growth geometry, or initial / boundary data (see Figure \ref{sixfig}). The long-time statistics have been written down explicitly for six of these subclasses which, arguably are the most important phenomenalogically.

\subsubsection{Exact statistics breakthrough for the KPZ equation}

One should be clear now that all of the mathematically rigorous results until recently involving exact statistics have been for models which are believed to be in the KPZ universality class, but not for the KPZ equation itself. In fact, before 2010 very little was known of the exact statistics of the solution to the KPZ equation itself.

One stumbling block is that the KPZ equation is mathematically ill-posed due to its non-linearity (the function $\mathcal{H}$ is believed to be locally Brownian and so it doesn't make sense to square its derivative). However, this issue was essentially resolved in the mid-1990s by Bertini-Giacomin \cite{BG} who provided a convincing interpretation for what it means to solve the KPZ equation. First they formally defined the {\it Hopf-Cole solution} to the KPZ equation as $$\mathcal{H}(T,X):=-\log \mathcal{Z}(T,X)$$
where $\mathcal{Z}(T,X)$ satisfies the well-posed SHE (stochastic heat equation)
\begin{equation}\label{SHE}
\partial_T \mathcal{Z} = \tfrac{1}{2}\partial_X^2 \mathcal{Z} -\mathcal{Z}\dot{\mathscr{W}}.
\end{equation}
They then proved that the discrete Hopf-Cole transform of the height function for a certain KPZ class model (the corner growth model) converged weakly as a space-time process to the solution to the stochastic heat equation. Their work required a special sort of scaling which is called {\it weakly asymmetric} scaling, and only applied to the growth model started close to its stationary distribution.\footnote{The stationary distribution is a height profile given by the trajectory of a simple symmetric random walk. Strictly speaking this is not stationary though its discrete derivative (interpreted as particles and holes) is stationary in time. Within the probability literature (as well as in \cite{BG}) this is called ``equilibrium'' initial data, though we attempt to avoid this term. The reason being that this system is a paradigm for non-equilibrium statistical physics and thus using equilibrium to describe initial data would likely cause un-needed confuse.} The interpretation of this result, however, is strikingly clear -- up to this odd transformation and scaling, the discrete growth model converges to the KPZ equation. This Hopf-Cole interpretation had been used previously (though without the justification provided by \cite{BG}) in the physics literature for some time. Other attempts at interpreting the KPZ equation (such as \cite{HOUZ}) have resulted in solutions which are considered to be physically irrelevant as their solutions have very different long-time scaling behavior \cite{TC}.

Since the KPZ equation is a fundamental mathematical / physical object, one would hope to be able to write down and prove the exact probability distribution of its solution and show that its long time limit is governed by the KPZ class scaling and statistics -- i.e., that the {\bf KPZ equation is in the KPZ universality class.}

This long-standing goal was achieved in the fall of 2009 (and posted in spring 2010). Two groups (independently and in parallel) derived the exact formula for the one-point statistics for the solution to the KPZ equation \cite{ACQ,SaSp1}. In addition to a derivation, \cite{ACQ} provided the highly non-trivial mathematical proof of the formula.

\begin{theorem}[\cite{ACQ}]\label{ACQmainthm}
For any $T>0$ and $X\in \R$, the Hopf-Cole solution to KPZ with {\it narrow wedge} initial data (given by $\mathcal{H}(T,X)= -\log\mathcal{Z}(T,X)$ with initial data $\mathcal{Z}(0,X)=\delta_{X=0}$) has the following probability distribution:
\begin{equation}
\PP(\mathcal{H}(T,X)-\frac{X^2}{2T} - \frac{T}{24} \geq -s) = F_{T}(s)
\end{equation}
where $F_T(s)$ does not depend on $X$ and is given by
\begin{equation}\label{fandef}
F_{T}(s) = \int_{C} \frac{d\mu}{\mu} e^{-\mu} \det(I-K_{\sigma_{T,\mu}})_{L^2(\kappa_T^{-1} s,\infty)}
\end{equation}
where $\kappa_T = 2^{-1/3} T^{1/3}$, $C$ is a contour positively oriented and going from $+\infty+\e i$ around $\R^+$ to $+\infty-i\e$, and $K_{\sigma}$ is an operator given by its integral kernel
\begin{equation}
K_{\sigma}(x,y) = \int_{-\infty}^{\infty} \sigma(t) \Ai(x+t)\Ai(y+t)dt, \qquad \textrm{and} \qquad \sigma_{T,\mu}(t)= \frac{\mu}{\mu - e^{-\kappa_T t}}.
\end{equation}
\end{theorem}

\begin{figure}
\begin{center}
\includegraphics[scale=1.5]{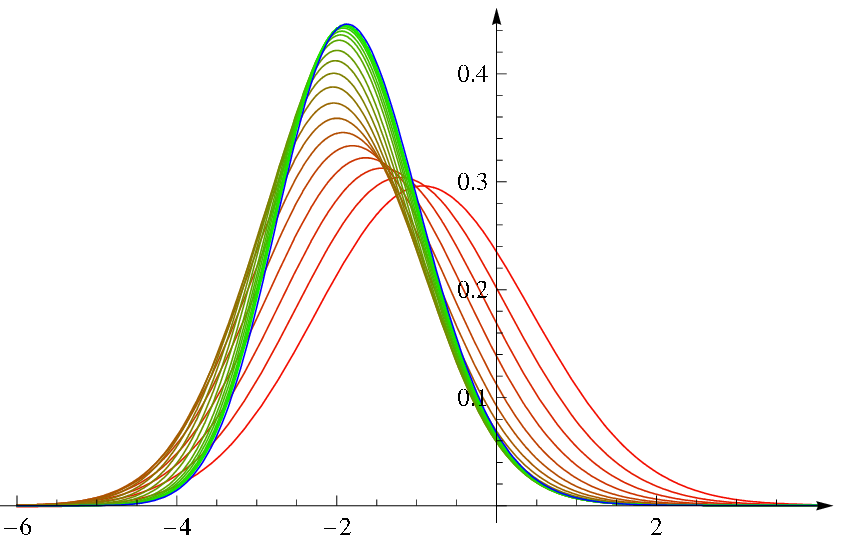}
\caption[Crossover distributions]{Properly scaled probability density functions for $F_{T}(s)$ showing the crossover to $F_{{\rm GUE}}$ in long time (green, upper curves) and Gaussian in short time (red, lower curves). These plots have been made by evaluating the Fredholm determinant expressions for $F_{T}(s)$ (in fact in a slightly different form than quoted in this survey -- see for instance the Gumbel convolution formula of \cite{ACQ} Theorem 1). Recent advances in numerical calculations of Fredholm determinants by F. Bornemann \cite{Born} have enabled fast calculations of these plots (as well as many other distribution functions arising in random matrix theory and the study of the KPZ universality class). Special thanks to S. Prolhac and H. Spohn for permission to reproduce these plots from \cite{ProS3}.}\label{ProSfig}
\end{center}
\end{figure}

Given this explicit formula it is then a  relatively easy corollary that as $T$ goes to infinity, under $T^{1/3}$ scaling, the statistics of the KPZ equation converge to those of Baik, Deift and Johansson (the Tracy-Widom $F_{{\rm GUE}}$ distribution) -- thus the KPZ equation is in the KPZ universality class! The short time ($T$ goes to zero) statistics scale like $T^{1/4}$ and converge to those of the EW (Edwards Wilkinson) class \cite{EW} -- which is governed by the {\it additive stochastic heat equation} and hence has a Gaussian one-point distribution. The short time results are more easily derived directly from the chaos expansion for the multiplicative stochastic heat equation. This shows, in fact, that the KPZ equation actually represents a mechanism for crossing over between two universality classes -- the KPZ class in long-time and the EW (Edwards Wilkinson) class in short-time (see also the plots of Figure \ref{ProSfig})

\begin{corollary}
The Hopf-Cole solution to the KPZ equation with narrow wedge initial data has the following long-time and short-time asymptotics:
\begin{equation}
F_{T}(2^{-1/3}T^{1/3} s) \stackrel{T\rightarrow\infty}{\longrightarrow} F_{{\rm GUE}}(s) \qquad
F_{T}(2^{-1/2}\pi^{1/4}T^{1/4} (s-\log\sqrt{2\pi T})) \stackrel{T\rightarrow 0}{\longrightarrow} G(s).
\end{equation}
\end{corollary}

The GUE distribution function was first discovered in \cite{TW} and therein was expressed in terms of a Fredholm determinant as well as in terms of the Hastings-McLeod solution to the non-linear ODE known as Painlev\'{e} II
\begin{equation}
F_{\rm{GUE}}(s)= \det(I-K_{\rm{Ai}})_{L^2(s,\infty)} = \exp\left(-\int_s^{\infty} (x-s)q^2(x)dx\right)
\end{equation}
where $q(x)$ solves the ODE
\begin{equation}
q''(x) = (x+2q^2(x))q(x)
\end{equation}
subject to $q(x)\sim \rm{Ai}(x)$ as $x\rightarrow \infty$. Also
\begin{equation}
K_{\rm{Ai}}(x,y) = \int_{-\infty}^{\infty} {\bf 1}_{t\geq 0} \Ai(x+t)\Ai(y+t) dt
\end{equation}
and hence is of the form of the $K_{\sigma}$ with $\sigma = {\bf 1}_{t\geq 0}$. $G(s)$ is the Gaussian distribution
\begin{equation}
G(s) = \int_{-\infty}^{s} \frac{1}{\sqrt{ 2\pi}} e^{-r^2/2}dr.
\end{equation}

It is important to note that the $T^{2/3}$ correlation scaling for the KPZ equation still remains unproved. This is because the methods used above for the KPZ equation presently only give one-point fluctuation results, and one would require control over the two-point correlation in order to show this scaling. A proof of this and moreover an exact expression for the multi-point correlation functions remains an important goal.

\subsubsection{The KPZ equation as a critical scaling object}\label{critscalingobject}
In order to properly understand the above crossover we shall see that the KPZ equation arises as the scaling limit of the corner growth model under weak asymmetry (we will soon specify in Section \ref{asymsym} exactly what this means). In Section \ref{asymsym} we will observe, however, for any positive asymmetry the corner growth model is in the KPZ class, while for zero asymmetry it is in the EW class.

This entire picture of the KPZ equation being the weakly asymmetric limit of the corner growth model, and having statistics which crossover between universality classes just further confirms a mantra in statistical physics that at critical points (in this case the interface between universality classes) one can expect to see universal scaling objects. In this case that object is the KPZ equation, though other examples abound. Analogous to our case is the fact that the scaling limit for critical percolation (on the hexagonal lattice) is governed by SLE6.

To motivate this whole pursuit and approach, it is useful to draw upon a rather simple analogy with the story of random walks and Brownian motion. There are a wide (infinite) variety of random walks. However, assuming the jumps have finite second moment, if you subtract off their linear drift, and scale properly, they all converge to a single (hence universal) continuum object -- the Brownian motion.

In that case, the benefit of having a single limiting object is enormous. As long as one is dealing with a suitably long time-scale, the above result shows that all such random walks essentially look similar. Thus, if we compute properties of the Brownian motion, these will (with a little care) imply similar (and exact in the limit) results for the random walks. This observation can be thought of as the basis for Gaussian statistics. Because of its importance, the properties of Brownian motion have been incessantly studied and results pertaining to its properties fills countless volumes of books and research articles. In particular, a calculus was developed for the Brownian motion.

When it comes to proving properties of the Brownian motion there are many approaches.  However, it is informative to recall that the earliest approach, taken by de Moivre, was to consider a particularly {\it solvable discrete model} -- the simple symmetric random walk. In that setting he took asymptotics of Binomial coefficients to compute the limiting Gaussian distribution. Even though now the techniques have evolved past relying on such simple models to prove properties of the Brownian motion, early on such solvable models were essential.

Presently, the state-of-the-art in solving for the statistics of the KPZ equation is analogous to taking asymptotics of Binomial coefficients. It is only through exactly solvable discrete models (which can be proved to converge to the KPZ equation) that one is able to compute {\it anything} about the KPZ equation rigorously. As far as universality of the KPZ equation as a scaling limit, few doubt it, but presently it is only known rigorously to arise in the particular case of the corner growth model. However, the KPZ equation is essentially equivalent to the free energy of a continuum directed random polymer model, and within the study of directed polymers, it is simpler to prove universality of the continuum model (see Section \ref{freeECDRP}).

\subsubsection{The proof of Theorem \ref{ACQmainthm} in two steps}
There are two main steps involved in the proof of the above theorem -- the first is a universality result and the second a solvability result (these are both explained in greater depth in Section \ref{talk1} and \ref{talk2} respectively). The first step is to rigorously relate the corner growth model to the Hopf-Cole solution to KPZ with narrow-wedge initial data. Bertini-Giacomin's earlier work \cite{BG} {\it does not apply} to the wedge geometry since it is very far from stationary initial data. In fact, a new term appears in the scaling and in \cite{ACQ} it is necessary to provide a proof that under this corrected scaling, the KPZ equation arises (the proof of this result certainly draws on the work of \cite{BG} though).

To fix notation, $h_{\gamma}(t,x)$ represents the height function for the corner growth model at time $t\geq 0$ above position $x\in \Z$. The $\gamma$ is the asymmetry and is given by $\gamma = q-p$ where $q$ is the growth rate and $p$ is the death rate (see Figure 0a or Section \ref{IPSsec}). Wedge initial condition corresponds to $h_{\gamma}(0,x) =|x|$.

This first step is summarized in the following:

\begin{theorem}[Theorem 1.14 of \cite{ACQ}]\label{ACQBGthm}
Fix wedge initial conditions. For $\e>0$ set $\gamma=\e^{1/2}$, $\nu_{\e}= p+q-2\sqrt{qp}\approx \tfrac{1}{2}\e$ and $\lambda_{\e} = \tfrac{1}{2}\log(q/p)\approx \e^{1/2}$ (recall $q+p=1$ and $q-p=\gamma$) and define the discrete Hopf-Cole transformed height function as
\begin{equation}
Z_{\e}(T,X) = \tfrac{\e^{-1}}{2} \exp\{-\lambda_\e h_{\gamma}(\tfrac{t}{\gamma},x)+\nu_\e \tfrac{t}{\gamma}\}.
\end{equation}
Let $\PP_{\e}$ denote the probability measure on space-time processes in $D_{u}([0,\infty);D_u(\R))$ (here $D_u$ refers to right continuous paths with left limits in the topology of uniform convergence on compact sets). Then $\PP_\e$, for $\e\in(0,1/4)$, are a tight family of measures and the unique limit point is supported on $C((0,\infty);C(\R))$ (continuous in both space and time) and corresponds to the solution of the SHE with delta function initial data $\mathcal{Z}(0,X) = \delta_{X=0}$.
\end{theorem}

The second step is to use the recently discovered exact formulas of Tracy and Widom \cite{TW1,TW2,TW3} for the corner growth model height function, and take their weakly asymmetric limit.

\begin{theorem}[Tracy-Widom ASEP formula \cite{TW3}]\label{TW}
Consider the corner growth model with wedge initial conditions and with $q>p$ such that $q+p=1$. Let $\gamma=q-p$ and $\tau=p/q$. For $m=\lfloor \tfrac{1}{2}(s+x)\rfloor$, $t\geq 0$ and $x\in \Z$
\begin{equation}
P(h_{\gamma}(t,x)\geq s) = \int_{S_{\tau^+}}\frac{d\mu}{\mu} \prod_{k=0}^{\infty} (1-\mu\tau^k)\det(I+\mu J_{t,m,x,\mu})_{L^2(\Gamma_{\eta})}
\end{equation}
where $S_{\tau^+}$ is a positively oriented circle centered at zero of radius strictly between $\tau$ and 1, and where the kernel of the determinant is given by
\begin{equation}
J_{t,m,x,\mu}(\eta,\eta')=\int_{\Gamma_{\zeta}} \exp\{\Psi_{t,m,x}(\zeta)-\Psi_{t,m,x}(\eta')\}\frac{f(\mu,\zeta/\eta')}{\eta'(\zeta-\eta)}d\zeta
\end{equation}
where $\eta$ and $\eta'$ are on $\Gamma_{\eta}$, a circle centered at zero of radius strictly between $\tau$
and $1$, and the $\zeta$ integral is on $\Gamma_{\zeta}$, a circle centered at zero of radius strictly between $1$ and $\tau^{-1}$, and where
\begin{eqnarray}
\nonumber f(\mu,z)&=&\sum_{k=-\infty}^{\infty} \frac{\tau^k}{1-\tau^k\mu}z^k\\
\Psi_{t,m,x}(\zeta) &=& \Lambda_{t,m,x}(\zeta)-\Lambda_{t,m,x}(\xi)\\
 \nonumber \Lambda_{t,m,x}(\zeta) &=& -x\log(1-\zeta) + \frac{t\zeta}{1-\zeta}+m\log\zeta
\end{eqnarray}
\end{theorem}

It is relatively straight-forward to guess the correct limiting formula (and hence guess the exact formula for the KPZ equation) as we demonstrate in Section \ref{WASEPsec}. However, to make these asymptotics rigorous requires a significant amount of work. Besides the usual complications involved in proving uniformity and tail bounds (for the purpose of trace-class convergence), there is a major technical issue one encounters in the asymptotics. Briefly put, the heart of the asymptotics involves Taylor expanding the function $\Psi$ around a critical point $\xi$ and showing that, in rescaled coordinates, along one contour ($\Gamma_{\eta'}$) $\Psi$ rapidly goes to $+\infty$ away from $\xi$ and along the other contour ($\Gamma_{\zeta}$) $\Psi$ rapidly goes to $-\infty$ away from $\xi$. To facilitate this, one generally wants to deform the two contours to be very far from each other (for instance one contour leaves $\xi$ at angle $\pm \pi/3$ and the other at $\pm 2\pi/3$). However, due to the weak asymmetry scaling, the two contours of interest can not be deformed as such since doing so one would cross a diverging number of poles as indicated by the restriction from the function $f$ that $\zeta/\eta'\in (1,\tau^{-1})$. Therefore, it is necessary to prove that despite staying extremely close together, one can still choose contours which have the desired opposing decay properties. A priori it is far from clear that this is possible. \cite{ACQ} explicitly shows it is possible and give the necessary contours and estimates.

\subsubsection{The chronology}\label{chrono}
Before going further we should remark on the chronology of the solution to the KPZ equation (see also Sasamoto and Spohn's contributions to StatPhys 24 \cite{SaSp4}). The formula was independently and in parallel discovered in mathematics by Amir-Corwin-Quastel \cite{ACQ} and mathematical physics by Sasamoto-Spohn \cite{SaSp1,SaSp2,SaSp3} in the fall of 2009 and posted early in 2010. The derivation (unsurprisingly) relied on a similar approach.

Amir-Corwin-Quastel's work contains a complete and rigorous proof of the exact formula (and in particular dealt carefully with the technically involved asymptotic analysis of the Tracy-Widom corner growth model solution mentioned above) and the necessary extension of the work of Bertini-Giacomin away from near-stationary initial data. Additionally, they showed that the solution satisfies a certain extension of the Painlev\'{e} II equation, suggesting a relationship to integrable systems.

Sasamoto-Spohn derive the exact formula as well as provide detailed numerical plots of the newly discovered statistics, compute the long-time correction to the Tracy-Widom GUE distribution, notice certain connections with $q$ combinatorial identities and explain the physical context of the KPZ equation in detail. Certain points in the derivation, however, proceed without rigorous mathematical justification. The first is in the application of Bertini-Giacomin's convergence theorem (which does not apply in this setting and must be replaced by \cite{ACQ}, Theorem 1.3 -- resulting in a $\log{\epsilon}$ shift in the height function). Instead, \cite{SaSp1} determine that this should be the necessary shift by using a first moment fitting argument. The second point is in freely deforming contours during asymptotic analysis (when in fact there is a very restrictive condition on contour manipulations imposed by the critical scaling of the formulas of Tracy and Widom). The deformation through many poles should introduce a diverging correction to their formula which, however in later manipulations, is luckily canceled and results in a correct final formula.

At the same time, a very different (and mathematically non-rigorous) approach was pursued by two groups of physicists, Dotsenko \cite{Dot} and Calabrese-Le Doussal-Rosso \cite{CDR}. The $n^{th}$ moment of $\mathcal{Z}$ satisfies a closed evolution equation which coincide with the quantum many body system known as the attractive $\delta$-Bose gas on $\R$ with $n$ particles. The eigenfunctions for this system were found by McGuire \cite{McGuire} in 1964, and the norms of these eigenfunctions were determined using the algebraic Bethe Ansatz (for instance \cite{KirKor,BogIzer,Slavnov}) in \cite{CalCaux} in 2007. Going from the moments of $\mathcal{Z}$ to its distribution (or the distribution of its logarithm) is mathematically unsound since one readily checks that the moments grow like $e^{cn^3}$ which means that the moment problem is ill-posed. This is a classic issue in replica trick calculations and necessitates summing highly divergent series and performing unjustified analytic continuations. In January 2010 Dotsenko \cite{DotKlum} posted an attempt at this, but admittedly got the wrong answer. Likewise, in February 2010 in the first posted version of \cite{CDR}, Calabrese, Le  Doussal and Rosso find formulas which do not have the Tracy-Widom GUE distribution as the long-time asymptotics. In late March both groups realized the correct way to sum the divergent series and get a Fredholm determinant and in the end they were able to show agreement of their formulas with those of Sasamoto-Spohn and Amir-Corwin-Quastel.\footnote{In fact, while both \cite{Dot} and \cite{CDR} write down generating functions for the KPZ equation distribution function and find the long-time asymptotics to be Tracy-Widom GUE, only \cite{CDR} extracted finite time statistics.} Despite the inherent lack of mathematical rigor in this approach, these methods have proved to be powerful and have lead to certain predictions which can not yet be confirmed rigorously (see Section \ref{replicasec} for more about this approach).

\subsection{Models in the KPZ universality class}\label{growthregimesec}
A growth model is considered to be in the KPZ universality class if its long time behavior is similar to that of the KPZ equation itself. There are a wide variety of systems which are predicted to fall into this class and herein we only consider a few of the simplest examples.

\subsubsection{Interacting particle systems and the simple exclusion process}\label{IPSsec}
In many cases, growth models can be thought of as integrated version of interacting particle systems (see \cite{Lig} for an overview of this field). Particles moving stochastically and interacting according to a specified rule system are effective ways of simulating and studying real world systems and gaining key insight into complex phenomena. Particle systems such as the ones we consider are important models for mass transport, traffic flow, queueing theory, driven lattice systems, and turbulence. We consider particle systems where individual particles attempt to orchestrate random walks (often with positive drift) with the caveat that they are influenced by their local environment (of other particles). Just as with growth processes, it is believed that varying the mechanisms for this interaction and for the random walks should not affect the limiting fluctuations of the integrated particle density (i.e., the height function). Furthermore these height fluctuations should be governed by the KPZ class scaling and statistics.

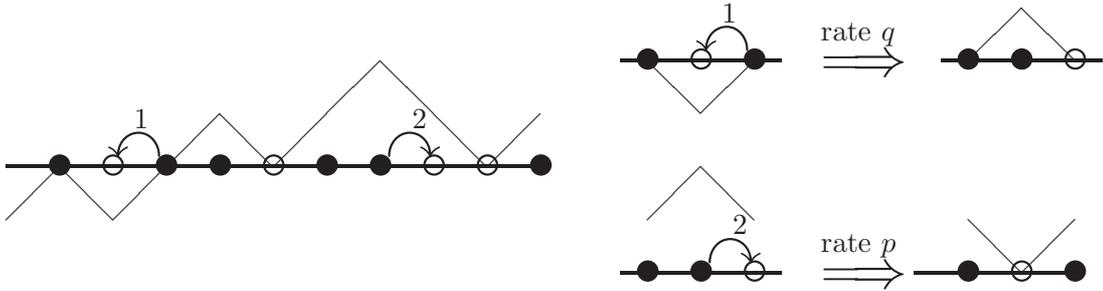
\begin{figure}
\setlength{\unitlength}{1pt}
\begin{picture}(200,100)(100,0)
\linethickness{1pt}
\put(0,40){\line(1,0){200}}
\put(0,20){\line(1,1){20}}
\put(20,40){\line(1,-1){20}}
\put(40,20){\line(1,1){20}}
\put(60,40){\line(1,1){20}}
\put(80,60){\line(1,-1){20}}
\put(100,40){\line(1,1){20}}
\put(120,60){\line(1,1){20}}
\put(140,80){\line(1,-1){20}}
\put(160,60){\line(1,-1){20}}
\put(180,40){\line(1,1){20}}

\put(15,40){\makebox(0,0)[l]{\huge{$\bullet$}}}
\put(35,40){\makebox(0,0)[l]{\huge{$\circ$}}}
\put(55,40){\makebox(0,0)[l]{\huge{$\bullet$}}}
\put(75,40){\makebox(0,0)[l]{\huge{$\bullet$}}}
\put(95,40){\makebox(0,0)[l]{\huge{$\circ$}}}
\put(115,40){\makebox(0,0)[l]{\huge{$\bullet$}}}
\put(135,40){\makebox(0,0)[l]{\huge{$\bullet$}}}
\put(155,40){\makebox(0,0)[l]{\huge{$\circ$}}}
\put(175,40){\makebox(0,0)[l]{\huge{$\circ$}}}
\put(195,40){\makebox(0,0)[l]{\huge{$\bullet$}}}

\put(38,48){\makebox(0,0)[l]{\huge{$\curvearrowleft$}}}
\put(48,58){\makebox(0,0)[l]{1}}

\put(142,48){\makebox(0,0)[l]{\huge{$\curvearrowright$}}}
\put(152,58){\makebox(0,0)[l]{2}}

\put(230,80){\line(1,0){60}}
\put(240,80){\line(1,-1){20}}
\put(260,60){\line(1,1){20}}
\put(235,80){\makebox(0,0)[l]{\huge{$\bullet$}}}
\put(255,80){\makebox(0,0)[l]{\huge{$\circ$}}}
\put(275,80){\makebox(0,0)[l]{\huge{$\bullet$}}}
\put(258,88){\makebox(0,0)[l]{\huge{$\curvearrowleft$}}}
\put(268,98){\makebox(0,0)[l]{1}}

\put(305,78){\makebox(0,0)[l]{\huge{$\Longrightarrow$}}}
\put(305,90){\makebox(0,0)[l]{rate $q$}}

\put(350,80){\line(1,0){60}}
\put(360,80){\line(1,1){20}}
\put(380,100){\line(1,-1){20}}
\put(355,80){\makebox(0,0)[l]{\huge{$\bullet$}}}
\put(375,80){\makebox(0,0)[l]{\huge{$\bullet$}}}
\put(395,80){\makebox(0,0)[l]{\huge{$\circ$}}}

\put(230,0){\line(1,0){60}}
\put(240,20){\line(1,1){20}}
\put(260,40){\line(1,-1){20}}
\put(235,0){\makebox(0,0)[l]{\huge{$\bullet$}}}
\put(255,0){\makebox(0,0)[l]{\huge{$\bullet$}}}
\put(275,0){\makebox(0,0)[l]{\huge{$\circ$}}}
\put(262,08){\makebox(0,0)[l]{\huge{$\curvearrowright$}}}
\put(272,18){\makebox(0,0)[l]{2}}

\put(305,-2){\makebox(0,0)[l]{\huge{$\Longrightarrow$}}}
\put(305,10){\makebox(0,0)[l]{rate $p$}}

\put(340,0){\line(1,0){60}}
\put(360,20){\line(1,-1){20}}
\put(380,0){\line(1,1){20}}
\put(355,0){\makebox(0,0)[l]{\huge{$\bullet$}}}
\put(375,0){\makebox(0,0)[l]{\huge{$\circ$}}}
\put(395,0){\makebox(0,0)[l]{\huge{$\bullet$}}}

\end{picture}
\caption[Corner growth model]{Integrating the spin variables of the SEP yields the corner growth model.}\label{cornergrowth}
\end{figure}

The simple exclusion process \cite{Kaw,Spi} (and its integrated version -- the corner growth model) is the poster child for all such particle systems and serve as a paradigm for non-equilibrium statistical mechanics (with thousands of articles in mathematics and physics written on it, alone). It is simple to state and in some ways exactly solvable, yet seems to contain all of the expected complexities and phenomena of a general system. Particles attempt continuous time simple random walks on $\Z$, jumping left at rate $q$ and right at rate $p=1-q$ with the caveat that jumps are suppressed if the destination site is already inhabited (see Figure \ref{cornergrowth}). The asymmetry or drift here is recorded as $$\gamma = q-p.$$ This process can be coupled to the corner growth model (see also Figure 1a) in terms of a height function $h_{\gamma}(t,x)$ given as
\begin{equation}\label{heightfunctiondef}
h_{\gamma}(t,x)= \begin{cases} 2N(t)+\sum_{0<y\leq x} \hat{\eta}(t,y), & x>0,\\
2N(t), &x=0,\\
2N(t)-\sum_{0<y\leq x} \hat{\eta}(t,y), & x<0,
\end{cases}
\end{equation}
where $N(t)$ records the net number of particles to cross from site $1$ to site $0$ in time $t$ and where $\hat{\eta}(t,x)$ equals $1$ if there is a particle at $x$ at time $t$ and $-1$ otherwise (these are called spin variables, though it is sometimes more natural to speak of occupation variables $\eta(t,x) = \tfrac{1}{2}(\hat{\eta}(t,x)+1)$ which are $1$ if there is a particle and $0$ otherwise). The dynamics for the corner growth model were previously illustrated in Figure 1a. To review them, a local valley is filled in to form a local hill at rate $q$ and likewise a local hill is removed to form a local valley at rate $p$.

Let us fix some notation. When $\gamma=1$ this process is totally asymmetric (TASEP), when $\gamma>0$ it is partial asymmetric (PASEP) and when $\gamma=0$ it is symmetric (SSEP). The case we will focus on and define soon is weakly asymmetric (WASEP). We will freely move between speaking about the particle process and growth process since they are exactly coupled as above. We will use capital letters $X$ and $T$ for variables in the continuum models (such as the KPZ equation and SHE) and small letters $x$ and $t$ for discrete models.

\subsubsection{Hydrodynamics}
Before studying fluctuations we must address the question of the long-time limit shape of the height function $h_{\gamma}(t,x)$ of the corner growth model (with $\gamma>0$) with fixed given initial conditions. The natural scaling to see a limit shape of $h$ is to take $t=\e^{-1}T$ and $x=\e^{-1}X$ and consider whether
$$\bar{h}(T,X)=\lim_{\e\rightarrow 0}\e h_{\gamma}(t/\gamma,x)$$ exists. If we assume that such a limit exists (in a suitable sense) for $T=0$ (and calling the limit $\bar{h}^0(X)$) then it is a theorem (see for instance \cite{Rez,Rost,Var}) that the limit exists for all $T>0$ and that $\bar{h}(T,X)$ is the unique weak solution to the PDE known as the inviscid Burgers equation
$$\partial_T \bar{h} =  \frac{1-(\partial_X \bar{h})^2}{2},$$
with initial data $\bar{h}^0(X)$ and subject to an entropy condition.

For instance, with a step initial condition (particles initially only at $x>0$), $\bar{h}^0(X)=|X|$ and after time $T$ the wedge shape has resolved itself in the window $x\in [-T,T]$ to equal the parabola (see Figure \ref{patrikTASEP_fig})
$$\bar{h}(T,X) = T\frac{1+(X/T)^2}{2},$$
and remains unchanged everywhere else. In particular, right above the origin it is at height $T/2$.

This is just a glimpse at the rich theory of model dependent hydrodynamical limits. It is widely believed that the fluctuations around these model dependent limits are universal (i.e. model independent) and fall into a few large universality classes characterized in terms of scaling exponents and limiting statistics.

\subsubsection{Fluctuations: asymmetry versus symmetry}\label{asymsym}

\begin{figure}
\begin{center}
\includegraphics[scale=.5]{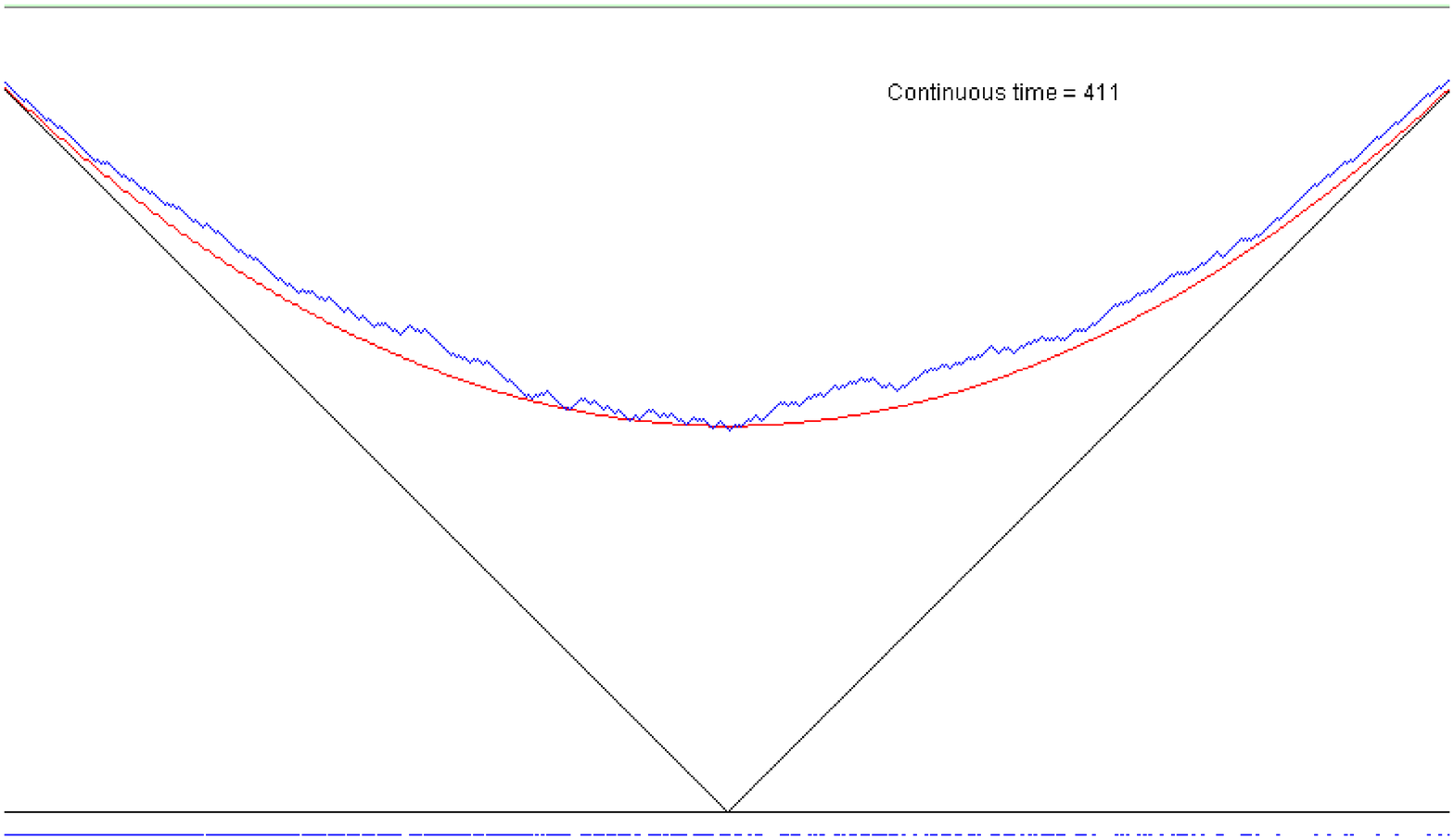}
\caption[TASEP simulation]{A simulation of the height function fluctuations for the $\gamma=1$ corner growth model started in the wedge initial condition. The curve represents the limit shape (a parabola) while the piecewise linear line represents the height function. Fluctuations live on the $t^{1/3}$ scale and are correlated spatially in the $t^{2/3}$ scale (as indicated by the box). Special thanks to Patrik Ferrari for the above simulation.}\label{patrikTASEP_fig}
\end{center}
\begin{picture}(0,0)(0,0)
\put(-15,197){\line(1,0){30}}
\put(-15,210){\line(1,0){30}}

\put(-15,197){\line(0,1){13}}
\put(15,197){\line(0,1){13}}
\end{picture}
\end{figure}

In 1985, based on physical methods known as mode-coupling, \cite{vBKS} argued that the simple exclusion process with positive asymmetry $\gamma>0$ should have height function fluctuations like $t^{1/3}$ and exhibit non-trivial spatial correlations on the $t^{2/3}$ scale. After the work of \cite{KPZ} this behavior became known as being in the KPZ universality class. On the other hand, the symmetric case $\gamma=0$ corresponds to the EW (Edwards-Wilkinson) class \cite{EW} and have fluctuations of scale $t^{1/4}$ with spatial correlation on the scale of $t^{1/2}$. For the EW class, the limiting fluctuation statistics were also predicted (and relatively easily proved -- see for example \cite{Spo}) to be Gaussian, however for the KPZ class the limiting fluctuation statistics were not found until the work of Baik, Deift and Johansson \cite{BDJ, KJ} and then Pr\"{a}hofer and Spohn \cite{PS2} (for the spatial correlation). These works only dealt with the totally asymmetric ($\gamma=1$) simple exclusion process for step initial condition (or equivalently the corner growth model in the wedge geometry) -- which is illustrated in Figure \ref{patrikTASEP_fig}. Tracy and Widom \cite{TW3} extended the one-point fluctuation results to $\gamma>0$ by way of their exact formula, recorded here as Theorem \ref{TW}. It is clear that in order to treat all values of $\gamma>0$ equivalently, we should speed up time to compensate for smaller growth asymmetry: we should take time like $t/\gamma$. The $t/2$ which is subtracted from $h_{\gamma}$ comes directly from the hydrodynamic theory. Putting together the one-point fluctuation results of \cite{KJ,TW3} we have:
\begin{theorem}\label{PASEP}
For all $\gamma\in (0,1]$ and for $\rho_-=0$ and $\rho_+=1$
$$\lim_{t\rightarrow\infty} \PP\left(\frac{h(\tfrac{t}{\gamma},0)-\tfrac{t}{2}}{2^{-1/3}t^{1/3}}\geq -s\right) = F_{\rm{GUE}}(s).$$
\end{theorem}

Thus we see a critical point: For any $\gamma>0$, fluctuations scale like $t^{1/3}$ and have limiting GUE statistics, while for $\gamma=0$, fluctuations scale like $t^{1/4}$ and have limiting Gaussian statistics. Thus, scaling $\gamma$ to zero with the other model parameters one would hope to find a scaling limit which interpolates between these two universality classes. As proved in \cite{ACQ}, that is exactly what the KPZ equation does.

When we later consider directed polymers, this same perspective of looking between universality classes will guide us towards looking at the {\it weak noise} or {\it intermediate disorder} regime -- where again we will encounter the KPZ equation (in the guise of the continuum directed random polymer).

\subsubsection{Six fundamental KPZ universality class geometries and fluctuation statistics}\label{sixkpzclass}

The rescaled initial conditions for the particles/height function plays an important role in the hydrodynamic theory. We turn now to fluctuations and ask how does the initial data affect the scalings and statistics for the long time fluctuations. Height function fluctuations detect differences in the initial data or geometry, though scalings remain governed by the characteristic $1/3$, $2/3$ exponents. There are six types of initial data on which we will focus. Each of these cases corresponds to a markedly different geometry, and the spatial limit process associated to each of these has been rigorously derived through asymptotic analysis of exact formulas. There are, however, also geometries in which the fluctuations of the initial data overwhelm the KPZ class fluctuations (see \cite{CFP1} for example of both Gaussian and KPZ type).


The following six initial conditions have been well-studied and are arguably the most important subclasses of the KPZ universality class (see Figure \ref{sixfig}):
\begin{enumerate}
\item {\bf Wedge (step):} $h_{\gamma}(0,x) = |x|$, or $\eta(0,x) = 0$ for $x\leq 0$ and $1$ for $x>0$.
\item {\bf Brownian (stationary):} $h_{\gamma}(0,x) = RW(x)$ (a simple symmetric random walk), or $\eta(0,x)$ are i.i.d. Bernoulli random variables with parameter $1/2$.
\item {\bf Flat (periodic):} $h_{\gamma}(0,x)$ linearly oscillates between 0 and 1, or $\eta(0,x) = x \mod 2$.
\item {\bf Wedge$\rightarrow$Brownian (half Bernoulli):} $h_{\gamma}(0,x)= -x$ for $x\leq 0$ and $RW(x)$ for $x>0$, or $\eta(0,x)$ are 0 for $x\leq 0$ and are i.i.d. Bernoulli random variables for $x>0$.
\item {\bf Wedge$\rightarrow$Flat (half periodic):} $h_{\gamma}(0,x)= -x$ for $x\leq 0$ and linearly oscillates between 0 and 1 for $x>0$, or $\eta(0,x)$ are 0 for $x\leq 0$ and $x \mod 2$ for $x>0$.
\item {\bf Flat$\rightarrow$Brownian (half periodic / half Bernoulli):}  $h_{\gamma}(0,x)$ linearly oscillates between 0 and 1 for $x\leq 0$ and equals $RW(x)$ for $x>0$, or $\eta(0,x)$ are $x \mod 2$ for $x\leq 0$ and are i.i.d. Bernoulli random variables with parameter $1/2$ for $x>0$.
\end{enumerate}

It is fairly simple to see, heuristically, why these are distinguished. Assume that we are considering a {\it fixed} (though possibly random) initial condition. Furthermore, assume that there is a hydrodynamic limit for the evolution of this initial condition (i.e., a law of large numbers). We are interested in the fluctuations of the height function around its deterministic limit shape. The KPZ scaling suggests that we set $t=\e^{-3/2}T$, $x=\e^{-1}X$ and we scale fluctuations by $\e^{-1/2}$ so as to look at
$$\e^{-1/2}\left(h_{\gamma}(\tfrac{t}{\gamma},x) - \bar{h}(T,0)\right),$$
where we have subtracted $\bar{h}(T,0)$ because the spatial scaling is now smaller than the temporal scaling (and hence only the asymptotic height above the origin plays a role). Since we have fixed the initial data and scaled diffusively, the only possible limits for the rescaled initial data are those which are invariant in distribution under diffusive scaling. The three simplest such limits are: constant zero, Brownian motion, and, infinity (except for $X=0$). As the limit can be different on the left and right of the origin these cases give six different limits, which correspond with the six geometries above. There are other possible limits, such as combinations of independent Brownian motions under the $(\max, +)$ algebra (see \cite{BBP} for instance).

We give two examples to illustrate this rescaling. Consider some finite pattern of particles and holes (for instance particle, particle, hole, hole). Fix the initial condition to be a tiling of this pattern. Under the diffusive scaling of the initial fluctuations this initial condition converges to the same limit (constant zero) as does the flat geometry. Thus, the statistics of the $T>0$ height function fluctuations {\it should} agree. For the second example consider initial conditions where $\eta(0,x)$ are determined by Bernoulli random variables with periodic parameters (i.e., parameter depending on $x\mod m$ for some $m\geq 1$) then the $T=0$ fluctuations converge to a scaling of a Brownian motion and likewise, the statistics of the $T>0$ height function fluctuations should agree with the usual Bernoulli case (Tracy and Widom called these {\it periodic Bernoulli initial conditions} \cite{TWperiodicbernoulli}).

Of course there are some serious mathematical issues which, in fact, are unresolved here. First of all, even if two different initial conditions converge to the same fluctuation initial data, does that mean that their positive $T$ fluctuation statistics coincide as well? The hydrodynamic theory showed this was true at the level of limit shapes, but fluctuations are much more complicated and such a theorem does not exist. Another complication comes from trying to make sense of the diffusive rescaling of initial data which seems to scale to infinity (except at $X=0$). There are many legitimate types of initial data which rescale to this under diffusive scaling.

\hskip2in
\begin{figure}
\setlength{\unitlength}{1pt}
\begin{picture}(200,540)(120,-300)
\linethickness{1pt}
\put(-10,230){\line(1,0){460}}
\put(-10,200){\line(1,0){460}}
\put(-10,170){\line(1,0){460}}
\put(-10,90){\line(1,0){460}}
\put(-10,10){\line(1,0){460}}
\put(-10,-70){\line(1,0){460}}
\put(-10,-150){\line(1,0){460}}
\put(-10,-230){\line(1,0){460}}
\put(-10,-310){\line(1,0){460}}

\put(-10,230){\line(0,-1){540}}
\put(110,200){\line(0,-1){510}}
\put(230,200){\line(0,-1){510}}
\put(345,200){\line(0,-1){510}}
\put(450,230){\line(0,-1){540}}

\put(15,215){\makebox(0,0)[l]{\huge{{\bf Six Geometries of KPZ Universality}}}}

\put(20,185){\makebox(0,0)[l]{\large{{\bf Geometry}}}}
\put(0,160){\makebox(0,0)[l]{$\bullet$ Wedge}}

\put(0,100){\line(1,0){100}}

\put(0,150){\line(1,-1){10}}
\put(10,140){\line(1,-1){10}}
\put(20,130){\line(1,-1){10}}
\put(30,120){\line(1,-1){10}}
\put(40,110){\line(1,-1){10}}

\put(50,100){\line(1,1){10}}
\put(60,110){\line(1,1){10}}
\put(70,120){\line(1,1){10}}
\put(80,130){\line(1,1){10}}
\put(90,140){\line(1,1){10}}

\put(8,99){\makebox(0,)[l]{$\circ$}}
\put(18,99){\makebox(0,)[l]{$\circ$}}
\put(28,99){\makebox(0,)[l]{$\circ$}}
\put(38,99){\makebox(0,)[l]{$\circ$}}
\put(48,99){\makebox(0,)[l]{$\circ$}}

\put(49,94){\makebox(0,)[l]{\tiny{0}}}

\put(59,99){\makebox(0,)[l]{$\bullet$}}
\put(69,99){\makebox(0,)[l]{$\bullet$}}
\put(79,99){\makebox(0,)[l]{$\bullet$}}
\put(89,99){\makebox(0,)[l]{$\bullet$}}
\put(99,99){\makebox(0,)[l]{$\bullet$}}

\put(0,80){\makebox(0,0)[l]{$\bullet$ Brownian}}

\put(0,40){\line(1,0){100}}

\put(0,50){\line(1,1){10}}
\put(10,60){\line(1,-1){10}}
\put(20,50){\line(1,-1){10}}
\put(30,40){\line(1,1){10}}
\put(40,50){\line(1,-1){10}}

\put(50,40){\line(1,-1){10}}
\put(60,30){\line(1,1){10}}
\put(70,40){\line(1,1){10}}
\put(80,50){\line(1,1){10}}
\put(90,60){\line(1,-1){10}}

\put(8,39){\makebox(0,)[l]{$\bullet$}}
\put(18,39){\makebox(0,)[l]{$\circ$}}
\put(28,39){\makebox(0,)[l]{$\circ$}}
\put(38,39){\makebox(0,)[l]{$\bullet$}}
\put(48,39){\makebox(0,)[l]{$\circ$}}

\put(49,34){\makebox(0,)[l]{\tiny{0}}}

\put(59,39){\makebox(0,)[l]{$\circ$}}
\put(69,39){\makebox(0,)[l]{$\bullet$}}
\put(79,39){\makebox(0,)[l]{$\bullet$}}
\put(89,39){\makebox(0,)[l]{$\bullet$}}
\put(99,39){\makebox(0,)[l]{$\circ$}}

\put(0,0){\makebox(0,0)[l]{$\bullet$ Flat}}
\put(0,-40){\line(1,0){100}}

\put(0,-30){\line(1,-1){10}}
\put(10,-40){\line(1,1){10}}
\put(20,-30){\line(1,-1){10}}
\put(30,-40){\line(1,1){10}}
\put(40,-30){\line(1,-1){10}}

\put(50,-40){\line(1,1){10}}
\put(60,-30){\line(1,-1){10}}
\put(70,-40){\line(1,1){10}}
\put(80,-30){\line(1,-1){10}}
\put(90,-40){\line(1,1){10}}

\put(8,-41){\makebox(0,)[l]{$\circ$}}
\put(18,-41){\makebox(0,)[l]{$\bullet$}}
\put(28,-41){\makebox(0,)[l]{$\circ$}}
\put(38,-41){\makebox(0,)[l]{$\bullet$}}
\put(48,-41){\makebox(0,)[l]{$\circ$}}

\put(49,-46){\makebox(0,)[l]{\tiny{0}}}

\put(59,-41){\makebox(0,)[l]{$\bullet$}}
\put(69,-41){\makebox(0,)[l]{$\circ$}}
\put(79,-41){\makebox(0,)[l]{$\bullet$}}
\put(89,-41){\makebox(0,)[l]{$\circ$}}
\put(99,-41){\makebox(0,)[l]{$\bullet$}}

\put(0,-80){\makebox(0,0)[l]{$\bullet$ Wedge$\rightarrow$Brownian}}
\put(0,-140){\line(1,0){100}}

\put(0,-90){\line(1,-1){10}}
\put(10,-100){\line(1,-1){10}}
\put(20,-110){\line(1,-1){10}}
\put(30,-120){\line(1,-1){10}}
\put(40,-130){\line(1,-1){10}}

\put(50,-140){\line(1,1){10}}
\put(60,-130){\line(1,-1){10}}
\put(70,-140){\line(1,-1){10}}
\put(80,-150){\line(1,1){10}}
\put(90,-140){\line(1,1){10}}

\put(8,-141){\makebox(0,)[l]{$\circ$}}
\put(18,-141){\makebox(0,)[l]{$\circ$}}
\put(28,-141){\makebox(0,)[l]{$\circ$}}
\put(38,-141){\makebox(0,)[l]{$\circ$}}
\put(48,-141){\makebox(0,)[l]{$\circ$}}

\put(49,-146){\makebox(0,)[l]{\tiny{0}}}

\put(59,-141){\makebox(0,)[l]{$\bullet$}}
\put(69,-141){\makebox(0,)[l]{$\circ$}}
\put(79,-141){\makebox(0,)[l]{$\circ$}}
\put(89,-141){\makebox(0,)[l]{$\bullet$}}
\put(99,-141){\makebox(0,)[l]{$\bullet$}}


\put(0,-160){\makebox(0,0)[l]{$\bullet$ Wedge$\rightarrow$Flat}}
\put(0,-220){\line(1,0){100}}

\put(0,-170){\line(1,-1){10}}
\put(10,-180){\line(1,-1){10}}
\put(20,-190){\line(1,-1){10}}
\put(30,-200){\line(1,-1){10}}
\put(40,-210){\line(1,-1){10}}

\put(50,-220){\line(1,1){10}}
\put(60,-210){\line(1,-1){10}}
\put(70,-220){\line(1,1){10}}
\put(80,-210){\line(1,-1){10}}
\put(90,-220){\line(1,1){10}}

\put(8,-221){\makebox(0,)[l]{$\circ$}}
\put(18,-221){\makebox(0,)[l]{$\circ$}}
\put(28,-221){\makebox(0,)[l]{$\circ$}}
\put(38,-221){\makebox(0,)[l]{$\circ$}}
\put(48,-221){\makebox(0,)[l]{$\circ$}}

\put(49,-226){\makebox(0,)[l]{\tiny{0}}}

\put(59,-221){\makebox(0,)[l]{$\bullet$}}
\put(69,-221){\makebox(0,)[l]{$\circ$}}
\put(79,-221){\makebox(0,)[l]{$\bullet$}}
\put(89,-221){\makebox(0,)[l]{$\circ$}}
\put(99,-221){\makebox(0,)[l]{$\bullet$}}

\put(0,-240){\makebox(0,0)[l]{$\bullet$ Flat$\rightarrow$Brownian}}
\put(0,-270){\line(1,0){100}}

\put(0,-260){\line(1,-1){10}}
\put(10,-270){\line(1,1){10}}
\put(20,-260){\line(1,-1){10}}
\put(30,-270){\line(1,1){10}}
\put(40,-260){\line(1,-1){10}}

\put(50,-270){\line(1,-1){10}}
\put(60,-280){\line(1,-1){10}}
\put(70,-290){\line(1,1){10}}
\put(80,-280){\line(1,-1){10}}
\put(90,-290){\line(1,1){10}}

\put(8,-271){\makebox(0,)[l]{$\circ$}}
\put(18,-271){\makebox(0,)[l]{$\bullet$}}
\put(28,-271){\makebox(0,)[l]{$\circ$}}
\put(38,-271){\makebox(0,)[l]{$\bullet$}}
\put(48,-271){\makebox(0,)[l]{$\circ$}}

\put(49,-276){\makebox(0,)[l]{\tiny{0}}}

\put(59,-271){\makebox(0,)[l]{$\circ$}}
\put(69,-271){\makebox(0,)[l]{$\circ$}}
\put(79,-271){\makebox(0,)[l]{$\bullet$}}
\put(89,-271){\makebox(0,)[l]{$\circ$}}
\put(99,-271){\makebox(0,)[l]{$\bullet$}}


\put(135,185){\makebox(0,0)[l]{\large{{\bf Limit shape}}}}

\put(120,160){\makebox(0,0)[l]{$\bar{h}(T,TX)=$}}
\put(120,125){\makebox(0,0)[l]{$\begin{cases} -X & X<-1 \\ T\frac{1+X^2}{2} & |X|\leq 1\\ X& X>1 \end{cases}$}}
\put(120,50){\makebox(0,0)[l]{$\bar{h}(T,TX)=T/2$}}
\put(120,-30){\makebox(0,0)[l]{$\bar{h}(T,TX)=T/2$}}
\put(120,-80){\makebox(0,0)[l]{$\bar{h}(T,TX)=$}}
\put(120,-115){\makebox(0,0)[l]{$\begin{cases} -X & X<-1 \\ T\frac{1+X^2}{2} & X\in [-1,0] \\ T/2 & X>0\end{cases}$}}
\put(120,-160){\makebox(0,0)[l]{$\bar{h}(T,TX)=$}}
\put(120,-195){\makebox(0,0)[l]{$\begin{cases} -X & X<-1 \\ T\frac{1+X^2}{2} & X\in [-1,0] \\ T/2 & X>0\end{cases}$}}
\put(120,-270){\makebox(0,0)[l]{$\bar{h}(T,TX)=T/2$}}

\put(250,185){\makebox(0,0)[l]{\large{{\bf Fluctuations}}}}

\put(235,155){\makebox(0,0)[l]{$\bullet$ One pt: $F_{\rm{GUE}}$}}
\put(240,140){\makebox(0,0)[l]{\cite{BDJ,KJ}}}

\put(235,120){\makebox(0,0)[l]{$\bullet$ Multi pt: $\rm{Airy}_2$}}
\put(240,105){\makebox(0,0)[l]{\cite{PS2,KJPNG,BF}}}

\put(235,75){\makebox(0,0)[l]{$\bullet$ One pt: $F_{0}$}}
\put(240,60){\makebox(0,0)[l]{\cite{BR,FS}}}

\put(235,40){\makebox(0,0)[l]{$\bullet$ Multi pt: $\rm{Airy}_{\rm{stat}}$}}
\put(240,25){\makebox(0,0)[l]{\cite{BFP}}}

\put(235,-5){\makebox(0,0)[l]{$\bullet$ One pt: $F_{\rm{GOE}}$}}
\put(240,-20){\makebox(0,0)[l]{\cite{BR2,BR3,FS2,Sa}}}

\put(235,-40){\makebox(0,0)[l]{$\bullet$ Multi pt: $\rm{Airy}_{1}$}}
\put(240,-55){\makebox(0,0)[l]{\cite{BFPS,BFPr}}}

\put(235,-85){\makebox(0,0)[l]{$\bullet$ One pt: $(F_{\rm{GOE}})^2$}}
\put(240,-100){\makebox(0,0)[l]{\cite{BR,BBP,PS1,BAC}}}

\put(235,-120){\makebox(0,0)[l]{$\bullet$ Multi pt: $\rm{Airy}_{2\rightarrow \rm{BM}}$}}
\put(240,-135){\makebox(0,0)[l]{\cite{ImSa,CFP1}}}

\put(235,-185){\makebox(0,0)[l]{$\bullet$ Multi pt: $\rm{Airy}_{2\rightarrow 1}$}}
\put(240,-200){\makebox(0,0)[l]{\cite{BFS}}}

\put(235,-265){\makebox(0,0)[l]{$\bullet$ Multi pt: $\rm{Airy}_{1\rightarrow \rm{BM}}$}}
\put(240,-280){\makebox(0,0)[l]{\cite{BFS2}}}

\put(360,185){\makebox(0,0)[l]{\large{{\bf KPZ/SHE}}}}

\put(350,155){\makebox(0,0)[l]{$\mathcal{Z}(0,X)=\delta_{\{X=0\}}$}}
\put(350,140){\makebox(0,0)[l]{$\bullet$ Converges:\cite{ACQ}}}
\put(350,125){\makebox(0,0)[l]{$\bullet$ One pt:\cite{ACQ,SaSp1}}}
\put(350,110){\makebox(0,0)[l]{(bounds and stats)}}

\put(350,75){\makebox(0,0)[l]{$\mathcal{Z}(0,X)=e^{B(X)}$}}
\put(350,60){\makebox(0,0)[l]{$\bullet$ Converges:\cite{BG}}}
\put(350,45){\makebox(0,0)[l]{$\bullet$ One pt:\cite{BQS,CQ}}}
\put(350,30){\makebox(0,0)[l]{(bounds, NO stats)}}

\put(350,-5){\makebox(0,0)[l]{$\mathcal{Z}(0,X)=1$}}
\put(350,-20){\makebox(0,0)[l]{$\bullet$ Converges:\cite{BG}}}
\put(350,-35){\makebox(0,0)[l]{$\bullet$ One pt: OPEN}}
\put(350,-50){\makebox(0,0)[l]{(NO bounds / stats)}}

\put(350,-80){\makebox(0,0)[l]{$\mathcal{Z}(0,X)=$}}
\put(350,-95){\makebox(0,0)[l]{$e^{B(X)}\mathbf{1}_{X\geq 0}$}}
\put(350,-110){\makebox(0,0)[l]{$\bullet$ Converges:\cite{CQ}}}
\put(350,-125){\makebox(0,0)[l]{$\bullet$ One pt:\cite{CQ}}}
\put(350,-140){\makebox(0,0)[l]{(bounds and stats)}}

\put(350,-165){\makebox(0,0)[l]{$\mathcal{Z}(0,X)=\mathbf{1}_{X\geq 0}$}}
\put(350,-180){\makebox(0,0)[l]{$\bullet$ Converges: here}}
\put(350,-195){\makebox(0,0)[l]{$\bullet$ One pt: OPEN}}
\put(350,-210){\makebox(0,0)[l]{(NO bounds / stats)}}

\put(350,-240){\makebox(0,0)[l]{$\mathcal{Z}(0,X)=$}}
\put(350,-255){\makebox(0,0)[l]{$\mathbf{1}_{X<0} + e^{B(X)}\mathbf{1}_{X\geq 0}$}}
\put(350,-270){\makebox(0,0)[l]{$\bullet$ Converges: here}}
\put(350,-285){\makebox(0,0)[l]{$\bullet$ One pt: OPEN}}
\put(350,-300){\makebox(0,0)[l]{(NO bounds / stats)}}

\end{picture}
\vskip .2in
\caption[Six KPZ universality classes]{Six fundamental geometries in the KPZ universality class. Fluctuations have been entirely classified for these geometries at the level of one-point and multi-point spatial statistics (the temporal evolution of these spatial processes is conjecturally described in \cite{CQ2}). Under suitable weak asymmetry the entire space-time evolution of KPZ class models converge as a process to the Hopf-Cole solution to the KPZ equation (note the logarithmic correction to the scaling in the wedge geometry). Exact formulas and moment/variance bounds for the KPZ equation solutions (stated in terms of the initial data for the associated stochastic heat equation) have only been derived and proved in the wedge, wedge to Brownian, and Brownian (only bounds) cases. Even in those cases, only one-point formulas and bounds exist.}\label{sixfig}
\end{figure}

It turns out that all of these issues can be resolved mathematically if you scale $\gamma$ to zero (like $\gamma =\e^{1/2}$) and deal with the Hopf-Cole transform of the fluctuations, rather than the fluctuations directly. We will return to this point very soon.

Having identified six fundamental growth geometries in the KPZ universality class, we should ask how does the geometry affect the associated fluctuation statistics (the scaling exponents, after all, should all be the same). Thanks to a great deal of effort in the last twelve years, the spatial fluctuations are entirely classified when $\gamma=1$ and exact formulas have been proved to describe the multi point (fixed $T$) statistics. For general $\gamma>0$, Tracy and Widom have extended the one-point statistic results for the wedge and wedge$\rightarrow$Brownian geometries (of course there is no doubt all of this holds for general positive $\gamma$). The chart (\ref{sixfig}) summarizes the six spatial limit processes. As can be predicted by the composite nature of their associated initial conditions, the second set of three limit processes are actually transition processes between combinations of the first three. Understanding the evolution of these spatial limit processes in time $T$ is a problem of great interest which only recently has begun to be answered (see \cite{CQ2}).

\subsubsection{Six fundamental initial data of the KPZ equation}\label{KPZgrowthregime}
We introduced above six fundamental geometries (or subclasses) of the KPZ universality class. Under weak scaling of the asymmetry $\gamma$ (like $\gamma=\beta \e^{1/2}$), we may now rigorously prove that all initial data which rescale diffusively to a limit given by one of these geometries, will have fluctuations which rescale (as $\e$ goes to zero) to the same space-time fluctuation process as in these distinguished and solvable cases.

This fact was proved in the Brownian and flat geometries by Bertini-Giacomin \cite{BG} in the mid 1990s (the following explanation is illustrated in Figure \ref{BGroute}). Rather than directly studying the height function fluctuations, they looked at its Hopf-Cole transform ($f\mapsto \exp\{-f\}$) which they called $Z_\e$. G\"{a}rtner \cite{G} had previously recognized that this transform linearizes the dynamics of the corner growth model into a discrete stochastic heat equation with multiplicative noise. Taking $\e$ to zero, the solution ($Z_\e$) to the discrete stochastic heat equation converges to the solution ($\mathcal{Z}$) to the continuum stochastic heat equation with multiplicative space-time white noise. Un-doing the Hopf-Cole transform gives (at least formally) the KPZ equation. As explained before, due to the ill-posedness of  the KPZ equation, it is better to deal at the level of the well-posed SHE. The initial data associated with the Brownian geometry becomes $\mathcal{Z}_0(X)=\exp\{-B(X)\}$ for a Brownian motion $B(X)$ independent of the white noise, while flat geometry becomes $\mathcal{Z}_0(X)=\exp\{-0\}=1$. Their work also applies to the flat$\rightarrow$Brownian geometry and gives $\mathcal{Z}_0(X)= \mathbf{1}_{X<0} + \exp\{-B(X)\}\mathbf{1}_{X\geq 0}$. Moreover, one does not need to start with the exact initial conditions specified above. It suffices that the initial condition fluctuations converge to the limit point associated with these three geometries.\newline

\begin{figure}
\vskip.2in
\hskip-1.5in
\setlength{\unitlength}{1pt}
\begin{picture}(200,200)(0  ,0)
\linethickness{1pt}
\put(-10,200){\line(1,0){460}}
\put(-10,170){\line(1,0){460}}
\put(-10,90){\line(1,0){460}}
\put(-10,0){\line(1,0){460}}
\put(-10,-60){\line(1,0){460}}
\put(-10,200){\line(0,-1){260}}
\put(220,200){\line(0,-1){260}}
\put(450,200){\line(0,-1){260}}

\put(40,185){\makebox(0,0)[l]{\Large{{\bf Continuum SPDE}}}}
\put(255,185){\makebox(0,0)[l]{\Large{{\bf Discrete model (SPDE)}}}}

\put(0,160){\makebox(0,0)[l]{{\bf Stochastic Heat Equation (SHE)}}}
\put(20,140){\makebox(0,0)[l]{$\partial_T \mathcal{Z} = \tfrac{1}{2} \partial_X^2 \mathcal{Z} - \mathcal{Z}\mathscr{\dot{W}}$}}
\put(40,130){\line(0,-1){35}}
\put(40,95){\line(1,1){10}}
\put(40,95){\line(-1,1){10}}
\put(60,115){\makebox(0,0)[l]{{\it Hopf-Cole solution} to KPZ}}
\put(60,100){\makebox(0,0)[l]{$\mathcal{H}(T,X)=-\log\mathcal{Z}(T,X)$}}

\put(230,160){\makebox(0,0)[l]{{\bf Discrete SHE}}}
\put(250,140){\makebox(0,0)[l]{$\partial_T Z_\e = \tfrac{1}{2}D_\e \Delta_{\e} Z_\e - \mathcal{Z}_\e dM_\e$}}
\put(360,150){\makebox(0,0)[l]{\tiny{(Martingale)}}}
\put(270,130){\line(0,-1){35}}
\put(270,130){\line(1,-1){10}}
\put(270,130){\line(-1,-1){10}}
\put(290,115){\makebox(0,0)[l]{{\it G\"{a}rtner's microscopic Hopf-Cole}}}
\put(290,100){\makebox(0,0)[l]{{\it transform}: $Z_\e = \exp\{-\hfluc\}$}}

\put(240,140){\line(-1,0){115}}
\put(125,140){\line(1,1){10}}
\put(125,140){\line(1,-1){10}}
\put(137,145){\makebox(0,0)[l]{Bertini-Giacomin}}

\put(0,80){\makebox(0,0)[l]{(c) {\bf KPZ Equation}}}
\put(20,60){\makebox(0,0)[l]{`` $\partial_T \mathcal{H} = \tfrac{1}{2} \partial_X^2 \mathcal{H} -  \tfrac{1}{2} \left(\partial_X \mathcal{H}\right)^2 +\mathscr{\dot{W}}$ ''}}
\put(20,45){\makebox(0,0)[l]{`` $\E[\mathscr{\dot{W}}(T,X)\mathscr{\dot{W}}(S,Y)]=\delta_{T=S}\delta_{X=Y}$ ''}}
\put(40,35){\line(0,-1){30}}
\put(40,5){\line(1,1){10}}
\put(40,5){\line(-1,1){10}}
\put(60,15){\makebox(0,0)[l]{$\mathcal{U}=\partial_X  \mathcal{H}$}}

\put(0,-10){\makebox(0,0)[l]{{\bf Stochastic Burgers Equation}}}
\put(20,-30){\makebox(0,0)[l]{`` $\partial_T \mathcal{U} = -\tfrac{1}{2} \partial_X \mathcal{U}^2 +  \tfrac{1}{2} \partial_X^2 \mathcal{U} +\partial_X\mathscr{\dot{W}}$ ''}}

\put(230,80){\makebox(0,0)[l]{{\bf Growth/Height Interface Models}}}
\put(250,60){\makebox(0,0)[l]{Corner growth model $h(t,x)$}}
\put(250,45){\makebox(0,0)[l]{(given by a discrete space SPDE)}}
\put(270,40){\line(0,-1){35}}
\put(270,40){\line(1,-1){10}}
\put(270,40){\line(-1,-1){10}}
\put(290,25){\makebox(0,0)[l]{Define $h$ such that}}
\put(290,10){\makebox(0,0)[l]{$h(t,x)-h(t,x-1)=\hat{\eta}(t,x)$}}

\put(230,-10){\makebox(0,0)[l]{{\bf Interacting Particle Systems}}}
\put(250,-30){\makebox(0,0)[l]{Simple exclusion process: $\hat{\eta}(t,x)=\pm1$}}
\put(250,-45){\line(1,0){180}}
\put(250,-45){\line(1,1){10}}
\put(250,-45){\line(1,-1){10}}
\put(430,-45){\line(-1,1){10}}
\put(430,-45){\line(-1,-1){10}}
\put(260,-45){\makebox(0,0)[l]{$\vdash$}}
\put(270,-45){\makebox(0,0)[l]{$\vdash$}}
\put(280,-45){\makebox(0,0)[l]{$\vdash$}}
\put(290,-45){\makebox(0,0)[l]{$\vdash$}}
\put(300,-45){\makebox(0,0)[l]{$\vdash$}}
\put(310,-45){\makebox(0,0)[l]{$\vdash$}}
\put(320,-45){\makebox(0,0)[l]{$\vdash$}}
\put(330,-45){\makebox(0,0)[l]{$\vdash$}}
\put(340,-45){\makebox(0,0)[l]{$\vdash$}}
\put(350,-45){\makebox(0,0)[l]{$\vdash$}}
\put(360,-45){\makebox(0,0)[l]{$\vdash$}}
\put(370,-45){\makebox(0,0)[l]{$\vdash$}}
\put(380,-45){\makebox(0,0)[l]{$\vdash$}}
\put(390,-45){\makebox(0,0)[l]{$\vdash$}}
\put(400,-45){\makebox(0,0)[l]{$\vdash$}}
\put(410,-45){\makebox(0,0)[l]{$\vdash$}}
\put(420,-45){\makebox(0,0)[l]{$\vdash$}}

\put(278,-45){\makebox(0,0)[l]{$\bullet$}}
\put(308,-45){\makebox(0,0)[l]{$\bullet$}}
\put(318,-45){\makebox(0,0)[l]{$\bullet$}}
\put(358,-45){\makebox(0,0)[l]{$\bullet$}}
\put(408,-45){\makebox(0,0)[l]{$\bullet$}}

\put(0,0){\line(1,0){140}}

\end{picture}\newline\newline\newline\newline\newline\newline
\caption[Discrete to continuum]{Illustration of the route between discrete and continuum models. The quotation marks imply that these are ill-posed SPDEs and interpreted via the solution to the stochastic heat equation.}\label{BGroute}
\end{figure}

Bertini-Giacomin's results do not extend to the three wedge-related geometries. It is clear from the above discussion what initial data should be associated with the two transition geometries -- wedge$\rightarrow$Brownian and wedge$\rightarrow$flat. The first should have $\mathcal{Z}_{0}(X) = \exp\{-B(X)\}\mathbf{1}_{X\geq 0}$ while the second should have $\mathcal{Z}_{0}(X) = \mathbf{1}_{X\geq 0}$. Of course, it does not make sense to un-do the Hopf-Cole transform for $T=0$, but due to a theorem of M\"{u}ller \cite{Muller}, $\mathcal{Z}$ is almost surely everywhere positive for $T>0$ and $X\in R$. The above convergence is made rigorous in \cite{ACQ,CQ} (and builds on Bertini-Giacomin) and will be explained in Section \ref{Gartnersec}.

For the wedge geometry the above arguments fail. Since $h_{\gamma}(0,x)=|\e^{-1}X|$, it is clear that $Z_e = \exp\{-\e^{1/2}|\e^{-1}X|\} = \exp\{-\e^{-1/2}|X|\}$. As $\e$ goes to zero, $Z_\e$ goes to 0 for all $X\neq 0$, and goes to $1$ for $X=0$ (this is analogous to the issue discussed earlier regarding initial data rescaling diffusively to infinity). Such initial data is massless and can not be right. The solution, however, is simple and was discovered in \cite{ACQ}. Bertini-Giacomin's scalings do not exactly apply for this geometry -- there is a logarithmic correction needed to the height function fluctuations when dealing with weak asymmetry. This can be seen by observing that putting a prefactor of $\e^{-1/2}/2$ in front of the above $Z_\e$, results in conservation of mass and in this modified scaling, the wedge geometry scales to a delta function at zero. Thus, with this additional correction, \cite{ACQ} shows that the wedge geometry corresponds to $\mathcal{Z}_0(X)= \delta_{X=0}$ (in fact one can show that any wedge-like initial condition has the same limit). It is exactly this $\e^{-1/2}/2$ prefactor which accounts for the necessary correction to Bertini-Giacomin's scalings.

This shows that, at least in the weakly asymmetric scaling, these six geometries correspond to large universality subclasses. The asymmetry (written before as $\gamma =\beta\e$) can be ratcheted up after taking $\e$ to zero, by increasing $\beta$. Due to scalings of the continuum equations, this amounts to taking $T$ to infinity. As such, one would hope to recover the $\gamma>0$ statistics and scales by simply studying the statistics of the long-time limit of the KPZ equation.

A proof of this fact (at the level of one point statistics) has now been given by \cite{ACQ} and \cite{CQ} for the wedge and the wedge$\rightarrow$Brownian initial data / geometry. This relies on taking the long time limit of the rigorously proved explicit formula for the finite $T$ solution to the KPZ equation (with these initial data).

Using the replica trick method and carefully summing divergent series \cite{Dot,CDR} have successfully (though highly non-rigorously) re-derived the above formula and large $T$ asymptotic result for the wedge geometry. Using this method, \cite{ProS2} gave a non-rigorous derivation of the large $T$ limit of the multi point statistics which matched up with the KPZ class statistics (the $\textrm{Airy}_2$ process) for the wedge geometry.

At a rigorous level, nothing is known about the multi-point (fixed time) statistics for the solution to the KPZ equation. Aside from the exact solutions of \cite{ACQ} and \cite{CQ}, the only other case with some success is that of the Brownian geometry -- \cite{BQS} and \cite{CQ} have derived bounds showing that the variance and moments of the solution to the KPZ equation with Brownian initial data are of the correct scale (\cite{CQ} also provides some tail decay bounds in one direction). The exact statistics here, and bounds or statistics for the other three cases of initial data remain entirely open.

\subsection{Discrete and continuum directed polymers in random media}
Until now we have dealt entirely with growth processes and interacting particle systems. Under suitable scaling we showed (at least for the simple exclusion process or corner growth model) that the KPZ equation arises as a continuum scaling limit. We also saw that the initial conditions for the discrete models translate readily into initial data for the KPZ equation. In fact, it was only at the level of the SHE where everything became rigorous and meaningful (especially with respect to the wedge geometry).

In this section we will observe that, via a version of the Feynman Kac formula, it is possible to interpret the solution to the SHE as the partition function for what is called the {\it continuum directed random polymer} (CDRP). Thus the free energy (essentially the logarithm of the partition function) for the CDRP corresponds with the Hopf-Cole solution to the KPZ equation. The initial data for the SHE corresponds to an initial potential which affects the starting position of the polymer. The fundamental solution for which we have presented the exact probability distribution formula above corresponds with fixing the departure position at zero.

Once we have introduced the CDRP we will provide a few approximation schemes which have the CDRP partition function (hence SHE) as a limit. Most significant of these schemes is that of approximation via a large class of discrete directed polymers in random media (explained in Section \ref{AKQDDPSEC}).

Polymer models are important from a number of perspectives. They were introduced in \cite{HuHe} to study the domain walls of Ising models with impurities \cite{LFC} and have been applied to other problems like vortices in superconductors \cite{BFG}, roughness of crack interfaces \cite{HHR}, Burgers turbulence \cite{FNS}, interfaces in competing bacterial colonies \cite{HHRN} (see also the physical review \cite{HHZ} or \cite{FH} for more applications).

Directed polymers in disordered environments provide a unified mathematical framework for studying a variety of different abstract and physical problems. They represent generalizations of path integrals through disordered potentials and are thus useful ways of solving stochastic (partial) differential equations; they are useful for solving and formalizing many optimization problems, including important problems in bio-statistics \cite{SM, HL, MMN,SAY} and operations research \cite{BBSSS}; they arise in the study of branching Brownian motions and random walks in random environments \cite{Bra}; they serve as paradigms for the study of other disordered systems \cite{DerSpohn}; and phenomena like pinning and wetting \cite{AV,Holl}. On top of that, their relationship to random growth models and interacting particle systems means that any exact solvability developed in the framework of polymers, can be translated (at least at the continuum level) into exact solvability of the associated continuum growth models or particle systems (SPDEs). Add to these topics the connections to combinatorics and integrable systems (see \cite{KJ,COSZ} for example) and it is clear why polymers have received so much interest in recent years.

\subsubsection{The free energy of the continuum directed random polymer}\label{freeECDRP}

Let us start at the level of physics and write down the following expression for the free energy of the CDRP:
\begin{equation}\label{CDRPFreeEnergy}
\mathcal{F}(T,X) = \log \E \left[ :\:\!\exp\!: \left\{-\int_{0}^{T} \dot{\mathscr{W}}(t,b(t))dt\right\}\right]
\end{equation}
where the expectation $\E$ is over Brownian bridges $b(\cdot)$ such that $b(0)=0$ and $b(T)=X$, and $:\!\exp\!:$ is known as the {\it Wick exponential}. Observe that $\mathcal{F}(T,X)$ is random with respect to the disorder (the Gaussian space-time white noise $\dot{\mathscr{W}}$). If one considers the white noise as a random potential in which paths arrange themselves according to a Boltzmann weight\footnote{As written the inverse-temperature parameter $\beta$ has been set to equal 1, though rescaling time like $\beta^4 T$ and space like $\beta^2 X$ effectively reintroduces this factor}, then $\exp\{\mathcal{F}\}$ is the partition function for this ensemble and $\mathcal{F}$ is the free energy.

Mathematically, equation (\ref{CDRPFreeEnergy}) requires some work to make rigorous and in Section \ref{approxSec} we give the following five schemes to define / approximate the free energy of the CDRP -- all of which are equivalent in that they can be shown to lead to the same object. They are:

\begin{enumerate}
\item Chaos series and time ordering,
\item Stochastic PDEs,
\item Spatial smoothing of the white noise,
\item Discrete directed polymers,
\item Growth processes and interacting particle systems.
\end{enumerate}

The fourth and fifth scheme are perhaps most interesting from a physical perspective. That growth processes and interacting particles systems would be related to the CDRP is initially unclear. However, one may recognize that the path integral in the definition of the free energy can be interpreted, via a version of the Feynman Kac formula, as the solution to the SHE with multiplicative white noise. Thus the free energy is closely related to the KPZ equation, which we have seen is a continuum growth model.

The connection between polymers and growth models was previously observed only at the level of zero-temperature polymers (last passage percolation) and totally asymmetric growth models (the corner growth model with $\gamma=1$).

It is worth noting that while everything so far has been focused on the case of the point-to-point polymer, it is possible to likewise consider other initial potentials in which the polymer optimizes its starting position. This corresponds with different initial data for the SHE, to different growth geometries for the growth models and to different initial conditions for the particle systems (see Section \ref{beyondptpt}).

\subsubsection{Exact solvability of the discrete and continuum directed polymer}

Having established above that the CDRP and its free energy / partition function is a truly universal object, it is of great interest to compute scaling exponents, exact and asymptotic statistics, and the answers to various other probabilistic questions related to this random processes. This direction has seen an explosion of progress in the past three years. There are essentially three approaches which have been recorded so far, and each one seems to uncover new structure and expand the scope of what can be computed. We list them chronologically in order of discovery. All but the replica trick provide rigorous means to access information about the CDRP (though the replica trick is useful none-the-less for deriving new formulas). These approaches are explained in Section \ref{solvablepolymers}.

\begin{enumerate}
\item The Tracy Widom ASEP (or equivalently corner growth model) formulas,
\item The replica trick and Bethe ansatz,
\item The solvable finite temperature polymers.
\end{enumerate}

The first approach greatly extends work of \cite{Sch}, while the third approach likewise extends the solvability of the zero temperature polymer models (last passage percolation) via the Robinson-Schensted-Knuth correspondence (studied at length since the work of \cite{BDJ,KJ} in the late 1990s) to the positive temperature setting. The replica trick has been used to study the CDRP since the work of Kardar \cite{K} but only recently developed to the point of yielding formulas for statistics.

\subsection{Some open problems}\label{opensec}
In this Section we briefly highlight a few important open questions (and any partial progress towards resolving them).

\begin{enumerate}
\item {\bf The full scope of solvability for the KPZ equation:}
Figure \ref{sixfig} highlighted the fact that for the KPZ equation we presently do not have one-point distribution results for four of the six fundamental initial data for which we expect solvability. Moreover, we do not have any rigorous results for spatial processes for the KPZ equation at finite time and it is not clear whether the formula derived in \cite{ProS1,ProS2} via replica trick is correct (as they note, it is under a false factorization assumption that they derive the finite $T$ formula and it is only in the large $T$ limit that the factorization seems to be correct and they recover the Airy$_2$ process).

In a related question one would like to understand Tracy and Widom's exact formulas through a more systematic and algebraic perspective in which their combinatorial identities have meaning. Such an understanding could enable extensions of their approach to more general initial data and multi-point distributions.

\item {\bf The multi-layer free energy process:}
The image of exponential weight matrices under the RSK correspondence can be projected to a measure on partitions where the lengths of the partition equal (in distribution) the eigenvalues of the Lageurre Unitary Ensemble of random matrix theory. It is only the top level of the partition (or eigenvalue) which we have focused on above. However, the entire collection of top eigenvalues rescales to a non-trivial point process (which can further be extended into a multi-layer extension of the Airy$_2$ process.

It was observed in \cite{OCon} (and then further discussed in \cite{OConWarren} that there exists a similar multi-layer process extension to the free energy of finite temperature polymers -- and in particular the CDRP. It is of interest to study the finite $T$ statistics of this multi-layer process and derive finite $T$ versions of many of the processes associated with classical random matrix theory.

\item {\bf Path properties:}
The entire multi-layer free energy process, as well as the multi-layer Airy$_2$ scaling limit, should have certain nice path properties. The first property is that every path should be locally absolutely continuous with respect to Brownian motion. This have been shown in \cite{CH} for the multi-layer Airy$_2$. Along the way to proving that result \cite{CH} proves that the entire ensemble of lines in the multi-layer Airy$_2$ process has a Brownian Gibbs resampling property -- which is to say that removing a portion of the $k^{th}$ line between $a$ and $b$, one can resample according to a Brownian bridge conditioned to start and end at the correct locations and to avoid the $(k+1)^{st}$ and $(k-1)^{st}$ line. The result of this resampling is a new configuration of lines which is distributionally equivalent to the original configuration. Such a property should hold for the finite temperature multi-layer free energy process.

\item {\bf Large deviations and the Painlev\'{e} II like expression of \cite{ACQ}:}
The tails of the Tracy and Widom $F_{{\rm GUE}}(s)$ \cite{TW} distribution decay differently --- the upper tail decays like $e^{-cs^{\tfrac{3}{2}}}$ while the lower tail decays like $e^{-cs^3}$. The lower tail exponent was particularly difficult to derive and required the connection with Painlev\'{e} II established in \cite{TW}. One would like to compute similar tail estimates for the finite $T$ solution to the KPZ equation as in \cite{ACQ}. In \cite{ACQ} an integro-differential equation generalizing the Painlev\'{e} II equation was derived. However, it has not yet been analyzed to the point of being able to prove tail estimates of the sort described above, nor has the connection with integrable systems been sufficiently understood (for instance whether it arises in relation to inverse scattering as in \cite{DT}). In \cite{CQ} some upper tail estimates are derived for various KPZ initial conditions.

\item {\bf Universality:}
The universality conjecture roughly says that (within reason) changing local rules of a model will not affect the global or scaling limit behavior. For instance, for random matrix theory changing the distribution of random variable matrix entries or the type of eigenvalue potential does not affect the limiting statistics for eigenvalues. Herein we have considered a number of types of systems including growth models, particle systems and polymers. We are particularly interested in universality of the KPZ equation (hence studying these systems under weak asymmetry). For polymers this goal is in progress in the work of \cite{AKQ}. For growth models and particle systems the only universality results are those of \cite{BG,ACQ,CQ} which deal only with the corner growth model / simple exclusion process.

It is of great interest to prove that for more general particle systems the KPZ equation still arises. For instance one could consider the exclusion process with jumps beyond just nearest neighbor. Or also jump rates which depend on the local environment. Work of \cite{JaraGon} deals with these more general processes and shows tightness of the fluctuation field within a class of so called ``energy solutions'' to the KPZ equation. This class certainly contains the unique Hopf-Cole solution (by virtue of the work of \cite{BG,ACQ,CQ}) but it is far from clear whether it is the only element of whether there are (possibly infinitely) many other energy solutions.

Another form of universality is with respect to the scaling exponents. A result which might reach beyond the limited class of exactly solvable models is \cite{BalKomSep} were exponents are derived for a class of zero range processes that are defined by the requirement that the slode of the nondecreasing, concave jump rate function decrease geometrically.

\item {\bf Other questions:} 
Without commenting on them, it is worth noting three other important questions. The first deals with proving the conjectural description and formulas of \cite{CQ2} with regards to the renormalization fixed point for the KPZ universality class. The second deals with making sense of the KPZ equation and CDRP in more than one spatial dimension (see \cite{IS,CSY} for some polymer results in higher dimensions). And the third deals with the still mysterious connections between the classes of models discussed above and random matrix theory.
\end{enumerate}

\subsection{Outline and acknowledgements}
The introductory section above contained the statement of a number of important results, a discussion of the KPZ universality class and subclasses, and an overview of the material which follows. In Section \ref{talk1} we present the rigorous connection between the WASEP and the KPZ equation and explain step by step how one goes about proving such a connection. In Section \ref{talk2} we show how, by using that connection and the explicit work of Tracy and Widom, one can derive the exact statistics for the KPZ equation. We also explain some of the important and highly non-trivial technical issues one must overcome to make such a derivation rigorous. In Section \ref{talk3} we review the theory of directed polymers in random media and show how the continuum directed random polymer (CDRP) is a universal scaling limit for a wide class of such models. We then review developments in the solvability of random polymers and the KPZ equation.

I wish to thank MSRI and the organizers of the Random Matrix Theory workshop during the fall of 2010. The notes on which this survey is based were developed during that semester and are loosely based on a MSRI/Evans lecture I had the opportunity to deliver in September 2010, and a series of lectures (organized by Freydoun Rezakhanlou) in the Berkeley math department I delivered during the fall of 2010. I would also like to thank all those who attended these lectures and encouraged me to turn these into a survey article. Parts of this survey were also written during my time at MSRI, the Fields Institute, and IMPA.

I would also like to thank many people for their collaboration, support or illuminating discussions (and often all three): Mark Adler, G\'{e}rard Ben Arous, Jinho Baik, Pasquale Calabrese, Percy Deift, Amir Dembo, Victor Dotsenko, Steven Evans, Patrik Ferrari, Timothy Halpin-Healy, Charles Newman, Neil O'Connell, Jeremy Quastel, Freydoun Rezakhanlou, Tomohiro Sasamoto, Timo Sepp\"{a}l\"{a}inen, Vladas Sidoravicius, Herbert Spohn and Pierre van Moerbeke.

Ongoing support from the NSF through the PIRE grant OISE-07-30136 is acknowledged and appreciated.

\section{Weakly asymmetric simple exclusion process approximation of the KPZ equation}\label{talk1}
In the original 1986 paper of Kardar, Parisi and Zhang \cite{KPZ}, it is predicted that the eponymous (KPZ) equation captures the behavior of the fluctuations of a wide class of discrete models. This prediction was rather vague, however, and did not explain exactly in what sense the continuum equation should be related to the discrete models. Still, it led to a prediction for the long time fluctuation scaling exponents for many physical systems and mathematical models -- which have since been shown to be startlingly accurate in certain cases. Even though the scaling exponents are universal through a wide class of models, \cite{KPZ} made no predictions about the exact statistics. Subsequent works in mathematics have revealed an additional level of universality -- that the statistics seem to only depend on certain model characteristics such as the growth geometry, and are otherwise universal across different models. There are now also reasonably well developed conjectures about large classes of models which are in the KPZ universality class (see for instance the review of \cite{KK}).

The purpose of this section is to explain how to make rigorous sense of the relationship between discrete models and the KPZ equation. We will focus mainly on a single discrete growth model called the {\it corner growth model} and its connection, under a particular type of weakly asymmetric scaling, to the KPZ equation. The below diagram records three levels of discrete and continuous processes and their relationship which is rigorously proved at the top level.
The relationship between the discrete corner growth and continuum KPZ equation exists only under a specific weakly asymmetric scaling for the corner growth model. It was the work of Bertini-Giacomin \cite{BG} that first rigorously proved this fact (they focused entirely on the model near stationary initial data). They considered this scaling in terms of hydrodynamics and interpreted the KPZ equation as the fluctuations one sees upon looking past the hydrodynamic scaling. In Section \ref{critscalingobject} we saw how one can also interpret this weak asymmetry as scaling into the critical point between two universality classes (the KPZ class and the EW or Edwards-Wilkinson class). The KPZ class arises when there is positive asymmetry and the EW class when their is zero asymmetry. The recent success in computing the statistics of the KPZ equation shows indeed that the KPZ equation really represents a crossover in scaling and statistics between these two extremes.

We now set out and briefly explain the four tasks of rigorously relating the corner growth model to the KPZ equation which will be the focus of this section (see Figure \ref{BGroute} for an illustration of these steps):
\begin{enumerate}
\item {\bf Define the corner growth model SPDE:} We provide an SPDE for the corner growth model and show how it arises from considering an interactive particle system called the simple exclusion process (SEP). We will also say a bit about the hydrodynamic theory (LLN) for this system.
\item {\bf Make rigorous sense of the KPZ equation:} We will address the question of how to define what it means to solve the ill-posed KPZ equation . We will provide a definition for the {\it Hopf-Cole solution to KPZ} which is defined as $\mathcal{H} = -\log\mathcal{Z}$ where $\mathcal{Z}$ solve the stochastic heat equation which is well-posed. Differentiating the KPZ equation yields the stochastic Burgers equation.
\item {\bf G\"{a}rtner's discrete Hopf-Cole transform:} We will show that, miraculously, applying the Hopf-Cole transform to the corner growth model results in a discrete version of the SHE.
\item {\bf Convergence of discrete SHE to continuum SHE under weak asymmetry:} We will follow the work of Bertini-Giacomin, and the extension of \cite{ACQ}, \cite{CQ}, to show how we can prove space-time process convergence of the discrete to continuum solution of the SHE.
\end{enumerate}

\subsection{Corner growth model stochastic PDE}
In the introductory section we informally defined the simple exclusion process (SEP), and then by integration, the corner growth model. Presently we will show how to rigorously construct the dynamics of the SEP, and in doing so we will write down a stochastic PDE for the corner growth model height function. The SEP is a continuous time, discrete space Markov process. Its state space can be written in terms of occupation variables $\eta_t\in \{0,1\}^{\Z}$ or spin variables $\hat{\eta}_t\in \{-1,1\}^{\Z}$ (related by $2\eta_t(x)-1 = \hat{\eta}_t(x)$):
\begin{equation}
\eta(t,x) = \begin{cases} 1 & \textrm{particle}\\ 0 & \textrm{hole}\end{cases} \qquad\qquad \hat{\eta}(t,x) = \begin{cases} 1 & \textrm{particle}\\ -1 & \textrm{hole}\end{cases}
\end{equation}

\subsubsection{The graphical construction and the jump processes $L$ and $R$}
One can define the process in terms of its Markov generator. However, we will not take this tact here, and will rather focus on a {\it graphical construction} of the process (see \cite{Lig} for more on this).

\setlength{\unitlength}{2pt}
\begin{figure}
\begin{picture}(200,100)(0,-20)
\put(0,0){\line(1,0){200}}

\put(75,-5){\makebox(0,0)[l]{$x$}}

\put(01,37){\makebox(0,0)[l]{$\longrightarrow$}}
\put(01,57){\makebox(0,0)[l]{$\longrightarrow$}}
\put(01,66){\makebox(0,0)[l]{$\longleftarrow$}}

\put(10,0){\line(0,1){70}}
\put(11,20){\makebox(0,0)[l]{$\longrightarrow$}}
\put(11,45){\makebox(0,0)[l]{$\longleftarrow$}}
\put(11,14){\makebox(0,0)[l]{$\longleftarrow$}}
\put(11,67){\makebox(0,0)[l]{$\longleftarrow$}}
\put(20,0){\line(0,1){70}}

\put(20,0){\line(0,1){70}}
\put(30,0){\line(0,1){70}}

\put(30,0){\line(0,1){70}}
\put(31,4){\makebox(0,0)[l]{$\longrightarrow$}}
\put(31,24){\makebox(0,0)[l]{$\longrightarrow$}}
\put(31,12){\makebox(0,0)[l]{$\longleftarrow$}}
\put(31,19){\makebox(0,0)[l]{$\longleftarrow$}}
\put(40,0){\line(0,1){70}}

\put(40,0){\line(0,1){70}}
\put(41,60){\makebox(0,0)[l]{$\longrightarrow$}}
\put(41,53){\makebox(0,0)[l]{$\longleftarrow$}}
\put(50,0){\line(0,1){70}}

\put(50,0){\line(0,1){70}}
\put(51,34){\makebox(0,0)[l]{$\longleftarrow$}}
\put(60,0){\line(0,1){70}}

\put(60,0){\line(0,1){70}}
\put(61,23){\makebox(0,0)[l]{$\longrightarrow$}}
\put(61,9){\makebox(0,0)[l]{$\longleftarrow$}}
\put(61,53){\makebox(0,0)[l]{$\longleftarrow$}}
\put(70,0){\line(0,1){70}}

\put(70,0){\line(0,1){70}}
\put(71,3){\makebox(0,0)[l]{$\longrightarrow$}}
\put(71,13){\makebox(0,0)[l]{$\longrightarrow$}}
\put(71,36){\makebox(0,0)[l]{$\longleftarrow$}}
\put(80,0){\line(0,1){70}}

\put(80,0){\line(0,1){70}}
\put(81,52){\makebox(0,0)[l]{$\longrightarrow$}}
\put(81,4){\makebox(0,0)[l]{$\longleftarrow$}}
\put(81,14){\makebox(0,0)[l]{$\longleftarrow$}}
\put(81,57){\makebox(0,0)[l]{$\longleftarrow$}}
\put(90,0){\line(0,1){70}}

\put(90,0){\line(0,1){70}}
\put(91,17){\makebox(0,0)[l]{$\longleftarrow$}}
\put(91,41){\makebox(0,0)[l]{$\longleftarrow$}}
\put(100,0){\line(0,1){70}}

\put(100,0){\line(0,1){70}}
\put(101,12){\makebox(0,0)[l]{$\longrightarrow$}}
\put(101,51){\makebox(0,0)[l]{$\longrightarrow$}}
\put(101,62){\makebox(0,0)[l]{$\longleftarrow$}}
\put(110,0){\line(0,1){70}}

\put(110,0){\line(0,1){70}}
\put(111,53){\makebox(0,0)[l]{$\longrightarrow$}}
\put(111,43){\makebox(0,0)[l]{$\longleftarrow$}}
\put(120,0){\line(0,1){70}}

\put(120,0){\line(0,1){70}}
\put(121,59){\makebox(0,0)[l]{$\longrightarrow$}}
\put(121,22){\makebox(0,0)[l]{$\longleftarrow$}}
\put(121,4){\makebox(0,0)[l]{$\longleftarrow$}}
\put(130,0){\line(0,1){70}}

\put(130,0){\line(0,1){70}}
\put(131,7){\makebox(0,0)[l]{$\longrightarrow$}}
\put(131,14){\makebox(0,0)[l]{$\longrightarrow$}}
\put(131,35){\makebox(0,0)[l]{$\longrightarrow$}}
\put(131,3){\makebox(0,0)[l]{$\longleftarrow$}}
\put(131,47){\makebox(0,0)[l]{$\longleftarrow$}}
\put(140,0){\line(0,1){70}}

\put(140,0){\line(0,1){70}}
\put(141,22){\makebox(0,0)[l]{$\longrightarrow$}}
\put(150,0){\line(0,1){70}}

\put(150,0){\line(0,1){70}}
\put(151,10){\makebox(0,0)[l]{$\longrightarrow$}}
\put(151,23){\makebox(0,0)[l]{$\longleftarrow$}}
\put(151,29){\makebox(0,0)[l]{$\longleftarrow$}}
\put(160,0){\line(0,1){70}}

\put(160,0){\line(0,1){70}}
\put(161,56){\makebox(0,0)[l]{$\longrightarrow$}}
\put(161,49){\makebox(0,0)[l]{$\longleftarrow$}}
\put(170,0){\line(0,1){70}}

\put(170,0){\line(0,1){70}}
\put(180,0){\line(0,1){70}}

\put(180,0){\line(0,1){70}}
\put(181,50){\makebox(0,0)[l]{$\longrightarrow$}}
\put(181,63){\makebox(0,0)[l]{$\longrightarrow$}}
\put(181,43){\makebox(0,0)[l]{$\longleftarrow$}}
\put(181,16){\makebox(0,0)[l]{$\longleftarrow$}}
\put(181,38){\makebox(0,0)[l]{$\longleftarrow$}}
\put(190,0){\line(0,1){70}}

\put(191,11){\makebox(0,0)[l]{$\longrightarrow$}}
\put(191,49){\makebox(0,0)[l]{$\longleftarrow$}}
\end{picture}
\vskip -.2in
\caption[Graphical constriction of SEP]{The graphical construction of the SEP (the horizontal axis represents space and the vertical axis represents time). Arrows are distributed independently in each column according to Poisson point processes with rates $q$ (left pointing arrows) and rates $p$ (right pointing arrows).}\label{graphicalconstruction}
\end{figure}
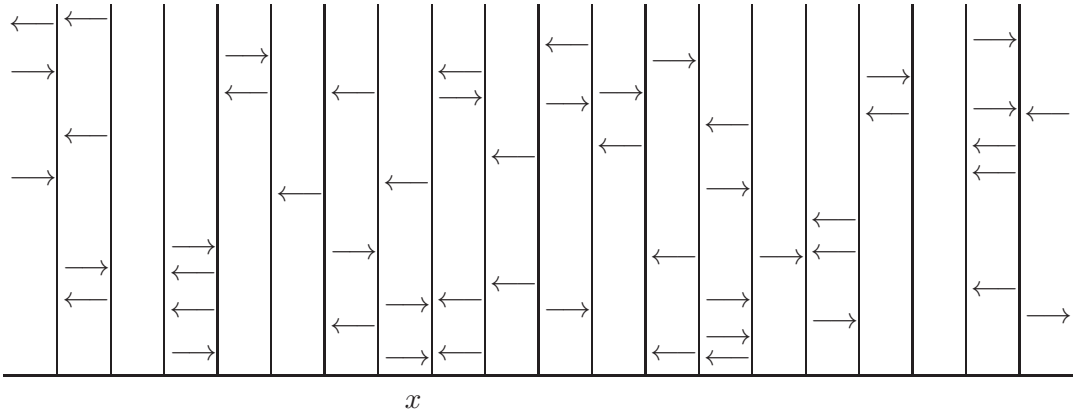

The figure above illustrates this construction. The horizontal axis represents space and the vertical axis represents time. Each vertical line will be called a ladder, and each region between two consecutive ladders will be called a column. We will construct the SEP up to a fixed time $T$. Consider two Poisson point processes $L$ and $R$. $L$ is a random point measure on $\Z\times [0,T]$ with the property that for each $x\in \Z$, $L(x)$ are i.i.d. poisson point processes (on $[0,T]$) with intensity $q$. In the above figure we represent the points of $L$ as left arrows in each column ($x$ identifies a column). Likewise, $R$ is defined but with intensity $p$, and its points are represented as right arrows. These point processes will be the randomness driving our SPDE. It is common to treat $L(x)$ and $R(x)$ as recording the net number of jumps above $x$ whereas $\partial_t L(x)$ and $\partial_t R(x)$ represent the random measure composed of adding delta functions for each jump.

Having fixed an environment, we now introduce particles at time 0 (at most one per ladder) and let them evolve as follows: particles stay in their ladder until an arrow at which point they attempt to follow the direction of the arrow. If their destination site is unoccupied the jump is achieved, otherwise the arrow is ignored. One should pause to ask whether such a dynamic is well-defined on all of $\Z$. Without going into the fine details, the saving grace here is the existence of gaps. In the above figure there are two columns without any arrows. Therefore, up to time $T$, the evolution of the SEP between the gaps is independent of everything else -- thus reducing our considerations to that of a finite state Markov chain. For any fixed $T$ such gaps will almost surely exist at a finite distance from the origin.

\subsubsection{Coupling with the corner growth model}
As explained in the introduction, the corner growth model can be viewed as the integrated spin variables of the simple exclusion process for a given asymmetry $\gamma = q-p$. Fixing that $h_{\gamma}(0,0)=0$ this coupling is given by setting, as in equation (\ref{heightfunctiondef})
\begin{equation}
h_{\gamma}(t,x)= \begin{cases} 2N(t)+\sum_{0<y\leq x} \hat{\eta}(t,y), & x>0,\\
2N(t), &x=0,\\
2N(t)-\sum_{0<y\leq x} \hat{\eta}(t,y), & x<0,
\end{cases}
\end{equation}
where $N(t)$ records the net number of particles to cross from site $1$ to site $0$ in time $t$ and where $\hat{\eta}(t,x)$ equals $1$ if there is a particle at $x$ at time $t$ and $-1$ otherwise.

The dynamics of the exclusion process are readily translated into corner growth model dynamics. In terms of spin variables:
\begin{itemize}
\item Local valley $\Rightarrow$ local hill at rate $q$: When $\hat{\eta}(t,x)=-1$ and $\hat{\eta}(t,x+1)=1$, then $h_{\gamma}(t,x)$ will increase by 2 at rate $qdt$ (and hence leaving $\hat{\eta}(t+dt,x)=1$ and $\hat{\eta}(t+dt,x+1)=-1$).
\item Local hill $\Rightarrow$ local valley at rate $p$: When $\hat{\eta}(t,x)=1$ and $\hat{\eta}(t,x+1)=-1$, then $h_{\gamma}(t,x)$ will decrease by 2 at rate $pdt$ (and hence leaving $\hat{\eta}(t+dt,x)=-1$ and $\hat{\eta}(t+dt,x+1)=1$).
\end{itemize}
This is illustrated in figure \ref{cornergrowth}. It is possibly to encode these dynamics as a stochastic PDE driven by the two Poisson point processes $L$ and $R$ (on $\Z\times \R^+$ of rates $q$ and $p$ respectively). Recall that $L$ corresponds to the running total of attempted jumps to the left (i.e., increases in the height by $2$). So all we need to do is multiply $\partial_t L$ by an indicator function for the event that there is no particle in position $x$ and there is a particle in position $x+1$. We do likewise for $R$, which yields:

\begin{eqnarray}\label{hspde}
\partial_t h_\gamma(t,x) &=& 2\left(\frac{1-\hat{\eta}(t,x)}{2}\right)\left(\frac{1+\hat{\eta}(t,x+1)}{2}\right)\partial_t L(t,x) \\&&- 2\left(\frac{1+\hat{\eta}(t,x)}{2}\right)\left(\frac{1-\hat{\eta}(t,x+1)}{2}\right)\partial_t R(t,x).
\end{eqnarray}

Of course, the above differential relation should be interpreted in weak integrated form. This is well-posed and has solutions due to our explicit construction of such a solution.

\subsubsection{Hydrodynamics}\label{hydrosec}
Before delving into the main focus of fluctuations, it is worth observing that from the SPDE (\ref{hspde}) it is possible to read off the PDE which governs the evolution of the limit shape for the corner growth model. Assume presently that $\gamma=q-p>0$ and define, according to {\it Eulerian scaling}
\begin{equation}
h_{\gamma}^{\e}(T,X) = \e h(\e^{-1}T,\e^{-1}X), \qquad \hat{\eta}_{\gamma}^{\e}(T,X) = \e \hat{\eta}_{\gamma}(\e^{-1}T,\e^{-1}X).
\end{equation}
Define $L^\e$ and $R^\e$ to be similarly rescaled versions of $L$ and $R$.

We will show, heuristically, that
$$\lim_{\e\rightarrow 0} h_{\gamma}^{\e}(T,X) = \bar{h}_{\gamma}(T,X)$$ exists and solves a PDE. We can rewrite the SPDE (\ref{hspde}) in terms of the scaled variables, as
\begin{eqnarray}
\partial_T h_{\gamma}^{\e}(T,X) &=& 2\left(\frac{1-\hat{\eta}_{\gamma}^{\e}(T,X)}{2}\right)\left(\frac{1+\hat{\eta}_{\gamma}^{\e}(T,X+\e)}{2}\right) \partial_T L^{\e}(T,X) \\&&- 2\left(\frac{1+\hat{\eta}_{\gamma}^{\e}(T,X)}{2}\right)\left(\frac{1-\hat{\eta}_{\gamma}^{\e}(T,X+\e)}{2}\right)\partial_T R^{\e}(T,X).
\end{eqnarray}
That equation can be written suggestively as
\begin{eqnarray}
\partial_T h_{\gamma}^{\e}(T,X) &=& \frac{1}{2}\left[1-\hat{\eta}_{\gamma}^{\e}(T,X)\hat{\eta}_{\gamma}^{\e}(T,X+\e)\right]\left(\partial_T R^{\e}(T,X)-\partial_T L^{\e}(T,X)\right) \\&&+\frac{1}{2}\left[\hat{\eta}_{\gamma}^{\e}(T,X+\e)-\hat{\eta}_{\gamma}^{\e}(T,X)\right]\left(\partial_T R^{\e}(T,X)+\partial_T L^{\e}(T,X)\right). \end{eqnarray}

Due to the scaling, a law of large numbers applies and shows that $\partial_T R^{\e}(T,X)-\partial_T L^{\e}(T,X)\rightarrow \gamma dT$ and that $\partial_T R^{\e}(T,X)+\partial_T L^{\e}(T,X)\rightarrow dT$. Recall now that $\hat{\eta}$ is the derivative of $h$. The term $\left[\hat{\eta}_{\gamma}^{\e}(T,X+\e)-\hat{\eta}_{\gamma}^{\e}(T,X)\right]$ then approximates $\e \partial_X^2 \bar{h}$ as $\e$ goes to zero. The tricker term is the non-linear one $\left[1-\hat{\eta}_{\gamma}^{\e}(T,X)\hat{\eta}_{\gamma}^{\e}(T,X+\e)\right]$. A reasonable assumption (which can be shown rigorously as well) is that of {\it local equilibrium} which says that the limit can be taken inside this non-linear function of the $\hat{\eta}^{\e}$, which yields $1-(\partial_X \bar{h}(T,X))^2$. Putting this all together we find that
\begin{equation}\label{viscPDE}
\frac{\partial}{\partial_T} \bar{h} = \frac{\gamma}{2}\left(1-\left(\frac{\partial}{\partial X}\bar{h}\right)^2\right) + \e \frac{1}{2} \frac{\partial^2}{\partial_X^2}\bar{h}.
\end{equation}
Since $\gamma>0$ the last term clearly drops off and we are left with a PDE known of as the inviscid Burgers equation. There are, however, multiple ways of solving this PDE and we need to know which one is the correct one. However, the above heuristic also gives us that. The Laplacian term which disappeared in this scaling as $\e$ went to zero is known as a viscosity. For $\e$ positive, there exists a unique weak solutions to equation (\ref{viscPDE}). The solution to the inviscid Burgers equation we want is the $\e=0$ limit of the unique weak solution to the viscous Burgers equation. This is equivalent to imposing the so called {\it entropy condition} and the resulting solution can be solved via the method of characteristics (see Chapter 3 of \cite{Ev98}).

If $\gamma$ we taken to go to zero like $\e$, and time were sped up like $\e^{-2}T$ (not $\e^{-1}T$) then the above argument shows that the resulting PDE for the limit shape is, in fact, the viscous Burgers equation. Likewise, if $\gamma=0$ then only the viscosity term (the Laplacian) remains and under the same $\e^{-2}T$ time scaling the heat equation governs the evolution of the limit shape. It is important to note that the weak scaling of $\gamma=\e$ above is actually different than the weak scaling necessary to see the KPZ equation. In that setting $\gamma=\e^{1/2}$ and time is like $\e^{-3/2}T/\gamma = \e^{-2}T$ and space is like $\e^{-3/2}X$.

As explained in the introduction, the idea of the hydrodynamic PDE is that if the time zero height profile converges to a given limit shape, then the positive $T$ evolution of that initial height profile will likewise converge to the PDE evolution of the limiting profile.

\subsubsection{Fluctuations}

The hydrodynamic PDE associated with a particle system or growth model is highly dependent on the local update rules. The fluctuations around this limit shape, however, are expected to be universal and highly independent of model specifics. This belief originates from the physics work of Kardar, Parisi and Zhang \cite{KPZ} in 1986. Dealing not with this growth model, but with a few other physically inspired models, KPZ proposed that a certain continuum growth model -- the stochastic PDE which know bears their initials -- should be studied and used to determine the long time scalings of these growth models. The believe that this equation should (somehow) govern the scaling properties of the fluctuations of a wide variety of models is termed universality.

From this SPDE, KPZ predicted that in a large time $t$, fluctuations should live on the scale of $t^{1/3}$ and should display non-trivial correlation to spatial scales of $t^{2/3}$. That is to say that time : space : fluctuations should scale like 3 : 2 : 1. Thereafter any growth model which displayed such scalings was said to be in the {\it KPZ universality class}. These exponent predictions (made via non-rigorous dynamic renormalization group methods) are based on the fact that the (formal) derivative of the solution to the KPZ equation solves the stochastic Burgers equation. The work of Forster, Nelson and Stephen \cite{FNS} from 1977 provides the relevant results then to deduce these exponents. It was not until the work of \cite{BQS} that these exponents were proved rigorously for the KPZ equation with equilibrium initial data. The work of \cite{ACQ} has rigorously proved the $1/3$ exponent for the KPZ equation with narrow wedge initial data. The techniques used in these rigorous proofs were very different than the non-rigorous approach of \cite{FNS} and the approach used in \cite{ACQ} rigorously computes the associated probability distribution for the narrow wedge initial data fluctuations.

Two points should be emphasized about the KPZ paper. The first is that, despite giving a heuristic explanation for each of the terms in their equation, KPZ gave no hint as to how one would actually derive the KPZ equation from a discrete growth model. The second point is that they gave absolutely no predictions about the exact nature of the statistics associated with the fluctuations.

As people began to realize the far reaching importance of the KPZ universality class, these two questions became active subjects of research within the physics community. There have been many heuristic derivations of the KPZ equation from microscopic models \cite{Nagatani, Park}. However, the only derivation which makes mathematical sense is that of Bertini and Giacomin \cite{BG}.
The exact nature of the statistics of models in the KPZ universality class was studied by a battery of methods and too much avail (see the discussion at the end of Section \ref{KPZpred}). However, despite all of this work, exact an analytic formulas for the statistics of models in the KPZ universality class went unknown until mathematicians entered the story.

In 1999, Baik, Deift and Johansson \cite{BDJ} and Johansson \cite{KJ} gave the first rigorous proofs of the $1/3$ fluctuation exponent for growth models. Moreover, the also proved exact formulas for the statistics of the fluctuations (in the $t^{1/3}$ scale). It came as a surprise to both the physics and mathematics communities that the KPZ class statistics were identical to the statistics describing the largest eigenvalue of a random (Gaussian) Hermitian matrix. These statistics had been discovered by Tracy and Widom \cite{TW} in the early 1990s and are now often called the Tracy-Widom GUE distribution and notated as $F_{\rm{GUE}}$ or $F_2$ (we will stick with the first). The spatial correlation exponent as well as exact statistics for multi point distributions was first worked out by Pr\"{a}hofer and Spohn \cite{PS2} in 2002 for the wedge geometry. Since then the spatial correlation statistics for the other five geometries of the KPZ universality class have similarly been worked out (see Section \ref{growthregimesec} of the introduction for a chart recording all of these statistics and geometries).

Until the recent work of Tracy and Widom, however, all of the exact statistics results pertained to growth models were for totally asymmetric models with only growth (and no means for interfaces to recede). With regards to the corner growth model, this means that $\gamma=1$. Slightly earlier, \cite{BSaom} rigorously proved the $1/3$ and $2/3$ exponents for the corner growth model with any $\gamma>0$ (see also work \cite{QV,BS}). Returning to the question of exact statistics, to illustrate the type of result one has, further set  $t=\e^{-3/2}T$ and $x_k= 2^{1/3}t^{2/3}X_{k}$ for $k=1,\ldots m$ (note that these definitions are temporary as we will soon redefine $t$ and $x$ with slightly different constants for the purpose of the rest of this work). Then
\begin{equation}
\lim_{\e\rightarrow 0} \PP\left(\bigcap_{k=1}^{m} \left\{ \frac{h_{\gamma}(\tfrac{t}{\gamma},x) - \tfrac{t}{2}}{t^{1/3}} \geq 2^{-1/3}(X_k^2 -s_k)\right\}\right) = \PP\left(\bigcap_{k=1}^{m} \{\mathcal{A}_2(X_k)\leq s_k\}\right)
\end{equation}
where $\mathcal{A}_2$ represents the $\textrm{Airy}_2$ process which is stationary, continuous \cite{KJPNG}, locally absolutely continuous to Brownian motion \cite{CH}, has $F_{\rm{GUE}}$ as its one point marginal and has finite dimensional distributions which can be written in terms of Fredholm determinants \cite{PS2,KJPNG} or certain PDEs \cite{AvM,TWairy}.

The above results deal only with the spatial correlation structure of the limit of the corner growth model -- but say nothing about the joint distribution of the fluctuations for two different values of $T$. For instance, how are the height function fluctuations above the origin at time $\e^{-3/2}T_1$ and $\e^{-3/2}T_2$ asymptotically correlated. Presently, all that is known rigorously of this important question is that if $T_2 = T_1 + o(1)$ then they two fluctuations converge in probability (as $\e$ goes to zero) \cite{CFP2}.

One might guess that the limiting space-time process for the corner growth model with $\gamma=1$ is the KPZ equation. This, however, is not true (as definitively proved by \cite{ACQ}), and the true limiting process has only recently been described in \cite{CQ2} as the so-called {\it KPZ renormalization fixed point}.

The KPZ equation does, however, arise as the scaling limit of the height function fluctuations of the corner growth model under weakly asymmetric scaling of $\gamma=\e^{1/2}$. This fact was first observed and rigorously proved (for certain geometries -- not including the wedge) in 1995 by Bertini-Giacomin \cite{BG}.\footnote{The term weakly asymmetric has been used in a few other contexts with regards to the simple exclusion process. In the study of fluctuations, \cite{DeMasi,Ditt,GoncJaraCross} deal with weaker asymmetry than $\gamma=\e^{1/2}$ such as $\gamma = \e$ and find Gaussian limiting fluctuation processes which correspond to the EW universality class. In the study of hydrodynamical theory as in Section \ref{hydrosec} (see also \cite{Spo}) one finds the vicious Burgers equation by tuning the asymmetry correctly.} The rest of this section aims to give a clear explanation of their methods and results, as well as detail extensions (in \cite{ACQ,CQ}) to the other geometries (such as the wedge). Before jumping into this, we first give a heuristic explanation of \cite{BG} (see also the discussion of Section \ref{sixkpzclass}).

\begin{definition}\label{kpzscalings}
We will use the following scalings through the rest of this work:
\begin{equation}
\gamma=\e^{1/2}, \qquad t=\e^{-3/2}T, \qquad x=\e^{-1}X,\qquad \textrm{flucuations}=\e^{-1/2}.
\end{equation}
\end{definition}

The object which should converge to the KPZ equation is the height function fluctuations $h^{\rm{fluc}}_{\e}(T,X)$:
\begin{equation}\label{fluceqn}
h^{\rm{fluc}}_{\e}(T,X)=\e^{1/2}\left(h_{\gamma}\left(\tfrac{t}{\gamma},x\right)-\tfrac{t}{2}\right).
\end{equation}

Bertini-Giacomin make the following assumption on the initial data:
\begin{equation}
h^{\rm{fluc}}_{\e}(0,X)=\e^{1/2} h_{\gamma}(0,x) = \mathcal{H}_0(X) + o(1)
\end{equation}
where $\mathcal{H}_0$ is a (possibly random) function which has at most linear growth and is H\"{o}lder $<1/2$ (see Section \ref{BGeqnsec} for the exact statement of these assumptions). Then Bertini-Giacomin show that for $T>0$
\begin{equation}
h^{\rm{fluc}}_{\e}(T,X) = \mathcal{H}(T,X) + o(1)
\end{equation}
and $\mathcal{H}(T,X)$ evolves according to the KPZ equation.

Two geometries for which is applies are Brownian and flat. For the Brownian geometry, $h(0,x)$ is a simple symmetric random walk. Thus $\e^{1/2}h(0,\e^{-1}X)$ converges to a Brownian motion, i.e., $\mathcal{H}_0(X)=B(X)$ for some Brownian motion $B(X):\R\rightarrow \R$ with $B(0)=0$. Likewise for the flat geometry, $\e^{1/2}h(0,x) = 0 +o(1)$ and hence $\mathcal{H}_0(X)=0$.

There are, however, some very important geometries which are not covered by Bertini-Giacomin's results -- in particular, anything wedge-like. For instance, the initial condition for KPZ in the wedge$\rightarrow$Brownian geometry should look like $\e^{1/2}$ times $|\e^{-1}X|$ for $X<0$ and like a Brownian motion for $X\geq 0$. The problem is, however, that the limit for $X<0$ is infinity. Rather weird initial data, no doubt?

Well things get even weirder for the wedge geometry. Since $h_{\gamma}(0,x) = |x|$ it follows that $\e^{1/2}|\e^{-1}X|$ should  be the KPZ initial data. However, this is infinity for all $X\neq 0$ and $0$ for $X=0$.

Clearly one would not have much luck solving a stochastic PDE with this type of initial data. The answers to this confusion, however, comes from understanding what it means to solve the KPZ equation. After all, it turns out that as written, this equation is ill-posed due to the non-linearity. Thus the first step towards proving the results of Bertini-Giacomin, and towards extending them to include wedge-like geometries, is to rigorously define what it means to solve the KPZ equation.

\subsection{Hopf-Cole solution to the KPZ equation}
At this point we turn our attention towards making rigorous the connection between the height function fluctuations of the weakly asymmetric corner growth model and the KPZ equation. Despite its popularity and importance, there is presently no satisfactory way of making sense of the KPZ equation directly. The reason is that due to the white noise one expects spatial regularity like that of the Brownian motion (in fact for a fixed time the solution should be locally absolutely continuous with respect to Brownian motion -- like the $\textrm{Airy}_2$ process). Thus the non-linearity (the square of the first derivative) is ill-defined. One might hope that Wick ordering of the non-linearity can lead to well defined solutions. While such an approach allows one to give a definition for solving the KPZ equation, it can be shown that the scalings for those solutions are all wrong \cite{TC}. Hence, from the physical perspective we are taking, that is the wrong definition.

As we shall see below (and later as well in Section \ref{talk3}) the physically relevant way to define the solution to the KPZ equation is via the definition which we give momentarily.

\subsubsection{The Hopf-Cole solution to KPZ}
Bertini-Giacomin provide the following definition for what it means to solve the KPZ equation with initial data $\mathcal{H}_0(X)$.

\begin{definition}\label{HopfColedef}
$\mathcal{H}(T,X)$ is the {\it Hopf-Cole solution to KPZ} if $\mathcal{Z}(T,X) = \exp\{-\mathcal{H}(T,X)\}$ solves the SHE (see equation (\ref{SHE}) or Definition \ref{SHEdef} below) with $\mathcal{Z}_0(X) = \exp\{-\mathcal{H}_0(X)\}$.
\end{definition}

As long as $\mathcal{Z}_0(X)$ is a sigma-finite positive measure, then with probability one, $\mathcal{Z}(T,X)$ is positive for all $T>0$ and all $X\in \R$ \cite{Muller}. Owing to this fact, $\mathcal{H}(T,X)$ will be almost-surely a well-defined process. In fact, since it is really the SHE that we will be focusing on, initial data is more naturally thought of as initial data for SHE. This perspective opens the door for dealing with geometries such as the wedge, for which the KPZ initial data seems hopeless.

Bertini-Giacomin give some evidence for the physical relevance of this definition since they prove that it is, in fact, the scaling limit of a discrete growth model. \cite{ACQ} provides further evidence by proving that this Hopf-Cole solution has the scaling properties predicted in \cite{KPZ}.

While the KPZ equation is expected to be the weakly asymmetric scaling limit of a whole class of growth models, the approach used in \cite{BG} only works for the corner growth model. \cite{JaraGon} study a broad class of such models but can not prove that they converge to the unique Hopf-Cole solution to KPZ. Instead, they define a class of {\it energy solutions} to KPZ (which contains the Hopf-Cole solution) and show that along subsequential limits, these models converge to such solutions. However, it is far from clear if this energy condition is a strong enough condition to impose uniqueness. Uniqueness is important since one would like to use the exact solution of \cite{ACQ} to show that (at least under weak asymmetry) the statistics for this broad class of models are all the same.

\subsubsection{The Stochastic Heat equation}\label{SHEsec}

Since the KPZ equation is now defined in terms of the SHE, it is important to understand what is means to solve the SHE. The lecture notes of Walsh \cite{W} are a good reference for the study of linear stochastic PDEs and even though he does not deal with multiplicative noise, most of the theorems he states can be immediately translated into the setting we require. Walsh defines It\^{o} integration with respect to space-time white noise in a very hands on manner. For a more abstract approach of may refer to the lecture notes of Hairer \cite{Hairer} on SPDEs (though again he does not deal with multiplicative noise). For the type of initial data with which we are concerned, we actually rely on results of Bertini-Cancrini \cite{BC} which are recounted below.

Let us take the opportunity to state equation (\ref{SHE}) precisely:  $\mathscr{W}(T)$, $T \ge 0$ is the cylindrical Wiener process, i.e. the continuous Gaussian process taking values in  $H^{-1/2-}_{\rm loc} (\mathbb R)=\cap_{\alpha<-1/2} H^{\alpha}_{\rm loc}(\mathbb R)$  with
\begin{equation*}
\E\left[ \langle \varphi ,\mathscr{W}(T)\rangle \langle \psi, \mathscr{W}(S)\rangle \right] =\min(T,S) \langle \varphi, \psi\rangle
\end{equation*}
for any $\varphi,\psi\in C_c^\infty(\mathbb R)$, the smooth functions with compact support in $\mathbb R$.
Here $H^{\alpha}_{\rm loc}(\mathbb R)$, $\alpha<0$,  consists of distributions $f$ such that for any $\varphi\in C_c^\infty(\mathbb R)$, $\varphi f$ is in the standard Sobolev space $H^{-\alpha}(\mathbb R)$, i.e. the dual of $H^{\alpha}(\mathbb R)$ under the $L^2$ pairing. $H^{-\alpha}(\mathbb R)$ is the closure of $C_c^\infty(\mathbb R)$ under the norm $\int_{-\infty}^{\infty} (1+ |t|^{-2\alpha}) |\hat f(t)|^2dt$ where $\hat f$ denotes the Fourier transform. The distributional time derivative $\dot{\mathscr{W}}(T,X)$ is the space-time white noise,
\begin{equation*}
\E\left[ \dot{\mathscr{W}}(T,X) \dot{\mathscr{W}}(S,Y)\right] = \delta(T-S)\delta(Y-X).
\end{equation*}
Note the mild abuse of notation for the sake of clarity; we write $\dot{\mathscr{W}}(T,X)$ even though it is a distribution on $(T,X)\in [0,\infty)\times \mathbb R$ as opposed to  a classical function of $T$ and $X$.
Let $\mathscr{F}(T)$, $T\ge 0$, be the natural filtration, i.e. the smallest $\sigma$-field with respect to which  $\mathscr{W}(S)$ are measurable for all $0\le S\le T$.

The stochastic heat equation (\ref{SHE}) with initial data $\mathcal{Z}_0(X)$ is shorthand for its integrated (mild) version
\begin{equation}\label{mildform}
\mathcal{Z}(T,X) = \int_{-\infty}^{\infty} p(T,X-Y)\mathcal{Z}_0(dY) - \int_0^{T} \int_{-\infty}^{\infty} p(T-S,X-Y)\mathcal{Z}(S,Y)\mathscr{W}(dS,dY)
\end{equation}
where $p(T,X)$ is the standard heat kernel
\begin{equation}
p(T,X) = \frac{1}{\sqrt{2\pi T}}e^{-X^2/2T},
\end{equation}
and where the integral with respect to $\mathscr{W}(dS,dY)$ is an It\^{o} stochastic integral against white noise \cite{W}, so that, in particular,
if $f(T,X)$ is any non-anticipating integrand,
\begin{eqnarray}\label{itoexp}
&
\E\left[ (\int_0^T \int_{-\infty}^\infty f(S,Y)\mathscr{W}(dY, dS) )^2\right]=\E\left[ (\int_0^T \int_{-\infty}^\infty f^2(S,Y)dY dS\right].&
\end{eqnarray} If $\mathcal{Z}_0=\delta_{X=0}$ then the first term above simply becomes $p(T,X)$. 

One should note that the slightly awkward notation here is inherited from stochastic partial differential equations: $\mathscr{W}$ represents the cylindrical Weiner process, $\mathscr{\dot{W}}$ represents space-time white noise (the temporal distributional derivative of $\mathscr{W}$), and the stochastic integrals are with respect to space-time white noise $\mathscr{W}(dS,dY)$. For more details into the properties of white noise and the cylindrical Weiner process and It\^{o} stochastic integrals with respect to these processes, see \cite{BC,BG,W}.

Following \cite{BC} we make the following:
\begin{definition}\label{SHEdef}
The process $\mathcal{Z}(T,X)$ is a {\it mild solution} of the stochastic heat equation with initial data $\mathcal{Z}_0$ if $\mathcal{Z}(T,\cdot)$, $T>0$ is a continuous $\mathcal{F}_t$-adapted process such that for any $U>0$
\begin{equation}
\sup_{T\in (0,U]} \sup_{X\in \R} \int_{0}^T dS\int_{0}^{S} dS' \int_{\R}dy \int_{\R} dy' p(T-S,X-Y)^2 p(S-S',Y-Y')^2 \E(\mathcal{Z}(S',Y')^2) <\infty,
\end{equation}
for any $T>0$, equation (\ref{mildform}) is satisfied, and if for any uniformly bounded $f\in C^{0}(\R)$
\begin{equation}
\lim_{T\rightarrow 0^+} \int_{\R} dX f(X) \mathcal{Z}(T,X) = \int_{\R} d\mathcal{Z}_0(X)f(X).
\end{equation}
\end{definition}

In order to have existence of uniqueness of the solution to the SHE as above we must impose a level of regularity on the initial data (however not so much as to rule our working with delta functions).

\begin{definition}\label{regdef}
We call initial data $\mathcal{Z}_0(X)$ {\it regular} if it is a (possibly random) positive Borel measure on $\R$ such that for any $T>0$ it satisfies
\begin{equation}
\sup_{S\in(0,T]} \sup_{X\in \R} \E\left[\sqrt{S}\left(\int_{\R}dY p(S,X-Y)\mathcal{Z}_0(Y)\right)^2\right]<\infty.
\end{equation}
\end{definition}

Iterating the above definition leads to an alternative expression for $\mathcal{Z}(T,X)$ in terms of a convergent chaos series (see Section \ref{chaos_sec}).

The class of initial data we consider is fairly large and allows for singularities of order $T^{-1/2}$ as $T\rightarrow 0^+$ in the initial data. One can therefore check that $\mathcal{Z}_0(X)=\delta_{X=0}$ falls into the class of {\it regular} initial data. The work of Bertini-Giacomin which we will draw upon in Section \ref{discSHEtoctnsSHE} uses a slightly different class of initial data, however, as we will observe, their class of initial data is strictly contained in the above class.

We can now record a few basic results the solution to the SHE which are proved in \cite{BC} Theorem 2.2 (Claim (2) is actually from earlier work of M\"{u}ller \cite{Muller}).
\begin{proposition}\label{SHEdefprop}
Fix regular initial data (Definition \ref{regdef}), then there exists a unique solution to the stochastic heat equation as in Definition \ref{SHEdef}. Furthermore, with probability one the following two events occur: (1) $\mathcal{Z}(T,X)$ is H\"{o}lder continuous with exponent $<1/2$ in space and $<1/4$ in time; (2) For all $T>0$ and $X\in \R$, $\mathcal{Z}(T,X)>0$.
\end{proposition}

The mild solution to the SHE is really a strong solution. Given the initial data and a white noise, one can construct the solution $\mathcal{Z}(T,X)$ for all $T>0$ and $X\in R$. Under more restrictive conditions, Walsh \cite{W} provides a Picard fixed point scheme for such a construction. That approach applies to a general class of equation and for the SHE more direct approaches are available. In particular in Section \ref{chaos_sec} we will provide the chaos series for the SHE and then in Section \ref{spatialsmoothing_sec} an approximation scheme involving spatial smoothing of the white noise. Both provide ways to construct the solution of the SHE in a path-wise manner.

In order to prove the convergence of the discrete version of the SHE (which we will derive in the coming section) to the continuous SHE above we will go through the technology of martingale problems. We will hold of on a discussion of this formulation of the SHE until necessary in Section \ref{discSHEtoctnsSHE}.

\subsection{G\"{a}rtner's discrete Hopf-Cole transform}\label{Gartnersec}
We want to show that the fluctuations in equation (\ref{fluceqn}) converge to the solution of the KPZ equation. In this section we will see how to prove such a result. It turns out that for near-stationary initial data (for ASEP) this convergence is essentially true (up to a finite shift in the solution to the KPZ equation) as shown in \cite{BG}. However, for initial data corresponding to wedge geometry (or step initial condition for ASEP) the scaling is equation (\ref{fluceqn}) is more severely wrong -- it is necessary to shift the height by a constant times $\log\e$ (which becomes divergently large as $\e$ goes to zero). This disparity in scalings was not predicted and only discovered in the convergence result proof of \cite{ACQ}.

Presently we only know how to define what it means to solve the KPZ equation via Definition \ref{HopfColedef} which is in terms of the SHE. Therefore, all convergence results must be lifted to the level of the SHE by way of applying the Hopf-Cole transform to both the KPZ equation and the discrete height function. For general growth models, the Hopf-Cole transform of the height function is governed by nasty non-linear SPDEs and seeing how these converge to the SHE is difficult (and presently beyond our reach). However, for the corner growth model G\"{a}rtner \cite{G} observed that the Hopf-Cole transform of the height function actually satisfies a linear SPDE -- a discrete SHE with multiplicative noise (given by a rather messy martingale though). As this is the starting point for both Bertini-Giacomin's work as well as our own, we will now review G\"{a}rtner's transform.

We will fix our scalings just as in Definition \ref{kpzscalings}. The choice of weak asymmetry $\gamma = \e^{1/2}$ will be justified near the end of these computations, so for now we leave that as a parameter in the following expressions. Define
\begin{equation}\label{ZepsDef}
Z_{\e}(T,X) = c_{\e} \exp\left\{-\lambda_\e h_{\gamma}\left(\frac{\e^{-3/2}T}{\gamma}, \e^{-1}X\right) + \nu_{\e} \frac{\e^{-3/2}T}{\gamma}\right\}
\end{equation}
where $c_{\e}$, $\lambda_{\e}$ and $\nu_{\e}$ are $\e$ dependent parameters which will be determined in time. In order to fit with the believe that KPZ class models have $3 : 2 : 1$ scalings, we expect that $\lambda_{\e}\approx \e^{1/2}$ and $\nu_{\e} \approx \tfrac{1}{2}\e$. Under these scalings we approximately recover the expression of equation (\ref{fluceqn}). The parameter $c_{\e}$ will depend on the growth regime and is necessary in order to have $Z_{\e}(0,X)$ converge to a non-trivial limit as $\e$ goes to zero. In particular, in the work of Bertini-Giacomin, $c_{\e}=1$ while for the wedge geometry, \cite{ACQ} must take $c_{\e} = \e^{-1/2}/2$ (see Sections \ref{KPZgrowthregime} or \ref{extendBGsec} for a discussion of this point).

In order to prove convergence of $Z_{\e}$ to $\mathcal{Z}$ we need to understand what SPDE $Z_{\e}$ solves. There are three factors which affect the value of $Z_{\e}(T,X)$ in an instant $dT$:
\begin{enumerate}
\item $h_{\gamma}(t,x)$ increases by 2 at rate $\tfrac{q}{4}(1-\hat{\eta}_t(x))(1+\hat{\eta}_t(x+1))$, which correspond with:
\begin{equation}
Z_{\e}(T,X)\mapsto Z_{\e}(T,X) e^{-2\lambda_\e} \qquad \textrm{at rate}\qquad r_{\e}^{-}(T,X) = \frac{\e^{-3/2}}{\gamma} \frac{q}{4} (1-\hat{\eta}_t(x))(1+\hat{\eta}_t(x+1)).
\end{equation}
\item $h_{\gamma}(t,x)$ decreases by 2 at rate $\tfrac{q}{4}(1+\hat{\eta}_t(x))(1-\hat{\eta}_t(x+1))$, which correspond with:
\begin{equation}
Z_{\e}(T,X)\mapsto Z_{\e}(T,X) e^{2\lambda_\e} \qquad \textrm{at rate}\qquad r_{\e}^{+}(T,X) = \frac{\e^{-3/2}}{\gamma} \frac{p}{4} (1+\hat{\eta}_t(x))(1-\hat{\eta}_t(x+1)).
\end{equation}
\item Continuous growth due to the $\nu_{\e} \tfrac{\e^{-3/2}T}{\gamma}$ term:
\begin{equation}
Z_{\e}(T,X)\mapsto Z_{\e}+Z_{\e}\nu_{\e} \frac{\e^{-3/2}T}{\gamma} dT.
\end{equation}
\end{enumerate}

Putting these three pieces together, and separating the martingales from the drifts we find that $Z_{\e}$ satisfies the following SPDE:
\begin{equation}
dZ_{\e} = \Omega_{\e}Z_{\e} dT + Z_{\e} dM_{\e}
\end{equation}
where the drift term is given by
\begin{equation}\label{omegaterm}
\Omega_{\e}(T,X) = \nu_{\e} \frac{\e^{-3/2}}{\gamma} + (e^{-2\lambda_{\e}}-1)r_{\e}^{-}(T,X) + (e^{2\lambda_{\e}}-1)r_{\e}^{+}(T,X).
\end{equation}
The term $dM_{\e}(\cdot,X)$ is a martingale for each $X$ and is given by
\begin{equation}\label{marteqn}
dM_{\e}(T,X) = (e^{-2\lambda_{\e}}-1)dM_{\e}^{-}(T,X) + (e^{2\lambda_{\e}}-1)dM_{\e}^{+}(T,X)
\end{equation}
where
\begin{equation}
dM_{\e}^{\pm}(T,X) = dP_{\e}^{\pm}(T,X) - r_{\e}^{\pm}(\frac{\e^{-3/2}}{\gamma}T,\e^{-1}X)dT
\end{equation}
and $P_{\e}^{\pm}(T,X)$ (for $X\in \e\Z$) are independent Poisson processes with rates $r_{\e}^{\pm}(\frac{\e^{-3/2}}{\gamma}T,\e^{-1}X)$.

At this point it is not clear how this is related to the SHE. The martingale term is a little messy but an analysis of its quadratic variation reveals that it is not so different than space-time white noise. Right now the major stumbling block towards establishing convergence to the SHE is the fact that $\Omega_{\e}Z_{\e}dT$ appears to be nothing like a Laplacian. In fact, up to this point one could have taken any growth rule / jumping rule and written a similar (all be it more involved) SPDE.

The critical observation of G\"{a}rtner \cite{G} is that we can choose values of $\lambda_{\e}, \nu_{\e}$ and a new parameter $D_{\e}$ so that our equation is linearized via
\begin{equation}\label{gartnerlinearized}
\Omega_{\e}Z_{\e} = \frac{1}{2}D_{\e}\Delta_{\e}Z_{\e}
\end{equation}
where $\Delta_{\e}$ is the $\e$ discrete Laplacian defined via
\begin{equation}
\Delta_{\e}f(X) =\e^{-2}[f(X+\e)-2f(X)+f(X-\e)].
\end{equation}

In order to show equation (\ref{gartnerlinearized}) we must calculate $\Delta_{\e}Z_{\e}(T,X)$. Recall that by the definition of our height function, the value of $h_{\gamma}(t,x+1)$ is given by $h_{\gamma}(t,x)$ plus $\eta_{t}(x+1)$. Likewise the value of $h_{\gamma}(t,x-1)$ is given by $h_{\gamma}(t,x)$ minus $\eta_{t}(x)$. This translates into the fact that
\begin{equation}
Z_{\e}(T,X+\e) = Z_{\e}(T,X) e^{-\lambda_{\e}\hat{\eta}_t(x+1)}, \qquad \textrm{and} \qquad Z_{\e}(T,X-\e) = Z_{\e}(T,X) e^{\lambda_{\e}\hat{\eta}_t(x)}.
\end{equation}

This shows that
\begin{equation}\label{deltaeps}
\frac{1}{2}D_{\e}\Delta_{\e}Z_{\e}(T,X) = \frac{1}{2}\e^{-2}D_{\e} (e^{-\lambda_{\e}\hat{\eta}_t(x+1)}-2+e^{\lambda_{\e}\hat{\eta}_t(x)}) Z_{\e}(T,X).
\end{equation}

Observe that the expression in equations (\ref{omegaterm}) and (\ref{deltaeps}) only depend on the values of $\hat{\eta}_t(x)$ and $\hat{\eta}_t(x+1)$. Therefore, in order to show the validity of the equality in equation (\ref{gartnerlinearized}) we must choose $\lambda_{\e},\nu_{\e}$ and $D_{\e}$ such that for all four possible values of the pair $\hat{\eta}_t(x)$ and $\hat{\eta}_t(x+1)$, we have equality. A priori one should not expect that this can be done since there are four rather non-linear equation and only three unknowns to play with. Two of the equations are identical, and surprising, despite the non-linearity, the remaining three equations can be solved. The equality of the first two equations can be understood as a consequence of the symmetry of $\Delta_{\e}$ and the fact that no changes can occur when $\hat{\eta}_t(x)$ and $\hat{\eta}_t(x+1)$ are the same. The four equations are recorded in the following chart.\newline
\begin{center}
\begin{tabular}{|c|c|c|c|}
  \hline
  $\rule{0cm}{.5cm} \hat{\eta}_t(x)$&$\hat{\eta}_t(x+1)$&$\frac{1}{2}D_{\e}\Delta_{\e}Z_{\e}(T,X)$&$\Omega_{\e}(T,X)$\\\hline
  $\rule{0cm}{.5cm} 1 $ & $1 $ & $\frac{1}{2}\e^{-2}D_{\e}\left(e^{-\lambda_{\e}}-2+e^{\lambda_{\e}}\right)  $& $\frac{\e^{-3/2}}{\gamma}\nu_{\e}$ \\\hline
  $\rule{0cm}{.5cm} -1$ & $-1$ & $\frac{1}{2}\e^{-2}D_{\e}\left(e^{-\lambda_{\e}}-2+e^{\lambda_{\e}}\right)  $& $\frac{\e^{-3/2}}{\gamma}\nu_{\e}$ \\\hline
  $\rule{0cm}{.5cm}1 $ & $-1$ & $\frac{1}{2}\e^{-2}D_{\e}\left(e^{\lambda_{\e}}-2+e^{\lambda_{\e}}\right)   $& $\frac{\e^{-3/2}}{\gamma}\left(\nu_{\e}+(e^{2\lambda_{\e}}-1)p\right)$ \\\hline
  $\rule{0cm}{.5cm} -1$ & $1 $ & $\frac{1}{2}\e^{-2}D_{\e}\left(e^{-\lambda_{\e}}-2+e^{-\lambda_{\e}}\right) $& $\frac{\e^{-3/2}}{\gamma}\left(\nu_{\e}+(e^{-2\lambda_{\e}}-1)q\right)$ \\
  \hline
\end{tabular}\newline
\end{center}

The three equations can be solved by setting
\begin{equation}
\lambda_{\e} = \frac{1}{2}\log(\frac{q}{p}), \qquad \nu_{\e}= p+q-2\sqrt{pq}, \qquad D_{\e} = \frac{\e^{1/2}}{\gamma} 2\sqrt{pq}.
\end{equation}

With these parameters we find that $Z_{\e}$ satisfies the following discrete space, continuous time version of the SHE:
\begin{equation}\label{DSHE}
\partial_T Z_{\e}(T,X) = \frac{1}{2}D_{\e} \Delta_{\e} Z_{\e} + Z_{\e}dM_{\e},
\end{equation}
which is really short hand for the integrated equation
\begin{equation}\label{integratedDSHE}
Z_{\e}(T,X) = \e \sum_{Y\in \e \Z}p_{\e}(T,X-Y)Z_{\e}(0,Y) + \int_{0}^{T} \e \sum_{Y\in \e \Z} D_{\e} p_{\e}(T-S,X-Y)Z_{\e}(S,Y) dM_{\e}(S,Y).
\end{equation}

Now the choice of $\gamma=\e^{1/2}$ becomes clear. In order for the above equation to converge to the SHE as in equation (\ref{SHE}) we require $D_{\e}$ to look like $1$ as $\e$ goes to zero. In fact, by choosing $q-p=\gamma=\e^{1/2}$ we find that
\begin{equation}
\lambda_{\e} = \frac{1}{2} \log\left(\frac{1+\gamma}{1-\gamma}\right), \qquad \nu_{\e} = 1- \sqrt{1-\gamma^2}, \qquad D_{\e} = \sqrt{1-\gamma^2},
\end{equation}
which after Taylor expansion yields
\begin{equation}
\lambda_{\e} = \e^{1/2} + \frac{1}{3}\e^{3/2} + O(\e^{5/2}), \qquad \nu_{\e} = \frac{1}{2}\e + \frac{1}{8}\e^{2} + O(\e^3), \qquad D_{\e} = 1-\frac{1}{2}\e + O(\e^2).
\end{equation}

These scalings correspond exactly with the choices of $\lambda_{\e}$ and $\nu_{\e}$ suggested by the KPZ scaling theory discussed right after equation (\ref{ZepsDef}).

The martingale $M_{\e}$ is explicitly defined above in equation (\ref{marteqn}) in terms of independent (for each value of $X$) Poisson processes $P_{\e}(X,T)$. As such it is easy to calculate its quadratic variation
$d\langle M_{\e}(X,T),M_{\e}(Y,T)\rangle$ to be:
\begin{equation}
{\bf 1}(X=Y) \left((e^{2\lambda_{\e}}-1)^2 d\langle  P_{\e}^{+}(X,T),P_{\e}^{+}(Y,T)\rangle +  (e^{-2\lambda_{\e}}-1)^2 d\langle P_{\e}^{-}(X,T),P_{\e}^{-}(Y,T)\rangle\right).
\end{equation}
Using the definition of the Poisson processes we find that
\begin{equation}
d\langle  P_{\e}^{+}(X,T),P_{\e}^{+}(Y,T)\rangle = r^{+}_{\e}(T,X) = \e^{-2} \frac{p}{4} (1+\hat{\eta}_t(x))(1-\hat{\eta}_t(x+1))
\end{equation}
where we recall that $x=\e^{-1}X$ and $t=\tfrac{\e^{3/2}}{\gamma}T$.
Similarly
\begin{equation}
d\langle  P_{\e}^{-}(X,T),P_{\e}^{-}(Y,T)\rangle = r^{-}_{\e}(T,X) = \e^{-2} \frac{q}{4} (1-\hat{\eta}_t(x))(1+\hat{\eta}_t(x+1)).
\end{equation}
This can be combine in a suggestive way so that
\begin{equation}
d\langle M_{\e}(X,T),M_{\e}(Y,T)\rangle = \e^{-1}  {\bf 1}(X=Y) b_{\e}(\tau_{-\e^{-1}X}\eta_t)dT.
\end{equation}
Here $\tau_x\eta(y) = \eta(y-x)$ and
\begin{equation}\label{beps}
b_{\e}(\eta) = 1 -\hat{\eta}_t(0)\hat{\eta}_t(1) + \hat{b}_{\e}(\eta)
\end{equation}
where $\hat{b}_{\eta}(\eta)$ is explicit (see \cite{ACQ}). Moreover it follows easily that (regardless of the type of initial data) for some constant $C<\infty$,
\begin{equation}
\hat{b}_{\e} \leq C\e^{1/2}
\end{equation}
and for sufficiently small $\e>0$,
\begin{equation}
b_{\e}\leq 3.
\end{equation}

The bound for $\hat{b}_{e}$ implies that the martingale with which we are dealing has quadratic variation essentially given by (up to the $\tau$ shift) $1-\hat{\eta}(0)\hat{\eta}(1)$. If one starts the underlying exclusion process at (or near) stationary initial data, then one expects that this non-linearity should behave roughly like its expectation -- zero. The analysis of this non-linearity requires care and follows from a key estimate in \cite{BG} which we restate in Section \ref{BGeqnsec}.

The bound for $b_{\e}$ is crude but turns out to be sufficient to analyze the short time behavior of the discrete SHE. We rely on this analysis in Section \ref{extendBGsec} where we show that Bertini-Giacomin's near-stationary initial data results may be extended to a much broader class of initial conditions -- in particular the wedge geometry.

\subsection{Convergence of discrete SHE to continuum SHE under weak asymmetry}\label{discSHEtoctnsSHE}
We now recount the ingredients of Bertini-Giacomin's proof of convergence of the discrete SHE (\ref{DSHE}) to the continuous SHE (\ref{SHE}). Their work applies only for near-stationary initial data as is clear from the hypothesis given below in Section \ref{BGeqnsec}. However, in Section \ref{extendBGsec} we recount the approach of \cite{ACQ} which enables one to extend far from these settings.

Recall the following function spaces and their topologies. $C(\R^+)$ (or perhaps more precisely $C(\R,R^+)$) is the space of continuous functions from $\R$ with range supported on the positive real numbers and endowed with the topology of uniform convergence on compact subsets. This topology is, furthermore, metrizable. $C([0,T];M)$ is the space of continuous trajectories subject to the uniform topology. $D([0,T];M)$ is the space of CADLAG (continuous from the left and with a limit from the right) functions from $[0,T]$ to $M$ (a metric space) endowed with the Skhorohod topology which allows for wiggles in both time and space $M$ variables. $D_{u}([0,T];M)$ is the space of CADLAG functions from $[0,T]$ to $M$ endowed with the topology of uniform convergence on compact subsets.

\subsubsection{Bertini-Giacomin's near-stationary initial data result}\label{BGeqnsec}
In this discussion of near-stationary initial data results we take $c_{\e}=1$ (see equation (\ref{ZepsDef})).

\begin{definition}\label{neareqdef}
The class of {\it near-stationary} initial data with which Bertini-Giacomin deal satisfies the following three hypothesis:
\begin{enumerate}
\item There exists $\mathcal{Z}_0(X)$, a random function in $C(\R^+)$, such that $Z_{\e}(0,\cdot)\Rightarrow \mathcal{Z}_0(\cdot)$ as $\e\to 0$ in the topology of $C(\R^+)$.
\item For all $p>0$ there exists $a=a(p)>0$ such that
\begin{equation}
\sup_{X\in \e \Z} e^{-a|X|}\E\left[(Z_{\e}(0,X))^p\right] <\infty.
\end{equation}
\item For all $p>0$ there exist positive $a=a(p)$ and $c=c(p)$ such that
\begin{equation}
\E\left[(Z_{\e}(0,X)-Z_{\e}(0,Y))^{2p}\right] \leq c e^{a(|X|+|Y|)} |X-Y|^p.
\end{equation}
\end{enumerate}
\end{definition}
The above definition is based on Definition 2.2 of \cite{BG} though that is written in terms of height function fluctuations and not $Z_\e$ (the assumptions are just translated from one setting to the other, however).

The first hypothesis is about convergence of the initial data to $\mathcal{Z}_0$; the second is an a priori bound which ensures that $\log Z_\e$ grows at most linearly; and the third is an a priori estimate which says that the initial data is effectively H\"{o}lder with any exponent $<1/2$ (just as Proposition \ref{SHEdefprop} shows the solution to the SHE is in space).

We may now state the main near-stationary convergence result. Observe that because height changes in the corner growth model are bound in size and intensity, we can consider $Z_{\e}$ as an element of $D_u([0,T];C(\R^+))$ instead of $D([0,T];C(\R^+))$.

\begin{theorem}[Theorem 2.3 of \cite{BG}]
Let $Q_{\e}$ denote the law of $Z_{\e}(\cdot,\cdot)$ on $D_{u}([0,T],C(\R^+))$. Then the family $\{Q_\e\}_{\e>0}$ is tight, concentrates on $C([0,T],C(\R^+))$, and furthermore has a unique limit point which coincides in law with the solution to the SHE with initial data $\mathcal{Z}_{0}$.
\end{theorem}

We will now briefly review the main steps in the proof of this result.

The tightness and fact that limit points concentrate on $C([0,T],C(\R^+))$ follow from bounds which show that the H\"{o}lder continuity properties of the continuous SHE (see Proposition \ref{SHEdefprop}) are also present for this discrete SHE. Specifically, under Hypotheses 1,2 and 3, \cite{BG} Lemma 4.1 shows that for all $p>0$ there exists $a=a(p)>0$ and $c=c(p)>0$ such that
\begin{equation}
\sup_{S\in [0,T]}\sup_{X\in \e \Z} e^{-a|X|} \E\left[(Z_{\e}(S,X))^p\right] <c.
\end{equation}
Then \cite{BG} Lemma 4.2 shows that Hypothesis 3 can likewise be extended uniformly over times $S\in [0,T]$ so that, using Kolmogorov's Continuity Theorem, it is possible (as done in \cite{BG} Lemma 4.5) to prove that $Z_{\e}$ is uniformly (for $S\in [0,T]$) H\"{o}lder continuous with exponent $<1/2$ in space (as measured in terms of $L^p$ of the probability measure on space-time evolutions). Likewise \cite{BG} Lemma 4.6 shows temporal H\"{o}lder continuity with exponent $<1/4$. A slight nuance of these estimates is that they must be done on temporally linearized versions of the space-time process (otherwise the discrete height jumps mess up the estimates). The difference between the original and linearized processes is negligible as shown in \cite{BG} Lemma 4.7.

The above estimates imply tightness and the claimed concentration of limit points. The proofs of the above estimates are essentially discrete versions of the arguments one uses for the continuous SHE to prove the analogous statements in Proposition \ref{SHEdefprop}. Iteration of the integrated version of the equation along with the initial data hypotheses and applications of inequalities such as Gronwall's and Burkholder-Davis-Gundy yield the necessary estimates.

The uniqueness of the limit of the laws $Q_{\e}$ and the coincidence with the law of the solution to the SHE is shown by way of the method of the martingale problem \cite{KShiga}. In the following $f_T(X)$ denotes the canonical coordinate in $C([0,\infty),C(\R^+))$ (i.e., the space-time function $f_T(X)$ is a function of $\omega\in \Omega$ corresponding to the probability measure $Q$ below) and $(g,h)=\int_{\R} g(X)h(X)dX$.
\begin{definition}
Let $Q$ be a probability measure on $C([0,\infty],C(\R^+))$ such that for all $T>0$
\begin{equation}
\sup_{S\in[0,T]}\sup_{X\in \R} e^{-a|X|} Q\left[(f(S,X))^2\right] <\infty
\end{equation}
for some $a>0$. The measure $Q$ {\it solves the martingale problem for the SHE with initial data $\mathcal{Z}_0$} if $Q[f(0,\cdot) \in A] = \PP[\mathcal{Z}_0(\cdot)\in A]$ for all Borel sets $A\in C(\R^+)$ and if for all $\varphi\in \mathcal{D}(R)$ (the space of smooth test functions)
\begin{eqnarray}
M_T(\varphi) &=& (f_T,\varphi) - (f_0,\varphi) - \tfrac{1}{2} \int_0^{T} dS (f_S,\varphi''),\\
\Lambda_T(\varphi) &=& M_T(\varphi)^2 - \int_0^{T} dS \int_{\R} dX f_S(X)^2 \varphi(X)^2
\end{eqnarray}
are $Q$-local martingales.
\end{definition}

The key fact is that the above martingale properties uniquely characterize the solution to the SHE as shown by
\begin{proposition}[Proposition 4.11 of \cite{BG}]
For every (possibly random) $\mathcal{Z}_0\in C(\R^+)$ satisfying hypothesis 1, the martingale problem for the SHE with initial data $\mathcal{Z}_0$ has a unique solution $Q$. Moreover $Q$ coincides with the law of the process $\mathcal{Z}$ which exists and is unique via Proposition \ref{SHEdefprop}.
\end{proposition}

The discrete SHE also satisfies an analogous martingale problem. The task of proving the convergence theorem, therefore, reduces to proving that the solutions $Q_{\e}$ to the martingale problem for $Z_{\e}$ converge to the unique solution to the continuous SHE. Convergence of the linear martingales to $M_T$ is fairly straight-forward. However the convergence of the quadratic variation martingales to $\Lambda_T$ requires a key estimate about the quadratic term which showed up in equation (\ref{beps}) for $b_{\e}(\eta)$. Lemma 4.8 of \cite{BG}, translated to our context, says that for any $0<\delta<T_0<\infty$ and $\rho>0$,
there are $a,C>0$ such that for all $\delta\le S< T\le T_0$ and $\e>0$,
\begin{equation}\label{gradbd}
E [| E[  (Z_\e(T, X+\e)-Z_\e(T, X))  (Z_\e(T, X)-Z_\e(T, X-\e)) | \mathscr{F}(S)]| ]\le C \e^{1/2 -\rho} |T-S|^{-1/2} e^{ a |X| },
\end{equation}
where $\mathscr{F}(S)$ is the sigma field generated by $Z_{\e}$ up to time $S$. From this estimate it is possible to deduce the convergence of the quadratic variation martingales to $\Lambda_T$ and complete the proof of the discrete to continuous SHE convergence.

To illustrate these results consider two initial data (a) Brownian, (b) flat. Then in case (a), (letting $RW$ denote a two-sided simple symmetric random walk with $RW(0)=0$, and $B$ likewise a two-sided Brownian motion)
\begin{equation}
Z_{\e}(0,X) = \exp\{-\e^{1/2} RW(\e^{-1}X)\} \Rightarrow \exp\{-B(X)\},
\end{equation}
and in case (b), (letting $mod_2$ denote the parity function)
\begin{equation}
Z_{\e}(0,X) = \exp\{-\e^{1/2} mod_2(\e^{-1}X)\} \Rightarrow 1.
\end{equation}
Note that there is an $O(\e^{1/2})$ error in the exponentials above which is inconsequential in the $\e\to 0 $ limit.

\subsubsection{Extending beyond the near-stationary case}\label{extendBGsec}
While Bertini-Giacomin's results apply directly to Brownian, flat and flat$\to$Brownian geometries (as well as many others) it is clear that the hypotheses of Definition \ref{neareqdef} given in the previous section do not apply in any wedge-like geometries. In \cite{ACQ} it was, therefore, necessary to prove a new result which shows how results like those of Bertini-Giacomin still apply out-side of the near-stationary regime. The key fact is that in the wedge geometry, the scaling of Bertini-Giacomin must be modified so that the initial data corresponds to a non-trivial limit. Plugging in the fact that $h_{\gamma}(0,x) = |x|$ and the scalings, we find that (recalling $\lambda_{\e} \approx \e^{1/2} + \frac{1}{3}\e^{3/2}$)
\begin{equation}
Z_{\e}(0,X) = c_{\e} \exp\{-\lambda_{\e} |\e^{-1}X|\} = c_{\e}\exp\{-\e^{-1/2}|X|\}.
\end{equation}
The only scaling for $c_{\e}$ which results in a non-trivial limit is $c_{\e} = \e^{-1/2}/2$ under which $Z_{\e}(0,X)$ converges (weakly in the sense of PDEs) to $\delta_{X=0}$. (This necessarily different choice for $c_{\epsilon}$ accounts for the $log{\epsilon}$ correction given below -- see also the discussion in Section \ref{chrono}). Note that there is an $O(\e^{1/2})$ error in the exponentials above which is inconsequential in the $\e\to 0 $ limit.

It is clear that even with this $c_{\e}$ scaling, hypotheses 2 and 3 of Definition \ref{neareqdef} do not apply to the wedge geometry initial data for $Z_{\e}$. Theorem \ref{ACQBGthm} above shows that \cite{ACQ} non-the-less find a way to prove convergence of the discrete SHE to continuous SHE for the wedge geometry. The approach splits into two parts. The first step involves explicit estimations which show that for any time $\delta>0$, $Z_{\e}(\delta,\cdot)$ satisfies the hypotheses of Definition \ref{neareqdef}. Essentially this is a direct result of the regularizing properties of the discrete SHE and is proved essentially through estimating the terms of the associate Chaos series.
The second step uses the above results of Bertini-Giacomin and show that by taking $\delta\to 0$, the consistent family of solutions after time $\delta$ converge to the solution to the SHE with the desired $\delta_{X=0}$ initial data. Details of this argument can be found in Section 3 of \cite{ACQ}.

This approach provides a general scheme for how to extend beyond Bertini-Giacomin's near-stationary setting. As an illustration, in the (a) wedge$\to$Brownian or (b) wedge$\to$flat regimes $c_{\e}=1$ but hypotheses 2 and 3 of Definition \ref{neareqdef} do not hold. The same approach as in \cite{ACQ} (worked out for case (a) in \cite{CQ} as well) yields that in case (a), (letting $RW$ denote a simple symmetric random walk and $B$ a one-side Brownian motion)
\begin{equation}
Z_{\e}(0,X) = \exp\{-\e^{1/2} RW(\e^{-1}X) {\bf 1}_{X\geq 0} -\e^{1/2}|\e^{-1}X|{\bf 1}_{X< 0}\} \Rightarrow \exp\{-B(X)\}{\bf 1}_{X\geq 0},
\end{equation}
whereas in case (b), (letting $mod_2$ denote the parity function)
\begin{equation}
Z_{\e}(0,X) = \exp\{-\e^{1/2} mod_2(\e^{-1}X) {\bf 1}_{X\geq 0} -\e^{1/2}|\e^{-1}X|{\bf 1}_{X< 0}\} \Rightarrow {\bf 1}_{X\geq 0}.
\end{equation}
Note that there is an $O(\e^{1/2})$ error in the exponentials above which is inconsequential in the $\e\to 0 $ limit.

\section{The exact formula for the one-point distribution of the KPZ equation with narrow wedge initial data}\label{talk2}

In section \ref{talk1} we developed the rigorous path between the weakly asymmetric simple exclusion process (WASEP) and the KPZ equation and explained how the initial conditions for the WASEP correspond to initial data for the KPZ equation (or really the SHE). For step initial condition -- which corresponds in its scaling limit to the fundamental solution to the SHE (i.e., narrow edge initial data for KPZ) -- Tracy and Widom \cite{TW1,TW2,TW3} derived an exact formula for the probability distribution of the location of a single particle in the ASEP with step initial condition (see section \ref{TWsec} for a review of this). Using this formula and the convergence results detailed already, \cite{ACQ} provided a rigorous derivation of the exact formula for the one-point probability distribution of the Hopf-Cole solution to the KPZ equation with narrow wedge initial data. This is stated herein as Theorem \ref{ACQmainthm}.


From the formula in Theorem \ref{ACQmainthm} it is relatively easy to see that $\lim_{s\rightarrow\infty}F_{T}(s)=1$. However, from the formula it is much harder to see that $\lim_{s\rightarrow -\infty}F_{T}(s)=0$ or that $F_{T}(s)$ is non-decreasing. However, we do know that almost surely, for all $T>0$ and $X\in \R$, $\mathcal{Z}(T,X)>0$, and hence $-\log \mathcal{Z}(T,X)$ is a well-defined random variable. The fact that $F_{T}(s)$ describes its statistics is due to the convergence result of Theorem \ref{ACQBGthm} which was explained in the previous section.

The formula in Theorem \ref{ACQmainthm} is not the first result of the rigorous asymptotic analysis of Tracy and Widom's formula and in fact comes after a fair amount of post-processing. Initially one finds the following version of the formula:

\begin{equation}\label{intintop}
F_T(s):=\lim_{\e\rightarrow 0}P(F_\e(T,X)+\tfrac{T}{4!}\leq s) = \int_{\mathcal{\tilde C}} e^{-\tilde\mu}\det(I-K^{\csc}_{s})_{L^2(\tilde\Gamma_{\eta})} \frac{d\tilde\mu}{\tilde\mu},
\end{equation}
where the contour $\mathcal{\tilde C}$, the contour $\tilde\Gamma_{\eta}$ and the operator $K_s^{\csc}$ are defined below in Definition \ref{thm_definitions}.

\begin{definition}\label{thm_definitions}
The contour $\mathcal{\tilde C}$ is defined as
\begin{equation*}
\mathcal{\tilde C}=\{e^{i\theta}\}_{\pi/2\leq \theta\leq 3\pi/2} \cup \{x\pm i\}_{x>0},
\end{equation*}
navigated starting at $\infty+i$ and going counter-clockwise.

The contours $\tilde\Gamma_{\eta}$, $\tilde\Gamma_{\zeta}$ are defined as
\begin{eqnarray*}
\tilde\Gamma_{\eta}&=&\{\frac{c_3}{2}+ir: r\in (-\infty,\infty)\}\\
\tilde\Gamma_{\zeta}&=&\{-\frac{c_3}{2}+ir: r\in (-\infty,\infty)\},
\end{eqnarray*}
where the constant $c_3$ is defined henceforth as
\begin{equation*}
 c_3=2^{-4/3},
\end{equation*}
and both contours are navigated from $\pm \frac{c_3}{2} + i\infty$ down.
The kernel $K_s^{\csc}$ acts on the function space $L^2(\tilde\Gamma_{\eta})$ through its kernel:
\begin{equation}\label{k_csc_definition}
 K_s^{\csc}(\tilde\eta,\tilde\eta') = \int_{\tilde\Gamma_{\zeta}} e^{-\frac{T}{3}(\tilde\zeta^3-\tilde\eta'^3)+2^{1/3}s(\tilde\zeta-\tilde\eta')}   \left(2^{1/3}\int_{-\infty}^{\infty} \frac{\tilde\mu e^{-2^{1/3}t(\tilde\zeta-\tilde\eta')}}{e^{t}-\tilde\mu}dt\right) \frac{d\tilde\zeta}{\tilde\zeta-\tilde\eta}.
\end{equation}
\end{definition}

It is very important to observe that our choice of contours for $\tilde\zeta$ and $\tilde\eta'$ ensure that $\re(-2^{1/3}(\tilde\zeta-\tilde\eta'))=1/2$. This ensures that the integral in $t$ above converges for all $\tilde\zeta$ and $\tilde\eta'$. In fact, the convergence holds as long as we keep $\re(-2^{1/3}(\tilde\zeta-\tilde\eta'))$  in a closed subset of $(0,1)$. The inner integral in (\ref{k_csc_definition}) can be evaluated and we find that following equivalent expression:
\begin{equation*}
 K_s^{\csc}(\tilde\eta,\tilde\eta') = \int_{\tilde\Gamma_{\zeta}} e^{-\frac{T}{3}(\tilde\zeta^3-\tilde\eta'^3)+2^{1/3}s(\tilde\zeta-\tilde\eta')} \frac{\pi 2^{1/3} (-\tilde\mu)^{-2^{1/3}(\tilde\zeta-\tilde\eta')}}{ \sin(\pi 2^{1/3}(\tilde\zeta-\tilde\eta'))} \frac{d\tilde\zeta}{\tilde\zeta-\tilde\eta}.
\end{equation*}
This serves as an analytic extension of the first kernel to a larger domain of $\tilde\eta$, $\tilde\eta'$ and $\tilde\zeta$. We do not, however, make use of this analytic extension, and simply record it as a matter of interest.

\subsection{Weakly asymmetric scaling limit of the Tracy-Widom formula}\label{WASEPsec}
Due to the process level convergence of WASEP to the stochastic heat equation, exact information about WASEP can be, with care, translated into information about the stochastic heat equation. Until recently, very little exact information was known about ASEP or WASEP. The work of Tracy and Widom in the past few years, however, has changed the situation significantly. In order to use their formula we must rewrite the probability that we are interested in, as the limit of probabilities for the WASEP.

Define
\begin{equation*}
 H_{\e}(T,X) =-\log(\e^{-1/2}/2) +\lambda_{\eps}h_{\gamma}(\tfrac{t}{\gamma},x)-\nu_\eps \eps^{-1/2}t.
\end{equation*}
We are interested in understanding the behavior of  $P(H_{\e}(T,X)-\tfrac{X^2}{2T} -\tfrac{T}{24}\geq - s)$  as $\e$ goes to zero. This probability can be translated into a probability for the height function, the current and finally the position of a tagged particle:

Since we are dealing with step initial conditions $h_{\gamma}$ is initially given by $h_{\gamma}(0,x)=|x|$.  It is easy to show that because of step initial conditions, the following three events are equivalent:
\begin{equation*}
\left\{h_{\gamma}(t,x)\geq 2m-x\right\}  =\{\tilde J_{\gamma}(t,x)\geq m\} = \{\mathbf{x}_{\gamma}(t,m)\leq x)
\end{equation*}
where $\mathbf{x}_{\gamma}(t,m)$ is the location at time $t$ of the particle which started at $m>0$ and where $\tilde J_{\gamma}(t,x)$ is a random variable which records the number of particles which started to the right of the origin at time 0 and ended to the left or at $x$ at time $t$. For this particular initial condition $\tilde J_{\gamma}(t,x) = J_{\gamma}(t,x) + x\vee 0$ where $J_{\gamma}(t,x)$ is the usual time integrated current which measures the signed number of particles which cross the bond $(x,x+1)$ up to  time $t$ (positive sign for jumps from $x+1$ to $x$ and negative for jumps from $x$ to $x+1$). The $\gamma$ throughout emphasizes the strength of the asymmetry.

Therefore we have the following string of manipulations:
\begin{eqnarray}\label{string_of_eqns}
 && \qquad P(H_{\e}(T,X)-\tfrac{X^2}{2T} -\tfrac{T}{24}\geq - s)= \\
\nonumber && P\left(-\log(\e^{-1/2}/2) +\lambda_{\eps}h_{\gamma}(\tfrac{t}{\gamma},x)-\nu_\eps \eps^{-1/2}t-\tfrac{X^2}{2T} -\tfrac{T}{24} \geq - s\right)=\\
 \nonumber&& P\left(h_{\gamma}(\tfrac{t}{\gamma},x) \geq \lambda_{\e}^{-1}[-s+\log(\e^{-1/2}/2)+\frac{X^2}{2T}+\nu_\eps \e^{-1/2}t+\tfrac{T}{4!}]\right)=\\
\nonumber&& P\left(h_{\gamma}(\tfrac{t}{\gamma},x) \geq \e^{-1/2}\left[-s+\log(\e^{-1/2}/2)+\frac{X^2}{2T}\right]+\frac{t}{2}\right)=\\
 \nonumber&& P(\tilde J_{\gamma}(\tfrac{t}{\gamma},x) \geq m) = P(\mathbf{x}_{\gamma}(\tfrac{t}{\gamma},m)\leq x),
\end{eqnarray}
where $m$ is defined as
\begin{equation}\label{m_eqn}
 m=\frac{1}{2}\left[\e^{-1/2}\left(-s+\log(\e^{-1/2}/2)+\frac{X^2}{2T}\right)+\frac{1}{2}t+x\right],
\end{equation}
and $[\cdot]$ above refers to the integer part.

Thus the proof of the exact formula for the solution to the KPZ equation amounts to two pieces of information. The first is

\begin{theorem}\label{epsilon_to_zero_theorem}
For all $s\in \R$, $T>0$ and $X\in\R$ we have the following convergence:
\begin{equation}
\lim_{\e\rightarrow 0}P(H_{\e}(T,X)-\tfrac{X^2}{2T} -\tfrac{T}{24}\geq - s) = F_{T}(s)
\end{equation}
where $F_{T}(s)$ is defined in equation (\ref{fandef}).
\end{theorem}

And the second was the convergence result of Theorem \ref{ACQBGthm} which as a corollary shows that
\begin{equation}
\lim_{\e\rightarrow 0} P(H_{\e}(T,X)-\tfrac{X^2}{2T} -\tfrac{T}{24}\geq - s) = P(\mathcal{H}(T,X) - \tfrac{X^2}{2T} -\tfrac{T}{24} \geq -s).
\end{equation}

\subsubsection{Tracy-Widom formula for the asymmetric simple exclusion process with step initial condition}\label{TWsec}
The starting point for asymptotics is the remarkable formula of Tracy and Widom first stated in \cite{TW3} in the form below, and developed in the three papers \cite{TW1,TW2,TW3}. We will apply it to the last line of  (\ref{string_of_eqns}) to give us an exact formula for $P(H_{\e}(T,X)-\tfrac{X^2}{2T} -\tfrac{T}{24}\geq - s)$, and then we will take asymptotics as $\e$ goes to zero.

We consider {\it only} step initial condition ASEP, where every positive integer lattice site is initial occupied by a particle (and zero and the negative sites are empty). Recall that $\mathbf{x}_{\gamma}(t,m)$ is the location at time $t$ of the particle which started at $m>0$. Consider $q>p$ such that $q+p=1$ and let $\gamma=q-p$ and $\tau=p/q$. For $m>0$, $t\geq 0$ and $x\in \Z$, it is shown in \cite{TW3} that,
\begin{equation}\label{TW_prob_equation}
P(\mathbf{x}(\gamma^{-1}t,m)\leq x) = \int_{S_{\tau^+}}\frac{d\mu}{\mu} \prod_{k=0}^{\infty} (1-\mu\tau^k)\det(I+\mu J_{t,m,x,\mu})_{L^2(\Gamma_{\eta})}
\end{equation}
where $S_{\tau^+}$ is a circle centered at zero of radius strictly between $\tau$ and 1, and where the kernel of the operator in the Fredholm determinant is given by
\begin{equation}\label{J_eqn_def}
J_{t,m,x,\mu}(\eta,\eta')=\int_{\Gamma_{\zeta}} \exp\{\Psi_{t,m,x}(\zeta)-\Psi_{t,m,x}(\eta')\}\frac{f(\mu,\zeta/\eta')}{\eta'(\zeta-\eta)}d\zeta
\end{equation}
where $\eta$ and $\eta'$ are on $\Gamma_{\eta}$, a circle centered at zero of radius $\rho$ strictly between $\tau$
and $1$, and the $\zeta$ integral is on $\Gamma_{\zeta}$, a circle centered at zero of radius strictly between $1$ and $\rho\tau^{-1}$ (so as to ensure that $|\zeta/\eta|\in (1,\tau^{-1})$), and where, for fixed $\xi$,
\begin{eqnarray*}
\nonumber f(\mu,z)&=&\sum_{k=-\infty}^{\infty} \frac{\tau^k}{1-\tau^k\mu}z^k,\\
\Psi_{t,m,x}(\zeta) &=& \Lambda_{t,m,x}(\zeta)-\Lambda_{t,m,x}(\xi),\\
 \nonumber \Lambda_{t,m,x}(\zeta) &=& -x\log(1-\zeta) + \frac{t\zeta}{1-\zeta}+m\log\zeta
 .
\end{eqnarray*}
All contours are counter-clockwise. Throughout the rest of the paper we will only include the subscripts on $J$, $\Psi$ and $\Lambda$ when we want to emphasize their dependence on a given variable.

Let us briefly remark on the approach Tracy and Widom used in deriving this formula in \cite{TW1,TW2,TW3}. It involves solving for an exact formula for the transition probability of a finite number of ASEP particles. Marginal distribution of a fixed particle are then compactly calculated using certain, so-called magical combinatorial formulas. The finite number of particles is taken to infinity and in the step initial condition, the marginals take particularly nice forms which can be turned into Fredholm determinants and manipulated into the form above. While this approach has been pushed far, it suffers from its ad hoc nature and, due to the lack of a proper algebraic framework or understanding for this derivation it is very difficult to extend Tracy and Widom's work to other initial conditions (aside from the half-Bernoulli case \cite{TW4}) as well as to calculations of multipoint distributions.

\subsubsection{Heuristic derivation of the one-point distribution}\label{heuristicsec}
We will now present a computation deriving the expressions given in Theorem \ref{epsilon_to_zero_theorem} for $F_T(s)$. After presenting the derivation, we will stress that there are a number of very important technical points necessary to overcome to make this computation into rigorous mathematics -- many of which require serious work to resolve. In \cite{ACQ} a rigorous proof of Theorem \ref{epsilon_to_zero_theorem} is given in which we deal with all of the possible pitfalls. Besides the usual issues of convergence of integrals, trace-class convergence and cutting of contours, the calculation presented here (similar also to that of Sasamoto and Spohn) is plagued by one very significant technical problem regarding the simultaneous deformation of the two contours along which we are preforming steepest descent. The fact is that one can not freely deform these contours since doing so would introduce a diverging (with $\e$) number of poles. However, a priori, it is not clear that one can actually perform steepest descent given the rather harsh constraint on the two main contours of interest.

\begin{definition}\label{quantity_definitions}
Recall the definitions for the relevant quantities in this limit:
\begin{eqnarray*}
&& p=\frac{1}{2}-\frac{1}{2}\e^{1/2},\qquad q=\frac{1}{2}+\frac{1}{2}\e^{1/2}\\
&& \gamma=\e^{1/2},\qquad \tau=\frac{1-\e^{1/2}}{1+\e^{1/2}}\\
&& x=\e^{-1}X,\qquad t=\e^{-3/2}T\\
&& m=\frac{1}{2}\left[\e^{-1/2}\left(-s+\log(\e^{-1/2}/2)+\frac{X^2}{2T}\right)+\frac{1}{2}t+x\right]\\
&& \left\{H_{\e}(T,X)-\tfrac{X^2}{2T} -\tfrac{T}{24}\geq - s\right\} = \left\{\mathbf{x}(\frac{t}{\gamma},m)\leq x\right\}.
\end{eqnarray*}
We also define the contours $\Gamma_{\eta}$ and $\Gamma_{\zeta}$ to be
\begin{equation*}
 \Gamma_{\eta}=\{z:|z|=1-\tfrac{1}{2}\e^{1/2}\} \qquad \textrm{and} \qquad \Gamma_{\zeta}=\{z:|z|=1+\tfrac{1}{2}\e^{1/2}\}
\end{equation*}

\end{definition}
The first term in the integrand of (\ref{TW_prob_equation}) is the infinite product $\prod_{k=0}^{\infty}(1-\mu \tau^k)$. Observe that $\tau\approx 1-2\e^{1/2}$ and that $S_{\tau^+}$, the contour on which $\mu$ lies, is a circle centered at zero of radius between $\tau$ and 1. The infinite product is not well behaved along most of this contour, so we will deform the contour to one along which the product is not highly oscillatory. Care must be taken, however, since the Fredholm determinant has poles at every $\mu=\tau^k$. The deformation must avoid passing through them. Observe now that
\begin{equation*}
 \prod_{k=0}^{\infty}(1-\mu \tau^k) = \exp\{\sum_{k=0}^{\infty} \log(1-\mu \tau^k)\},
\end{equation*}
and that for small $|\mu|$,
\begin{eqnarray}\nonumber
\sum_{k=0}^{\infty} \log(1-\mu (1-2\e^{1/2})^k)& \approx& \e^{-1/2} \int_0^{\infty} \log(1-\mu e^{-2 r}) dr \\ & \approx & \e^{-1/2}\mu \int_0^{\infty} e^{-2 r} dr = -\frac{\e^{-1/2}\mu}{2}.
\end{eqnarray}
With this in mind define
\begin{equation*}
\tilde\mu = \e^{-1/2}\mu,
\end{equation*}
from which we see that if the Riemann sum approximation is reasonable then the infinite product converges to $e^{-\tilde\mu/2}$. We make the $\mu \mapsto \e^{-1/2}\tilde\mu$ change of variables and find that the above approximations are reasonable if we consider a $\tilde\mu$ contour
\begin{equation*}
\mathcal{\tilde C}_{\e}=\{e^{i\theta}\}_{\pi/2\leq \theta\leq 3\pi/2} \cup \{x\pm i\}_{0<x<\e^{-1/2}-1}.
\end{equation*}
Thus the infinite product goes to $e^{-\tilde\mu/2}$.

Now we turn to the Fredholm determinant. We determine a candidate for the pointwise limit of the kernel.  That the combination of these two pointwise limits gives the actual limiting formula as $\e$ goes to zero is, of course, completely unjustified at this point. Also, the pointwise limits here disregard the existence of a number of singularities encountered during the argument.

The kernel $J(\eta,\eta')$ is given by an integral and the integrand has three main components: An exponential term
\begin{equation*}
 \exp\{\Lambda(\zeta)-\Lambda(\eta')\},
\end{equation*}
a rational function term (we include the differential with this term for scaling purposes)
\begin{equation*}
\frac{d\zeta}{\eta'(\zeta-\eta)},
\end{equation*}
and the term
\begin{equation*}
\mu f(\mu,\zeta/\eta').
\end{equation*}
We will proceed by the method of steepest descent, so in order to determine the region along the $\zeta$ and $\eta$ contours which affects the asymptotics we consider the exponential term first. The argument of the exponential is given by $\Lambda(\zeta)-\Lambda(\eta')$ where
\begin{equation*}
\Lambda(\zeta)=-x\log(1-\zeta) + \frac{t\zeta}{1-\zeta}+m\log(\zeta),
\end{equation*}
and where, for the moment, we take $m=\frac{1}{2}\left[\e^{-1/2}(-s+\frac{X^2}{2T})+\frac{1}{2}t+x\right]$. The real expression for $m$ has a $\log(\e^{-1/2}/2)$ term which we define in with the $s$ for the moment.
Recall that $x, t$ and $m$ all depend on $\e$. For small $\e$, $\Lambda(\zeta)$ has a critical point in an $\e^{1/2}$ neighborhood of -1. For purposes of having a nice ultimate answer, we choose to center in on
\begin{equation*}
\xi=-1-2\e^{1/2}\frac{X}{T}.
\end{equation*}
We can rewrite the argument of the exponential as $(\Lambda(\zeta)-\Lambda(\xi))-(\Lambda(\eta')-\Lambda(\xi))=\Psi(\zeta)-\Psi(\eta')$. The idea in \cite{TW3} for extracting asymptotics of this term is  to deform the $\zeta$ and $\eta$ contours to lie along curves such that outside the scale $\e^{1/2}$ around $\xi$, $\re\Psi(\zeta)$ is large and negative, and $\re\Psi(\eta')$ is large and positive. Hence we can ignore those parts of the contours. Then, rescaling around $\xi$ to blow up this $\e^{1/2}$ scale, we obtain the asymptotic exponential term. This final change of variables then sets the scale at which we should analyze the other two terms in the integrand for the $J$ kernel.

Returning to $\Psi(\zeta)$, we make a  Taylor expansion  around $\xi$ and find that in a neighborhood of $\xi$,
\begin{equation*}
 \Psi(\zeta) \approx -\frac{T}{48} \e^{-3/2}(\zeta-\xi)^3 + \frac{s}{2}\e^{-1/2}(\zeta-\xi).
\end{equation*}
This suggests the change of variables,
\begin{equation}\label{change_of_var_eqn}
 \tilde\zeta = 2^{-4/3}\e^{-1/2}(\zeta-\xi) \qquad \tilde\eta' = 2^{-4/3}\e^{-1/2}(\eta'-\xi),
\end{equation}
and likewise for $\tilde\eta$. After this our Taylor expansion takes the form
\begin{equation}\label{taylor_expansion_term}
 \Psi(\tilde\zeta) \approx -\frac{T}{3} \tilde\zeta^3 +2^{1/3}s\tilde\zeta.
\end{equation}
In the spirit of steepest descent analysis, we would like the $\zeta$ contour to leave $\xi$ in a direction where this Taylor expansion is decreasing rapidly. This is accomplished by leaving at an angle $\pm 2\pi/3$. Likewise, since $\Psi(\eta)$ should increase rapidly, $\eta$ should leave $\xi$ at angle $\pm\pi/3$. The $\zeta$ contour was originally centered at zero and of radius $1+\e^{1/2}/2$ and the $\eta$ contour of radius $1-\e^{1/2}/2$. In order to deform these contours without changing the value of the determinant, care must be taken since there are  poles of $f$ whenever $\zeta/\eta'=\tau^k$, $k\in \Z$. We ignore this issue for the calculation, and deal with it carefully in \cite{ACQ} by using different contours.

Let us now assume that we can deform our contours to curves along which $\Psi$ rapidly decays in $\zeta$ and increases in $\eta$, as we move along them away from $\xi$. If we apply the change of variables in (\ref{change_of_var_eqn}), the straight part of our contours become infinite at angles $\pm2\pi/3$ and $\pm \pi/3$ which we call $\tilde\Gamma_{\zeta}$ and $\tilde\Gamma_{\eta}$. Note that this is {\em not} the actual definition of these contours which one must use in the statement and proof the theorem because of the singularity problem mentioned above.

Applying this change of variables to the kernel of the Fredholm determinant changes the $L^2$ space and hence we must multiply the kernel by the Jacobian term $2^{4/3}\e^{1/2}$. We will include this term with the $\mu f(\mu,z)$ term and take the $\e\to0$  limit of that product.

As noted before, the term $2^{1/3}s\tilde\zeta$ should have been $2^{1/3}(s-\log(\e^{-1/2}/2))\tilde\zeta$ in the Taylor expansion above, giving
\begin{equation*}
 \Psi(\tilde\zeta) \approx -\frac{T}{3} \tilde\zeta^3 +2^{1/3}(s-\log(\e^{-1/2}/2))\tilde\zeta,
\end{equation*}
which would appear to blow up as $\e$ goes to zero. We now show how the extra $\log\e$ in the exponent can be absorbed into the $2^{4/3}\e^{1/2}\mu f(\mu,\zeta/\eta')$ term.
Recall
\begin{equation*}
 \mu f(\mu,z) = \sum_{k=-\infty}^{\infty} \frac{\mu \tau^k}{1-\tau^k \mu}z^k.
\end{equation*}
If we let $n_0=\lfloor \log(\e^{-1/2}) /\log(\tau)\rfloor$, then observe that for $1 < |z| < \tau^{-1}$,
\begin{equation*}
 \mu f(\mu,z) = \sum_{k=-\infty}^{\infty} \frac{ \mu \tau^{k+n_0}}{1-\tau^{k+n_0}\mu}z^{k+n_0} =z^{n_0} \tau^{n_0}\mu \sum_{k=-\infty}^{\infty} \frac{ \tau^{k}}{1-\tau^{k}\tau^{n_0}\mu}z^{k}.
\end{equation*}
By the choice of $n_0$, $\tau^{n_0}\approx \e^{-1/2}$ so
\begin{equation*}
 \mu f(\mu,z) \approx z^{n_0} \tilde\mu f(\tilde\mu,z).
\end{equation*}
The discussion on the exponential term indicates that it suffices to understand the behaviour of this function when $\zeta$ and $\eta'$ are within $\e^{1/2}$ of $\xi$. Equivalently, letting $z=\zeta/\eta'$, it suffices to understand $ \mu f(\mu,z) \approx z^{n_0} \tilde\mu f(\tilde\mu,z)$ for
\begin{equation*}
 z= \frac{\zeta}{\eta'}=\frac{\xi +  2^{4/3}\e^{1/2}\tilde\zeta}{\xi +  2^{4/3}\e^{1/2}\tilde\eta'}\approx 1-\e^{1/2}\tilde z, \qquad \tilde z=2^{4/3}(\tilde\zeta-\tilde\eta').
\end{equation*}
Let us now consider $z^{n_0}$ using the fact that $\log(\tau)\approx -2\e^{1/2}$:
\begin{equation*}
z^{n_0} \approx (1-\e^{1/2}\tilde z)^{\e^{-1/2}(\frac{1}{4}\log\e)} \approx e^{-\frac{1}{4}\tilde z \log(\e)}.
\end{equation*}
Plugging back in the value of $\tilde z$ in terms of $\tilde\zeta$ and $\tilde\eta'$ we see that this prefactor of $z^{n_0}$ exactly cancels the $\log\e$ term which accompanies $s$ in the exponential.

What remains is to determine the limit of $2^{4/3}\e^{1/2}\tilde\mu f(\tilde\mu, z)$ as $\e$ goes to zero, for $z\approx 1-\e^{1/2} \tilde z$. This can be found by interpreting the infinite sum as a Riemann sum approximation for a certain integral. Define $t=k\e^{1/2}$ and observe that
\begin{equation}\label{Riemann_limit}
 \e^{1/2}\tilde\mu f(\tilde\mu,z) = \sum_{k=-\infty}^{\infty} \frac{ \tilde\mu \tau^{t\e^{-1/2}}z^{t\e^{-1/2}}}{1-\tilde\mu \tau^{t\e^{-1/2}}}\e^{1/2} \rightarrow \int_{-\infty}^{\infty} \frac{\tilde\mu e^{-2t}e^{-\tilde z t}}{1-\tilde\mu e^{-2t}}dt.
\end{equation}
This used the fact that $\tau^{t\e^{-1/2}}\rightarrow e^{-2t}$ and that $z^{t\e^{-1/2}}\rightarrow e^{-\tilde z t}$, which hold at least pointwise in $t$. For (\ref{Riemann_limit}) to hold , we must have $\re \tilde z$ bounded inside $(0,2)$, but we disregard this difficulty for  the heuristic proof. If we change variables of $t$ to $t/2$ and multiply the top and bottom by $e^{-t}$ then we find that
\begin{equation*}
 2^{4/3}\e^{1/2}\mu f(\mu,\zeta/\eta') \rightarrow 2^{1/3} \int_{-\infty}^{\infty} \frac{\tilde\mu e^{-\tilde zt/2}}{e^{t}-\tilde\mu}dt.
\end{equation*}
As far as the final term, the rational expression, under the change of variables and zooming in on $\xi$, the factor of $1/\eta'$ goes to -1 and the $\frac{d\zeta}{\zeta-\eta'}$ goes to $\frac{d\tilde\zeta}{\tilde\zeta-\tilde\eta'}$.

Thereby we obtain from $\mu J$ the kernel $-K_{a'}^{\csc}(\tilde\eta,\tilde\eta')$ acting on $L^2(\tilde\Gamma_{\eta})$, where
\begin{equation*}
 K_{a'}^{\csc}(\tilde\eta,\tilde\eta') = \int_{\tilde\Gamma_{\zeta}} e^{-\frac{T}{3}(\tilde\zeta^3-\tilde\eta'^3)+2^{1/3}s'(\tilde\zeta-\tilde\eta')} \left(2^{1/3}\int_{-\infty}^{\infty} \frac{\tilde\mu e^{-2^{1/3}t(\tilde\zeta-\tilde\eta')}}{e^{t}-\tilde\mu}dt\right) \frac{d\tilde\zeta}{\tilde\zeta-\tilde\eta},
\end{equation*}
with $s'=s+\log2$. Recall that the $\log2$ came from the $\log(\e^{-1/2}/2)$ term.

We have the identity
\begin{equation}\label{cscid}
 \int_{-\infty}^{\infty} \frac{\tilde\mu e^{-\tilde zt/2}}{e^{t}-\tilde\mu}dt =(-\tilde\mu)^{-\tilde z/2}\pi \csc(\pi \tilde z/2),
\end{equation}
where the branch cut in $\tilde\mu$ is taken along the positive real axis, hence $(-\tilde\mu)^{-\tilde z/2} =e^{-\log(-\tilde\mu)\tilde z/2}$ where $\log$ is taken with the standard branch cut along the negative real axis. We may use the identity to rewrite the kernel as
\begin{equation*}
 K_{s'}^{\csc}(\tilde\eta,\tilde\eta') = \int_{\tilde\Gamma_{\zeta}} e^{-\frac{T}{3}(\tilde\zeta^3-\tilde\eta'^3)+2^{1/3}s'(\tilde\zeta-\tilde\eta')} \frac{\pi 2^{1/3}(-\tilde\mu)^{-2^{1/3}(\tilde\zeta-\tilde\eta')}}{ \sin(\pi 2^{1/3}(\tilde\zeta-\tilde\eta'))} \frac{d\tilde\zeta}{\tilde\zeta-\tilde\eta}.
\end{equation*}
Therefore we have shown (without mathematical justification) that
\begin{equation*}
 \lim_{\e\rightarrow 0} P(F_{\e}(T,X)+\tfrac{T}{4!}\leq s) := F_T(s) = \int_{\mathcal{\tilde C}}e^{-\tilde \mu/2}\frac{d\tilde\mu}{\tilde\mu}\det(I-K_{s'}^{\csc})_{L^2(\tilde\Gamma_{\eta})},
\end{equation*}
where $s'=s+\log 2$. To make it cleaner we replace $\tilde\mu/2$  with  $\tilde\mu$. This only affects the $\tilde\mu$ term above given now by $(-2\tilde\mu)^{-\tilde z/2}$$=$$(-\tilde\mu)^{-2^{1/3}(\tilde\zeta-\tilde\eta')} e^{-2^{1/3}\log2(\tilde\zeta-\tilde\eta')}$. This can be absorbed and cancels the $\log2$ in $s'$ and thus we obtain,
\begin{equation*}
F_T(s) = \int_{\mathcal{\tilde C}}e^{-\tilde \mu}\frac{d\tilde\mu}{\tilde\mu}\det(I-K_{s}^{\csc})_{L^2(\tilde\Gamma_{\eta})},
\end{equation*}
which, up to the definitions of the contours $\tilde\Gamma_{\eta}$ and $\tilde\Gamma_{\zeta}$, is the desired limiting formula.

\subsubsection{Important pitfall in making it math}\label{pitfallsec}
We now briefly note some of the problems and pitfalls of the preceding computation, all of which will be addressed in the real proof of \cite{ACQ}.

Firstly, the pointwise convergence of both the prefactor infinite product and the Fredholm determinant is certainly not enough to prove convergence of the $\tilde\mu$ integral. Estimates must be made to control this convergence or to show that we can cut off the tails of the $\tilde\mu$ contour at negligible cost and then show uniform convergence on the trimmed contour.

Secondly, the deformations of the $\eta$ and $\zeta$ contours to the steepest descent curves is {\it entirely} illegal, as it involves passing through many poles of the kernel (coming from the $f$ term). In the case of \cite{TW3} this problem could be dealt with rather simply by just slightly modifying the descent curves. However, in our case, since $\tau$ tends to $1$ like $\e^{1/2}$, such a patch is much harder and involves very fine estimates to show that there exists suitable contours which stay close enough together, yet along which $\Psi$ displays the necessary descent and ascent required to make the argument work. This issues also comes up in the convergence of (\ref{Riemann_limit}). In order to make sense of this we must ensure that $1 < |\zeta/\eta'| < \tau^{-1}$ or else the convergence and the resulting expression make no sense.

Finally, one must make precise tail estimates to show that the kernel convergence is in the sense of trace-class norm. The Riemann sum approximation argument can in fact be made rigorous. In \cite{ACQ} the proof proceeds, however, via analysis of singularities and residues.

\subsection{The one-point distribution of the KPZ equation with half Brownian initial data}
Based on the exact formula of \cite{TW4}, the convergence methods of \cite{BG,ACQ} and the asymptotic methods of \cite{ACQ}, it is likewise possible to rigorously derive the one-point (crossover) formula for the half Brownian initial data. We call this family {\it edge corssover} distributions as they represent the statistics near the edge of the rarefaction fan in the WASEP. Observe that $F_{T,X}^{{\rm edge}}(s)$ now depends on $X$ (in addition to $T$) which reflects the loss of spatial stationarity which was present in the case of \cite{ACQ}.

\begin{theorem}[\cite{CQ}]
For any $T>0$ and $X\in \R$, the Hopf-Cole solution to KPZ with {\it half Brownian} initial data (given by $\mathcal{H}(T,X)= -\log\mathcal{Z}(T,X)$ with initial data $\mathcal{Z}(0,X)=e^{B(X)}{\bf 1}_{X\geq 0}$) has the following probability distribution:
\begin{equation}
\PP(\mathcal{H}(T,X)-\frac{X^2}{2T} - \frac{T}{24} \geq -s) = F_{T,X}^{{\rm edge}}(s)
\end{equation}
where $F_{T,X}^{{\rm edge}}(s)$ and is given by
\begin{equation}\label{fedgedef}
F_{T,X}^{{\rm edge}}(s) = \int_{C} \frac{d\mu}{\mu} e^{-\mu} \det(I-K^{\Gamma}_{\sigma_{T,X,\mu}})_{L^2(\kappa_T^{-1} s,\infty)}
\end{equation}
where $\kappa_T = 2^{-1/3} T^{1/3}$, $C$ is a contour positively oriented and going from $+\infty+\e i$ around $\R^+$ to $+\infty-i\e$, and $K^{\Gamma}_{\sigma_{T,X,\mu}}$ is an operator given by its integral kernel
\begin{equation}
K^{\Gamma}_{\sigma_{T,X,\mu}}(x,y) = \int_{-\infty}^{\infty} \frac{\mu}{\mu-e^{-\kappa_T t}}\Ai^{\Gamma}(x+t,\kappa_T^{-1},-2^{-2/3}\kappa_T^{-1}X)\Ai_{\Gamma}(y+t,\kappa_T^{-1},-2^{-2/3}\kappa_T^{-1}X)dt.
\end{equation}
The $\Gamma$ deformed Airy functions are defined as follows
\begin{eqnarray}
\Ai^{\Gamma}(a,b,c) &=& \int_{C^{\Gamma}} \exp\{\tfrac{1}{3}z^3 - az\} \Gamma(-bz+c)dz,\\
\Ai_{\Gamma}(a,b,c) &=& \int_{C_{\Gamma}} \exp\{\tfrac{1}{3}z^3 - az\} \frac{1}{\Gamma(bz+c)}dz.
\end{eqnarray}
The contour $C^{\Gamma}$ comes from $\infty e^{\pi i /3}$ towards the origin but goes to the left of the pole at $z=c/b$ and then leaves in the direction $\infty e^{-\pi i /3}$. The contour $C_{\Gamma}$ comes from $\infty e^{\pi i /3}$ towards the origin and then leaves in the direction $\infty e^{-\pi i /3}$ (there are no issues with poles since $1/\Gamma$ has no poles).
\end{theorem}

As a corollary, \cite{CQ} prove that scaling $s$ like $T^{1/3}$ and $X$ like $T^{2/3}$, as $T$ goes to infinity, the one-point statistic converges to that of the so called BBP-transition \cite{BBP} first discovered in the study of perturbed Wishart (LUE) random matrices. \cite{CQ} is also able to give lower bounds on the growth (as $T$ goes to infinity) of the moments of the the KPZ equation one-point distribution as well as estimates for the decay of probability of the upper tail.

By developing an FKG inequality in the context of the stochastic heat equation, \cite{CQ} then extends these moment and tail bounds to the case of the KPZ equation started with a two-sided Brownian motion as the initial data. It appears that this FKG inequality is new.

\section{Directed polymers in random media}\label{talk3}

In this section we will focus on a class of models introduced first by Huse and Henley \cite{HuHe} which we will call directed polymers in a random media (DPRM). Such polymers are directed in what is often referred to as a {\it time} direction, and then are free to configure themselves in the remaining $d$ {\it spatial} dimensions. The probability of a given configuration of the polymer is then given (relative to an underlying path measure on paths $\pi(\cdot)$) as a Radon Nikodym derivative which is often written as a Boltzmann weight involving a Hamiltonian which assigns an energy to the path:
\begin{equation}
dP_{Q}^{\beta}(\pi(\cdot)) = \frac{1}{Z^{\beta}_{Q}} \exp\{\beta H_{Q}(\pi(\cdot))\} dP_0(\pi(\cdot)).
\end{equation}
In the above equation $dP_0$ represents the underlying path measure (which is independent of the Hamiltonian and its randomness). The parameter $\beta$ is known as the inverse temperature since modifying its value changes the balance between underlying path measure (entropy) and the energetic rewards presented by the disordered or random media in which the path lives. The term $H_{Q}$ represents the Hamiltonian which assigns an energy to a given path. The subscript $Q$ stands for {\it quenched} which means that this $H_{Q}(\pi(\cdot))$ is actually a random function of the disorder $\omega$ which we think of as an element of a probability space. Finally, $Z^{\beta}_{Q}$ is the quenched partition function which is defined as necessary to normalize $dP_{Q}^{\beta}$ as a probability measure:
\begin{equation}
Z^{\beta}_{Q}= \int \exp\{\beta H_{Q}(\pi(\cdot))\} dP_0(\pi(\cdot)).
\end{equation}
The measure $dP_{Q}^{\beta}$ is a quenched polymer measure since it is still random with respect to the randomness of the Hamiltonian $H_{Q}$. This is to say that $dP_{Q}^{\beta}$ is also a function of the disorder $\omega$. We denote averages with respect to the disorder $\omega$ by an overline, so that $\overline{Z^{\beta}_Q}$ represents the averaged value of the partition function. We also use $\var{\cdot}$ to denote the variance with respect to the disorder.

When $\beta=0$ the above model reduces to the original path measure $dP_0$. Let us now focus on the case when $dP_0$ is the path measure for a standard simple symmetric random walk (SSRW). This means that for $\beta=0$, $\pi(\cdot)$ will rescale diffusively to a Brownian motion (or Brownian bridge if we pin the endpoint). A general question of interest in the study of polymer is to understand the effect of a random Hamiltonian at positive $\beta$ on the behavior and energy of a $dP_{Q}^{\beta}$ typical path. This is generally recorded in terms of two scaling exponents: the transversal fluctuation exponent $\xi$ and the longitudinal fluctuation exponent $\chi$. When $dP_0$ is supported on $n$-step SSRWs, as $n$ goes to infinity, the first exponent describes the fluctuations of the endpoint of the path $\pi$: $\var{\pi(n)} \approx n^{2\xi}$. The second exponent likewise describes the fluctuations of the free energy: $\var{\beta^{-1}\log Z^{\beta}} \approx n^{\chi}$. On top of these scaling exponents it is of essential interest to understand the statistics for the properly scaled location of the endpoint and fluctuations of the free energy.

We will now focus entirely on Hamiltonians which take the form of a path integral through a space-time independent noise field. In the discrete setting of $dP_0$ as SSRW of length $n$, the noise field can be chosen as IID random variables $w_{i,x}$ and then $H_{Q}(\pi(\cdot)) = \sum_{i=0}^{n} w_{i,\pi(i)}$.

The first rigorous mathematical work on directed polymers was by Imbrie and Spencer \cite{IS} in 1988 where (by use of an elaborate expansion) they proved that in dimensions $d\ge 3$ and with small enough $\beta$, the walk is diffusive ($\xi=1/2$). Bolthausen \cite{Bolt} strengthened the result (under same same $d\ge 3$, $\beta$ small assumptions) to a central limit theorem for the endpoint of the walk. His work relied on the now fundamental observation that renormalized partition function (for $dP_0$ a SSRW of length $n$) $W_n=Z^{\beta}_{Q} / \overline{Z^{\beta}_Q}$ is a martingale.

By a zero-one law, the limit $W_\infty=\lim_{n\to \infty} W_n$ is either almost surely $0$ or almost surely positive. Since when $\beta=0$, the limit is 1, the term {\it strong disorder} has come to refer to the case of $W_{\infty}=0$ since then the disordered noise has, indeed, had a strong effect. The case $W_{\infty}>0$ is called {\it weak disorder}.

There is a critical value $\beta_c$ such that weak disorder holds for $\beta<\beta_c$ and strong for $\beta>\beta_c$.  It is known that $\beta_c=0$ for $d\in\{1,2\}$ \cite{CSY} and $0<\beta_c\le\infty$ for $d\ge 3$. In $d\ge 3$ and weak disorder the walk converges to a Brownian motion, and the limiting diffusion matrix is the same as for standard random walk \cite{come-yosh-aop-06}.

On the other hand, in strong disorder it is known (see \cite{CSY}) that there exist (random) points at which the path $\pi$ has a positive probability (under $dP_{Q}^{\beta}$) of ending. This is certainly different behavior than that of a Brownian motion.

The behavior of directed polymer when restricted to $d=1$ has drawn significant attention and the scaling exponents $\xi,\chi$ and fluctuation statistics are believed to be universal with respect to the underlying path measure and underlying random Hamiltonian. Establishing such universality has proved extremely difficult (see \cite{S} for a review of the progress so far in this direction).

The prototype for the class of $d=1$ directed polymers is the {\it Continuous Directed Random Polymer} CDRP in which the path measure is that of Brownian motion (or bridge) and the Hamiltonian is given by a path integral through space-time white noise. We shall focus mainly on the point-to-point CDRP in which the underlying path measure is that of a Brownian bridge.

We shall see below in Section \ref{approxSec} that the CDRP is not just the prototype of the $d=1$ polymer universality class, but is also a universal scaling limit in its own right. Moreover, it is related (via a Feynman Kac interpretation) to the SHE with multiplicative noise and hence to the Hopf-Cole solution to the KPZ equation. Given its universality, it is of essential physical and mathematical interest to study the exponents and statistics of the CDRP. This relies on the recently discovered exact solvability for the CDRP (as well as a few particular approximations to it) -- see Section \ref{solvablepolymers}.

\subsection{Approximating the free energy of the continuum directed random polymer}\label{approxSec}

Recall that at the level of physics the free energy of the CDRP is written as:
\begin{equation}
\mathcal{F}(T,X) = \log \E \left[ :\:\!\exp\!: \left\{-\int_{0}^{T} \dot{\mathscr{W}}(t,b(t))dt\right\}\right]
\end{equation}
where the expectation $\E$ is over Brownian bridges $b(\cdot)$ such that $b(0)=0$ and $b(T)=X$, and $:\!\exp\!:$ is known as the {\it Wick exponential}. Observe that $\mathcal{F}(T,X)$ is random with respect to the disorder (the Gaussian space-time white noise $\dot{\mathscr{W}}$). It requires some work to make sense of this object rigorously. For instance, it is not immediately clear how to interpret integrating a Brownian bridge through a space-time white noise, and nor is it clear how to exponentiate such an expression. The problem arises due to the fact that space-time white noise is not a function, but rather a random generalized function. We now introduce five different schemes through which to define and approximate the free energy of the CDRP.

\subsubsection{Chaos series and time ordering}\label{chaos_sec}
The correct way to interpret the Wick exponential is to Taylor expand the Wick exponential, time ordering the multiple It\^{o} integrals and then switching the order of integration with $\E$. Doing this results in a series of multiple stochastic integrals whose logarithm we define as $\mathcal{F}(T,X)$:
\begin{equation}
\exp\{\mathcal{F}(T,X)\} = \sum_{n=0}^{\infty} (-1)^{n} \int_{0\leq t_1<\cdots <t_n \leq T}\int_{\R^n}P_{BB}(t_1,\ldots,t_n;x_1\ldots,x_n) \dot{\mathscr{W}}(dt_1 dx_1)\cdots \dot{\mathscr{W}}(dt_n dx_n),
\end{equation}
where $P_{BB}(t_1,\ldots,t_n;x_1\ldots,x_n)$ represents the n-step transition probability of a Brownian bridge started at 0 at time 0 and ended at $X$ and time $T$ to go through positions $x_i$ at times $t_i$ for $i=1,\ldots, n$. The multiple stochastic integrals are those of It\^{o}. The series is convergent in $L^2(\dot{\mathscr{W}})$. This follows by applying equation (\ref{itoexp}) to each term in the Chaos series (since they orthogonal) and the estimate that
\begin{equation}
\int_{0\leq t_1<\cdots <t_n \leq T}\int_{\R^n}\left(P_{BB}(t_1,\ldots,t_n;x_1\ldots,x_n)\right)^2 \leq C (n!)^{-1/2}.
\end{equation}
The positivity of the right-hand side is not immediately clear but results from the next equivalent formulation, and work of \cite{Muller}.

Conceptually, the reason for the time ordering is to avoid self-interactions of the space-time white noise. By time ordering it is possible to use the independence of the field to define multiple stochastic integrals and hence to make sense of non-linear functions of the noise. This type of procedure is often referred to as self-energy renormalization in the physical literature.

In \cite{ACQ}, the Chaos series proves useful in calculating short time asymptotics of $\mathcal{F}(T,X)$, however it seems that fixed and large $T$ statistics are unattainable directly from the analysis of this series.

\subsubsection{Stochastic PDEs}\label{SHEpdesec}
Defining $\mathcal{Z}(T,X) = p(T,X)\exp\{\mathcal{F}(T,X)\}$ (where $p(T,X)$ is the standard Gaussian heat kernel and $\exp\{\mathcal{F}(T,X)\}$ is defined as above) we find that $\mathcal{Z}$ satisfies a simple recursion relation:
\begin{equation}
\mathcal{Z}(T,X) = p(T,X) - \int_0^T \int_{\R} p(T-S,X-Y) Z(S,Y) \dot{\mathscr{W}}(dS dY).
\end{equation}
This recursion is none other than the {\it mild} or {\it integrated} form of the stochastic heat equation with multiplicative noise (see equation \ref{mildform}); the initial data for $\mathcal{F}(0,X)$ is identically 1, hence the initial data for $\mathcal{Z}$ is $\delta_{X=0}$.  This is perhaps not surprising since one would expect (formally) that the path integral for $\exp\{\mathcal{F}(T,X)\}$ is a Feynman Kac formula for an imaginary time Sch\"{o}dinger equation with space-time white noise potential (the $p(T,X)$ accounts for the difference between Brownian bridges and Brownian motions).

As the Hopf-Cole solution to the KPZ equation is defined as $-\log \mathcal{Z}(T,X)$, we find that $\mathcal{F} = \log\sqrt{2\pi T} +\tfrac{X^2}{2T} - \mathcal{H}(T,X)$. In this way we see that $\mathcal{F}$ is effectively (up to shift by a deterministic function of $T,X$) equal to the Hopf-Cole solution to the KPZ equation. Differentiating $\mathcal{H}$ in space provides us (formally) with the solution to the stochastic Burgers equation \cite{BQS}.

\subsubsection{Spatial smoothing of the white noise}\label{spatialsmoothing_sec}
Smoothing space-time white noise by spatial convolution with an approximate delta function $\delta^{\kappa}$ (i.e., $\delta^{\kappa}(x) = \kappa^{-1} h(\kappa^{-1} x)$ where $h$ is an even, positive function of compact support and integral 1) results in a smooth noise field $\dot{\mathscr{W}}^{\kappa} = \delta^{\kappa}\ast \dot{\mathscr{W}}$ through which one can make sense of a path integral along the trajectory of a Brownian bridge. Let $M$ denote this, now well-defined path integral (see Section 2.1 of \cite{BC}) and define, for this smoothed noise, the Wick (or Girsanov) exponential as
\begin{equation}
:\!\exp\!:\{M\} = \exp\{M - \tfrac{1}{2}\langle M ,M\rangle\},
\end{equation}
The quadratic variation term $\langle M,M\rangle$ is necessary to keep $\E\left[:\!\exp\!:\{M\}\right]$ a martingale in $T$. Since the smoothed white noise is still temporally independent,
\begin{equation}
\langle M,M\rangle = \int_{0}^{T} \frac{1}{2} \langle \dot{\mathscr{W}}^{\kappa}(S,b(S))\rangle dS = \frac{T}{2} \left(\delta^k\ast\delta^k(0)\right)
\end{equation}
since $ \langle \dot{\mathscr{W}}^{\kappa}(S,b(S))\rangle = \left(\delta^k\ast\delta^k(0)\right)$. As the smoothing goes to zero ($\kappa\to 0$), Bertini-Cancrini \cite{BC} prove that the smoothed partition function (times $p(T,X)$) converges in $L^2(\dot{\mathscr{W}})$ to $\mathcal{Z}(T,X)$ with $\mathcal{Z}(0,X) = \delta_{X=0}$. The removal of the quadratic variation can be interpreted as another form of self-energy renormalization.

\subsubsection{Discrete directed polymers}\label{AKQDDPSEC}
There are many ways to discretize, or semi-discretize the CDRP. Perhaps the simplest approach is to replace the space-time white noise by an IID field of random variables, and to replace the Brownian bridge by a simple symmetric random walk bridge. Specifically let $w_{i,j}$ be IID random variables and define the polymer measure on $\Pi_{n,y}$ (the collection of simple symmetric random walk trajectories from 0 at time 0 to $y$ at time $n$) as
\begin{equation}
\PP^{\beta}_{n,y}(\pi) = \frac{1}{Z^{\beta}(n,y)} \exp\{ \beta T(\pi)\}
\end{equation}
where the discrete path integral $T(\pi)=\sum_{i=0}^{n} w_{i,\pi(i)}$ and the partition function is
\begin{equation}
Z^{\beta}(n,y) = \sum_{\pi\in \Pi_{n,y}} \exp(\beta \sum_{i=0}^{n} w_{i,\pi(i)}) = \sum_{\pi\in \Pi_{n,y}} \prod_{i=0}^{n} (1+\tilde{w}_{i,\pi(i)}).
\end{equation}
Here $\beta$ is called the inverse temperature, and $\tilde{w}_{i,j}$ are coupled to $w_{i,j}$ so that $e^{\beta w_{i,j}} = (1+\tilde{w}_{i,j})$. This means that for $\beta$ small, $\tilde{w}_{i,j} \approx \beta w_{i,j}$ by Taylor approximation.

At this point it is worth noting that it is possible to recover $\PP^{\beta}_{n,y}$ from partition functions. For instance if we let $Z^{\beta}(m,x;n,y)$ denote the quenched partition function for a polymer started from $x$ at time $m$ and ended at $y$ at time $n$, then
\begin{equation}
\PP^{\beta}_{n,y}(\pi(m)=x) = \frac{ Z^{\beta}(0,0;m,x) Z^{\beta}(m,x;n,y)}{Z^{\beta}(0,0;n,y)}.
\end{equation}
In fact, this also worked for the CDRP. In fact, we have not yet described what the quenched path measure is for the CDRP, though it can be built in terms of coupled solution to the SHE started and ended at different space-time locations (see for instance \cite{CQ2}).

If we denote the normalized partition function $\tilde Z^{\beta}(n,y) = 2^{-n}Z^{\beta}(n,y)$ then it clearly satisfies the recursion
\begin{equation}\label{tildezrec}
\tilde Z^{\beta}(n,y) = \tfrac{1}{2}\left(\tilde Z^{\beta}(n-1,y-1)+\tilde Z^{\beta}(n-1,y+1)\right)(1+\tilde{w}_{n,y}),
\end{equation}
which is a discrete form of the SHE with multiplicative noise.

As observed by Calabrese-Le Doussal-Rosso \cite{CDR} and Alberts-Khanin-Quastel \cite{AKQ}, under diffusive time-space scaling of $n=\e^{-4}T$, $y=\e^{-2}X$  and under {\it weak noise} scaling $\beta= \e \alpha$, the normalized partition function $\tilde Z^{\beta}(n,y)$ converges to $\mathcal{Z}(\alpha^4 T,\alpha^2 X)$ where $\mathcal{Z}$ solves the multiplicative noise stochastic heat equation with delta initial data. Alberts-Khanin-Quastel furthermore presented a scheme (which they are now working to make rigorous) through which this result can be proved for all $w_{i,j}$ which are IID (subject to centering and scaling which only depends on a certain number of finite moments) in time and space. The essential idea is that for small $\beta$, $\tilde{w}_{i,j}\approx \beta w_{i,j}$. Iterating the recursion (\ref{tildezrec}) and using this we get a finite series in $\beta$ which can be readily identified as a discrete version of the Chaos series of Section \ref{chaos_sec}. The observation is the basis for the proof. One should note that it is possible to try to make such a scheme work in the context of the discrete SHE for the corner growth model, yet the structure of the martingale is much more complex and this approach has not been successfully implemented.

This scheme can be extended to encompass more general noise and underlying random walk measures. Thus, in this weak scaling regime for the noise, the CDRP is a universal scaling object, and its statistics provide the asymptotic statistics for all models which scale to it! This is very much analogous to the weak asymmetry scaling for the corner growth model under which the KPZ equation arises.

\subsubsection{Growth processes and interacting particle systems}
The universality of the KPZ equation has already been discussed, and thus, owing to the equivalence of the KPZ equation to the free energy of the CDRP, one can consider growth models and interacting particle system height functions as approximation schemes for the free energy of the CDRP. Since these models are much further from the CDRP or SHE, this approximation scheme has only been rigorously implemented for the corner growth model (or equivalently SEP) in \cite{BG,ACQ,CQ} (or see Theorem \ref{ACQBGthm} above).

From above we know that weak noise discrete directed polymers and the weakly asymmetric corner growth model both scale to the KPZ equation. This connection is recent, though it is not the first connection observed between these models. {\it Last passage percolation} (LPP) is the zero temperature limit of the discrete directed polymer, in which the path measure concentrates entirely on the (in the case of continuously distributed weights) unique path which maximizes its path integral through the IID field (one can also look at the minimal path which is called first passage percolation). The free energy for this model is this maximal path integral -- known as the last passage time. If the path measure underlying the polymer is a nearest neighbor random walk, and the noise is formed from IID exponential random variable weights, then there is an exact coupling between LPP and the totally asymmetric corner growth model started with a wedge (see for example \cite{PS1} or \cite{BAC}).

\subsection{Beyond the point-to-point polymer}\label{beyondptpt}
The six fundamental initial data for the KPZ equation (Section \ref{KPZgrowthregime}) correspond with the CDRP with either final or initial potential \cite{CQ} and arise from appropriately tuned boundary data for discrete polymers as well (see \cite{CFP1} for example). By boundary data for a discrete polymer one means that the $w_{i,j}$ for $i=j$ or $i=-j$ have different distributions generally with a large mean. This creates an attractive potential along the boundary of the paths trajectory and hence causes it to remain straight much longer than a standard random walk. Critically tuning the type and strength of the boundary potential results in different subclasses of statistics \cite{BBP,BAC,CFP1}.

Under the Feynman Kac interpretation for the SHE with initial data $\mathcal{Z}_0$,
\begin{equation}
\mathcal{Z}(T,X) = \E_{X}\left[\mathcal{Z}_0(b(T)) :\:\!\exp\!: \left\{-\int_{0}^{T} \dot{\mathscr{W}}(t,b(t))dt\right\}\right]
\end{equation}
where $\E_{X}$ represents the expectation over Brownian motions which start at $X$ at time $0$. Thus $\log\mathcal{Z}_0$ can be interpreted as an a form of energetic cost/reward of ending at a certain location. Alternatively, due to time reversal properties of Brownian motion, it is possible to consider $\mathcal{Z}_0$ as the exponential of the integrated potential in which a CDRP chooses to depart from the time zero spatial axis.

We now review the six fundamental initial data for the KPZ equation (Section \ref{KPZgrowthregime}) and provide (1) a brief explanation for their polymer interpretation, (2) how they arise from weak noise scaling of discrete polymers with given terminating conditions, and (3) how they likewise arise from polymers with boundary conditions. Note that this dictionary between discrete polymers and initial data for the KPZ equation relies on the approach of \cite{AKQ} as well as the slow decorrelation results of \cite{CFP2}.

Below $RW$ represents a random walk independent of all else such that $RW(0)=n$.

\begin{itemize}
\item $\mathcal{Z}(0,X)=\delta_{X=0}$: (1) Point to point polymer (discussed at length above); (2) Limit of point to point discrete polymers; (3) Limit when boundary conditions are subcritical and not strong enough to attract the polymer path.
\item $\mathcal{Z}(0,X) = e^{B(X)}$: (1) Point to Brownian polymer; (2) Limit of polymers which go from a fixed point and terminate at the first space-time point $(\pi(m),m)$ at which $RW(\pi(m))=m$; (3) Limit when boundary conditions are random and critically tuned (in terms of mean). There is, in fact, a scaling window for the critical tuning and within that window one introduces a drift into the two sides of the Brownian motion $B(X)$. The long time scaling behavior of this whole window is understood (in terms of last passage percolation) in \cite{BP08,BFP,CFP1}.
\item $\mathcal{Z}(0,X) = 1$: (1) Point to line polymer; (2) Limit of polymers which go from a fixed point and terminate at a fixed time $n$ (without a fixed endpoint); (3) Limit when boundary conditions are deterministic and critically tuned.
\item $e^{B(X)}{\bf 1}_{X\geq 0}$: (1) Point to half Brownian polymer; (2) Limit of polymers which go from a fixed point and terminate at the first space-time point $(\pi(m),m)$ with $\pi(m)\geq 0$ at which $RW(\pi(m))=m$ (if no such point exists, then the path is not considered); (3) Limit when boundary conditions are random and the weights along the line $i=j$ are critically tuned (in terms of mean) while the weights along the line $i=-j$ are subcritical.
\item ${\bf 1}_{X\geq 0}$: (1) Point to half line polymer; (2) Limit of polymers which go from a fixed point and terminate at a fixed time $n$ subject to the condition that $\pi(n)\geq 0$; (3) Limit when boundary condition weights along line $i=j$ are deterministic and critically tuned, while along $i=-j$ the weights are subcritical.
\item  ${\bf 1}_{X<0} + e^{B(X)}{\bf 1}_{X\geq 0}$: Point to half line / Brownian polymer; (2) Limit of polymers which go from a fixed point and terminate either at the first space-time point $(\pi(m),m)$ with $\pi(m)\geq 0$ at which $RW(\pi(m))=m$ or terminates at time $n$ if $\pi(n)<0$; (3) Limit when boundary condition weights along line $i=j$ are random and critically tuned, while along $i=-j$ the weights are deterministic and critically tuned.
\end{itemize}

As recorded in the last column of Figure \ref{sixfig}, we are far from having a complete characterization of the statistics for the six types of initial data / polymers above. A goal, therefore, of developing the solvability of the CDRP is to be able to calculate one-point and multi-point statistics for each of the six different types of initial data. One expects to be able to do this owing to the fact that the long time limits (in terms of last passage percolation) are known and solvable in each of these cases.

We also remark that, as noted in the open question of Section \ref{opensec}, there exists a multi-layer extension to the polymer free energy whose statistics correspond to some sort of finite temperature version of the top eigenvalues of large LUE random matrices. See \cite{KJ, BP08, DW08} for explanations of how the zero-temperature polymer -- last passage percolation is related to LUE eigenvalues. Also, a similar relationship holds in the zero-temperature version of the Brownian polymer \cite{OCon} with LUE replaced by GUE.

For completeness let us briefly review the construction of the multi-layer extension as done in \cite{OConWarren}.
Define
\begin{equation}
\mathcal{Z}_n(T,X) = P(T,X)^n \sum_{k=0}^{\infty} (-1)^k \int_{0\leq t_1<\cdots <t_k \leq T}\int_{\R^k} R^{(n)}_k((t_1,x_1),\ldots,(t_k,x_k)) \dot{\mathscr{W}}(dt_1 dx_1)\cdots \dot{\mathscr{W}}(dt_k dx_k).
\end{equation}
Here $R^{(n)}_k$ is the $k$-point correlation function for a collection of $n$ non-intersecting Brownian bridges which all start at $0$ at time $0$ and end at $X$ at time $T$. These $\mathcal{Z}_n(T,X)$ form the partition function hierarchy. The multi-layer extension to the polymer free energy is achieved by setting
\begin{equation}
\mathcal{F}_n(T,X) = \log \frac{\mathcal{Z}_n(T,X)}{\mathcal{Z}_{n-1}(T,X)} - \log P(T,X),
\end{equation}
where by convention we set $\mathcal{Z}_{0}=1$.
Clearly $\mathcal{F}_1(T,X) = \mathcal{F}(T,X)$, as defined earlier. For a fixed $T$, these free energies form an ensemble of lines which may intersect (though have a tendency which grows with $T$ to avoid such overlapping). It is believed that they represent crossover, or finite temperature versions of the Airy line ensemble first introduced in \cite{PS2} and studied extensively since (see \cite{CH} in particular where the existence and various properties are proved). Thus, another goal of solvability is to calculate various statistics related to these additional levels of the multi-layer free energy process.

\subsection{Solvability of polymer models}\label{solvablepolymers}
The first suggestion that finite temperature polymer models may be solvable came from the solvability of the zero temperature analogous -- last passage percolation \cite{BDJ,KJ} with exponential or geometric weights (or a Poisson point environment). Likewise this corresponds with the solvability of the TASEP. The work of Tracy and Widom \cite{TW1,TW2,TW3} did not deal directly with polymers, however through the work of \cite{ACQ} their exact formulas for ASEP have been translated into exact one-point statistics for the free energy of the CDRP. In this section we discuss the solvability approaches which have presently been developed to date.

\subsubsection{Tracy Widom formulas}
Tracy and Widom \cite{TW1,TW2,TW3} extended earlier formulas of Sch\"{u}tz \cite{Sch} to give exact formulas for the transition probabilities for $N$ particles in the ASEP . Then, using two ``magical'' combinatorial formulas (whose origin and meaning are still a mystery) as well as a great deal of work manipulating formulas and Fredholm determinants, Tracy and Widom produced exact formulas for the transition probability of a single particle in the ASEP with step (and step Bernoulli) -- see Section \ref{sixkpzclass} for an explanation of these initial conditions. \cite{ACQ} (and \cite{CQ}) then took rigorous asymptotics of these formulas to derive the probability distribution for the free energy of the CDRP in the point-to-point (and half-Brownian) case -- or equivalently the KPZ equation with narrow-wedge or wedge$\rightarrow$Brownian initial data.

This approach was the first taken and has a few advantages. The pre-asymptotic formula (given by Tracy and Widom \cite{TW1,TW2,TW3} and restated in Theorem \ref{TW}) is already written as a Fredholm determinant and has undergone significant manipulations to make it amenable to asymptotics. On the other hand, taking rigorous asymptotics of this formula in the scaling necessary to correspond to the CDRP is a lengthy ordeal involving a number of highly non-trivial technical issues. Moreover, the derivation of the pre-asymptotic formulas is essentially ad-hoc and presently limited (to the extent needed for asymptotics) to the two cases mentioned above. Moreover, only one-point distributions have been computed by this method, rendering it presently impossible to find multi-point distributions in this manner for the KPZ equation.

\subsubsection{Replica trick and attractive delta-Bose gas Bethe ansatz}\label{replicasec}

The replica approach \cite{Dot,CDR} uses the polymer formulation of $\mathcal{Z}$ to express the moments of $\mathcal{Z}(T,X)$ (with respect to the disorder induced by the white noise potential) in terms of the solution to a quantum many body system governed by a the Lieb-Liniger Hamiltonian with two-body {\it attractive} delta interaction. Specifically if we write joint moments at time $T$ as
\begin{equation}
Z(X_1,\ldots,X_n;T) = \overline{\prod_{i=1}^{n}\mathcal{Z}(T,X_i)}
\end{equation}
then $Z(X_1,\ldots,X_n;T)$ satisfies the equation
\begin{equation}
-\partial_T  Z(X_1,\ldots,X_n;T) = H_n Z(X_1,\ldots,X_n;T),\qquad H_n = -\frac{1}{2}\sum_{j=1}^{n} \partial_{x_j}^2 - \frac{1}{2}\sum_{j\neq k}^{n} \delta(x_j-x_k),
\end{equation}
solved on the {\it bosonic subspace} (which means functions invariant under permutations of its entries).

Mathematically, it is convenient to restrict attention to the Weyl chamber $X_1<X_2<\cdots<X_n$. Then solving the above equation means that $Z$ solves the free equation (without $\delta$ interaction) inside the Weyl chamber, and at the boundary it satisfies the following boundary conditions
\begin{equation}
(\partial_{X_{i+1}}-\partial_{X_i} +1)\Psi(x)=0
\end{equation}
when $X_{i+1}=X_{i}$.

One may to solve this many body system is by demonstrating a complete basis of eigenfunctions (and normalizations) which diagonalize the Hamiltonian and respect the boundary condition.
The eigenfunctions we written down in 1963 for the {\it repulsive} delta interaction (the $+1$ in the boundary condition becomes $-1$) by Lieb and Liniger \cite{LL} by Bethe ansatz and their completeness was proved by Dorlas \cite{Dorlas} (on $[0,1]$) and Tracy and Widom \cite{TWBose} (on $\R$ as we are considering presently). For the attractive case, McGuire \cite{McGuire} wrote the eigenfunctions in terms of {\it string states} in 1964. However, the norms of these states were not derived until 2007 in \cite{CalCaux} using ideas of algebraic Bethe ansatz (see \cite{KirKor,BogIzer,Slavnov}). Dotsenko \cite{Dot} later worked these norms out very explicitly through combinatorial means. According to \cite{vanDiej} completeness of the eigenfunctions for the attractive case on $\R$ was shown in the Ph.D. thesis of Oxford \cite{Oxford}. The author, however, has not been able to access this Ph.D. thesis to determine to what degree this work is mathematically rigorous. Moreover, this work apparently is unknown to several experts. In \cite{ProSpoComp} the formula of \cite{TWBose} is analytically continued to consider an attractive delta potential and in this way the completeness for this case is proved.

As discussed in Section \ref{chrono} this replica approach is plagued by the fact that moments of $\mathcal{Z}$ grow far to rapidly to uniquely identity the Laplace transform of $\mathcal{Z}$ or the distribution of $\log \mathcal{Z}$. It took some time for \cite{Dot} and \cite{CDR} to work how to re-sum divergent series and analytically continue functions (only a priori defined on the integers) as necessary to derive distribution formulas. It was not, in fact, until after the work of \cite{SaSp1,SaSp2,SaSp3,ACQ} that the replica approach was able to recover the newly discovered distribution functions.

One should note that even though recovering the Laplace transform via the replica trick is non-rigorous mathematically, the formulas for moments derived via the replica trick are essentially rigorous (up to the question of completeness) and can also be rigorously calculated in terms of local-time as in \cite{BC}. Furthermore, even with regards to the full distribution function, the correctly applied replica trick has some benefits. Prolhac and Spohn \cite{ProS1,ProS2} used this approach to derive a conjectural form of the spatial process for the KPZ equation in the geometry corresponding to growth in a narrow wedge and confirm that the long time limit of this spatial process is the ${\rm Airy}_2$ process.\footnote{This and the work of \cite{CQ2} make a critical factorization assumption in the form of the moments of $\mathcal{Z}$ which may not be true at finite $t$ but which appears to hold in the long-time limit.} In \cite{CQ2} this approach is used to derive a conjectural form for the transition probabilities for the random non-linear semi-group which governs the KPZ renormalization fixed point.

Very recently \cite{CD} employed the replica trick to derive a conjectural formulas for the crossover one-point distributions for the flat $\mathcal{Z}(0,X)=1$ case. Also, \cite{ImSahalfbrown} used this approach to rederive the formula for the statistics correspond to $\mathcal{Z}(0,X)=\e^{B(X)}{\bf 1}_{X\geq 0}$, which \cite{CQ} had rigorously derived a year earlier based on work of \cite{TW4}.

\subsubsection{Solvable finite temperature polymers}
The solvability of LPP \cite{KJ,BDJ} relies on the combinatorial Robinson-Schensted-Knuth (RSK) correspondence. RSK maps a matrix of positive weights onto a pair of Young Tableaux -- or equivalently {\it Gelfand Zetlin (GZ) patterns} -- from which one can immediately read off information like the last passage time for the original matrix. When matrix weights are chosen as independent exponential random variables, the resulting measure on GZ patterns is given by the Schur measure and in this case it is possible to write exact formulas for the probability distribution for the last passage time  (and many other quantities) -- hence the solvability.

Last passage percolation represents a zero temperature polymer model, and thus in order to access the statistics of the KPZ equation, it is necessary to find solvable finite temperature polymers. This is accomplished in \cite{OCon} for semi-discrete {\it Brownian polymer} and \cite{COSZ} for a particular discrete polymer involving inverse-gamma weights (which was initially studied in \cite{S} in which a Burkes type theorem was uncovered -- first hinting at solvability). In the case of \cite{COSZ}, the solvability relies on a few facts. The first is the existence of a finite temperature version of RSK -- the {\it tropical RSK correspondence} introduced by A.N. Kirillov \cite{Kir} in the context of tropical combinatorics. The classical RSK correspondence can be encoded as a combinatorial algorithm over the $(\max,+)$ algebra -- gRSK amounts to formally replacing: $\max\mapsto +$ and $+\mapsto \times$. The image of a matrix of positive entries under gRSK is a triangular array from which one may immediately read off the polymer partition function associated with the original matrix. This can be done with any weight matrix, however, in the case of inverse-gamma distributed IID weights, a miracle occurs and, due to a certain {\it intertwining relation} and the theory of Markov functions \cite{PR}, it is possible to show that a wide variety of projections of the triangular array have the structure of Markov chains (with respect to the filtration formed by their own history). The transition kernels for these chains can be diagonalized in terms of ratios of Whittaker functions (which plays the tropical analogue of Schur functions) and an integrate-out lemma (the Bump-Stade identity) originally developed in the study of automorphic forms enables the exact calculation of the Laplace transforms of of polymer's partition function.

The algebraic structure associated with the finite temperature polymer solvability is much better understood than in the context of Tracy and Widom's work with the ASEP. Therefore, much more information about the KPZ equation should be accessible through this approach. Since the expression naturally (and rigorously) derived by this approach is the Laplace transforms of the partition function (which is a positive random variable), we can invert and recover expressions the distribution function of the free energy.


\end{document}